\documentclass[11pt]{article}
\usepackage{comment} % 引入comment宏包
\usepackage{geometry}
\usepackage{amsmath}
\usepackage{graphicx}
\usepackage{amsmath,amssymb,amsthm}
\usepackage[utf8]{inputenc}
\usepackage[greek,english]{babel}
\usepackage{enumerate}
\usepackage{mathrsfs}
\usepackage{hyperref}
\usepackage[active]{srcltx}
\numberwithin{equation}{section}
\usepackage[T1]{fontenc}
\usepackage{color}
\newtheorem{theo}{Theorem}[section]
\newtheorem{lem}{Lemma}[section]
\newtheorem{defi}{Definition}[section]

\newtheorem{prop}{Proposition}[section]
\newtheorem{rmk}{Remark}[section]

\newcommand{\R}{\mathbb{R}}

\newcommand{\N}{\mathbb{N}}

\newcommand{\me}{\mathrm{e}}

\newcommand{\bqs}{\begin{equation*}}
	\newcommand{\eqs}{\end{equation*}}
\newcommand{\bqq}{\begin{equation}}
	\newcommand{\eqq}{\end{equation}}

\DeclareMathAlphabet\mathbfcal{OMS}{cmsy}{b}{n}

\date{}

\begin{document}
\title{Effects of temporal variations on wave speeds of bistable traveling waves for Lotka-Volterra competition systems}	
\author{Weiwei Ding\footnote{School of Mathematical Sciences, South China Normal University, Guangzhou 510631, China} \qquad Zhanghua Liang\footnote{School of Mathematical Sciences, South China Normal University, Guangzhou 510631, China}  }
\maketitle

\begin{abstract}
This paper investigates the bistable traveling waves for two-species Lotka-Volterra competition systems in time periodic environments. We focus especially on the influence of the temporal period, with existence results established for both small and large periods.
We also show the existence of, and derive explicit formulas for, the limiting speeds as the period tends to zero or infinity, and provide estimates for the corresponding rates of convergence. Furthermore, we analyze the sign of wave speed.  
Assuming that both species share identical diffusion rates and intraspecific competition rates, we obtain a criterion for determining the sign of wave speed by comparing the intrinsic growth rates and interspecific competition strengths.
More intriguingly, based on our explicit formulas for the limiting speeds, we construct an example in which the sign of wave speed changes with the temporal period. This example reveals that temporal variations can significantly influence competition outcomes, 
enabling different species to become dominant under different periods.


\end{abstract}

%%%%%%%%%%%%%%%%%%%%%%%%%%%%%%%%%%%%%%%%%%%%%%%%%%%%%%%%%%%%%%%%%%%%%%%%%%%%%%%%%%%%%%%%%%

\section{Introduction and main results}
\subsection{Introduction}
In this paper, we study the propagation phenomena the following Lotka-Volterra competition-diffusion model in a time periodic environment:
\begin{equation}\label{eq:main}
\left\{\begin{array}{ll}
\partial_t u_1=d_{1, T}(t) \partial_{xx}u_1+u_1\left(r_{1, T}(t)-a_{1, T}(t) u_1 -k_{1, T}(t)u_2\right), & \hbox{for }\,\, (t,x) \in \R \times\R,  \vspace{5pt}\\
\partial_t u_2=d_{2, T}(t) \partial_{xx}u_{2}+u_2\left(r_{2, T}(t)-a_{2, T}(t) u_2 -k_{2, T}(t) u_1\right), &\hbox{for }\,\, (t,x)\in \R\times\R,
\end{array}\right.
\end{equation}
with the time period $T>0$. For each $i=1,2$, the functions $d_{i, T}, r_{i, T},a_{i, T}, k_{i, T}$ are given by
$$
d_{i, T}(t)=d_i(t / T), \quad r_{i, T}(t)=r_i(t / T), \quad a_{i, T}(t)=a_i(t / T),\quad k_{i, T}(t)=k_i(t / T),
$$
where $d_i,r_i,a_i,k_i \in C^{2}(\mathbb{R})$ are positive 1-periodic functions. 

Model \eqref{eq:main} describes the evolution of interaction and diffusion between two competing species, where $u_1(t, x)$ and $u_2(t,x)$ represent the population densities of each species. For $i=1,2$, $d_{i, T}$ denotes of the diffusion rate of species $u_i$;  $r_{i, T}$ and $a_{i,T}$  stand for its intrinsic growth rate and intraspecific competition coefficient, respectively;  
and $k_{i, T}$  measures the strength of interspecific competition. The periodic dependence of the coefficients more closely
captures temporal variations in natural environments, such as daily or seasonal changes.

It is clear that, for each $T>0$, the corresponding kinetic system of \eqref{eq:main}
\begin{equation}\label{eq:main-kinetic}
	\left\{\begin{array}{l}
\displaystyle		\frac{ d u_1}{dt}=u_1\left(r_{1, T}(t)-a_{1, T}(t) u_1 -k_{1, T}(t)u_2\right),  \vspace{5pt}\\
\displaystyle		\frac{d u_2}{dt}=u_2\left(r_{2, T}(t)-a_{2, T}(t) u_2 -k_{2, T}(t) u_1\right), 
	\end{array}\right.
\end{equation}
admits a trivial steady state $(0,0)$ and two semi-trivial periodic states $(p_{1,T},0)$ and $(0,p_{2,T})$, where for each $i=1,2$,  $p_{i,T}(\cdot)$ is the unique positive solution of the logistic equation
\begin{equation}\label{T-periodic equilibria}
	\left\{\begin{array}{l}
\displaystyle	\frac{d }{d t}p_{i,T}(t)=p_{i,T}(t)\left(r_{i, T}(t)-a_{i, T}(t) p_{i,T}(t)\right),\quad t\in\R,  \vspace{5pt}\\
	p_{i,T}(t)= p_{i,T}(t+T), \quad t\in\R.
\end{array}\right.
\end{equation} 
We emphasize that the function $p_{i,T}(t)$ may not be $T$-rescalings of a common function $p_i$.

We are interested in time periodic traveling waves connecting the two semi-trivial states, 
which is a key notion for describing spatiotemporal evolution in time periodic environments.
Specifically, for a given period $T>0$, a time periodic traveling wave of \eqref{eq:main} connecting $(0,p_{2,T})$ and $(p_{1,T},0)$ is a classical solution $\mathbf{u}_T(t, x)=\left(u_{1, T}(t, x), u_{2, T}(t, x)\right)$ of the type
$$
u_{i, T}(t, x)=\phi_{i, T}\left(x-c_T t, t\right), \quad i=1,2,
$$
where the functions $\phi_i: \mathbb{R} \times \mathbb{R} \rightarrow[0,1]$ are $T$-periodic in their second variable, and satisfy
$$
\lim _{\xi \rightarrow+\infty}\left(\phi_{1, T}(\xi, t), \phi_{2, T}(\xi, t)\right)=(0,p_{2,T}(t)), \quad \lim _{\xi \rightarrow-\infty}\left(\phi_{1, T}(\xi, t), \phi_{2, T}(\xi, t)\right)=(p_{1,T}(t),0),
$$
with the convergences taking place in the topology of $L^{\infty}(\mathbb{R})$. The constant $c_T \in \mathbb{R}$ is called the wave speed. The notion of periodic traveling waves naturally extends the classical one of traveling waves in the autonomous case, where wave propagation phenomena have been widely discussed in the literature (see e.g., the review articles \cite{gw12,girardin19} and references therein).

In the present paper, we focus on time periodic bistable traveling waves, where the limiting states $(p_{1,T},0)$ and $(0,p_{2,T})$ are locally stable solutions of the kinetic system \eqref{eq:main-kinetic}. Before proceeding, let us first review some well known results in the autonomous case, where 
the coefficients $d_i,r_i,a_i,k_i$ are positive constants. In such a homogeneous setting, system \eqref{eq:main} reduces to 
\begin{equation}\label{eq:homogeneous}
	\left\{\begin{array}{ll}
		\partial_t u_1=d_{1} \partial_{xx}u_1+u_1\left(r_{1}-a_{1} u_1 -k_{1}u_2\right), &\hbox{for }\,\, (t,x) \in \R \times\R,  \vspace{5pt}\\
		\partial_t u_2=d_{2} \partial_{xx}u_{2}+u_2\left(r_{2}-a_{2} u_2 -k_{2} u_1\right), &\hbox{for }\,\, (t,x)\in \R\times\R,
	\end{array}\right.
\end{equation}   
and the two semi-trivial equilibria are explicitly given by $(r_1/a_1,0)$ and $(0, r_2/a_2)$. Under the assumption 
\begin{equation}\label{homo-strong}
r_1a_2<k_1r_2\quad \hbox{and} \quad r_2a_1<k_2r_1,
\end{equation}
the kinetic system associated with \eqref{eq:homogeneous} exhibits strong competition dynamics, i.e., the nonnegative equilibria  consist of an unstable node $(0,0)$, two stable nodes $(r_1/a_1,0)$ and $(0, r_2/a_2)$, and a saddle point  corresponding to a coexistence state. 
In this situation, there is a unique speed $c\in\R$ for which system \eqref{eq:homogeneous} admits a  traveling wave connecting $(0, r_2/a_2)$ and $(r_1/a_1,0)$ (see \cite{C2,G1,K1}). Moreover, this wave is unique up to shifts and is globally stable (\cite{K2}).

When the coefficients $d_i,r_i,a_i,k_i$ truly depend on $t$, periodic traveling waves  connecting $(0,p_{2,T})$ and $(p_{1,T},0)$ are known to exist under certain strong restrictions. This problem was first studied by Zhao and Ruan \cite{zr} for the monostable case (where one semi-trivial periodic states is stable and the other is unstable, also referred to as the strong-weak competition case), and subsequently addressed in the bistable setting by Bao and Wang \cite{B}. Specifically, under the assumptions  
\bqq\label{wang-a2}
\bar{r}_1<\min_{s\in\R}\left(\frac{k_1(s)}{a_2(s)}\right)\bar{r}_2,\quad \bar{r}_2<\min_{s\in\R}\left(\frac{k_2(s)}{a_1(s)}\right)\bar{r}_1,
\eqq
and
\bqq\label{wang-a3}
\bar{r}_1+\bar{r}_2>\max_{s\in\R}\left(\frac{k_2(s)}{a_1(s)}\right)\bar{r}_1,\quad \bar{r}_1+\bar{r}_2>\max_{s\in\R}\left(\frac{k_1(s)}{a_2(s)}\right)\bar{r}_2,
\eqq
where $\bar{r}_i=\int_0^1 r_i(t)dt$ for $i=1,2$, the authors of \cite{B} established 
the existence, uniqueness and stability of time periodic bistable traveling waves. The existence result is proved by applying the abstract theory developed for bistable monotone semiflows in \cite{fz}; the proof of  uniqueness and stability relies on a dynamical systems approach from \cite{zhaobook}, combined with a sub- and super-solution technique.
It should be noted that, \eqref{wang-a3} serves as a technical assumption ensuring that the eigenvalue problem related to the linearized system of \eqref{eq:main-kinetic} at $(p_{1,T},0)$ and $(0,p_{2,T})$  admits a positive eigenvalue and the corresponding eigenfunction is positive. This condition was later removed in \cite{ymo}. In contrast,  assumption \eqref{wang-a2} remains essential in \cite{B,ymo}. On the one hand, it implies 
\begin{equation}\label{0-1-stable}
	\int_0^T \left[a_{1,T}(t) p_{1,T}(t)-k_{1,T}(t)p_{2,T}(t)\right]dt<0, 
	\quad
	\int_0^T \left[a_{2,T}(t)p_{2,T}(t)-k_{2,T}(t)p_{1,T}(t)\right]dt <0, 
\end{equation}
which guarantees the local stability of the semi-trivial states $(p_{1,T},0)$ and $(0,p_{2,T})$ (see the discussion in Section 2.2 below). On the other hand, it also ensures the existence of a unique $T$-periodic coexistence state (see also \cite{hess}). In other words, under assumption \eqref{wang-a2}, the time periodic system \eqref{eq:main-kinetic} exhibits a simple bistable structure analogous to that of the kinetic system associated with the homogeneous system \eqref{eq:homogeneous} under condition \eqref{homo-strong}. Without this assumption, however, it remains unclear that whether the temporal heterogeneity could lead to multiple coexistence states. Such a possibility may prevent the existence of periodic traveling waves connecting the semi-trivial states $(0,p_{2,T})$ and $(p_{1,T},0)$.

A similar challenge also arises for the Lotka-Volterra competition system in spatially periodic media, where the existence of traveling waves has been established under a general bistable structure (see \cite{dhy, dlw, girardin, wo}): In addition to the stability of the two semi-trivial states, any coexistence state in the space of spatially periodic functions must be unstable. With such conditions, the abstract existence theory developed in \cite{fz} remains applicable. Explicit and sufficient conditions ensuring this general bistable structure were provided by Girardin \cite{girardin} in the study of segregation phenomena induced by intense interspecific competitions,
and by Ding, Huang, Yu \cite{dhy} for media exhibiting rapid spatial oscillations. Moreover, a similar approach was previously employed in \cite{dhz, ducrot} to show the existence of bistable traveling waves for scalar spatially periodic equations, covering both rapidly and slowly oscillating media. We also note that, the limits of wave speed were determined in \cite{dhz,dll,hps} when the spatial period tends to $0^+$, and in \cite{dhl} when it tends to $+\infty$. 
These limiting results suggest that spatial oscillations have a notable influence on the speed of wave propagation in scalar bistable equations.

Inspired by the aforementioned results, our first objective is to consider the dependence of the existence of periodic traveling waves for system \eqref{eq:main} with respect to the temporal period $T$. We will show the existence of such waves in both rapidly oscillating environments (i.e., when $T > 0$ is small) and slowly oscillating environments (i.e., when $T > 0$ is large). 
Unlike the works \cite{B,ymo} requiring condition \eqref{wang-a2}, we assume instead that, in the rapidly oscillating case a homogenized system derived from \eqref{eq:main} exhibits strong competition dynamics, whereas in the slowly oscillating case, system \eqref{eq:main} with coefficients frozen at an arbitrary time is of the strong competition type. 
Furthermore, we determine the limiting wave speeds as $T\to 0^+$ and as $T\to +\infty$, and provide estimates for the convergence rates.
These estimates show that, compared to the limiting speeds, the deviations of the wave speed for periodic traveling waves are at most linear in certain sense.

Our second objective is to investigate the sign of wave speed for bistable traveling waves in system \eqref{eq:main}. This question is of considerable biological importance in population dynamics, as the sign determines which species becomes dominant and eventually occupies the entire habitat. Numerous results in this direction have been obtained for Lotka-Volterra competition model and other competition models with bistable nonlinearities; see, e.g., \cite{ccs,dl,fz2,gl,gno,gn,K1,lz,mho,myh,tw,wo,xiao} and references therein.
However, the understanding of wave speed sign is far from complete for system \eqref{eq:main} in temporally periodic environments. To our knowledge, only some sufficient conditions for a species to prevail have been provided in a recent paper \cite{myh} by a sub- and super-solution method. In the present work, we focus on a simplified case where the two species share identical diffusion rates and intraspecific competition rates, and derive an explicit criterion for the wave speed sign by comparing the intrinsic growth rates and interspecific competition strengths. Furthermore, based on our explicit formulas for the limiting speeds, we construct an example demonstrating that the wave speed sign varies with the period $T$. This finding reveals that even  under our simplified ecological conditions, temporal heterogeneity can significantly influence competition outcomes, allowing different species to dominate under different periods.

We mention that apart from the aforementioned studies on bistable traveling waves, the existence and qualitative properties of monostable traveling waves connecting two semi-trivial states have also been extensively investigated for Lotka-Volterra competition systems in both homogeneous and heterogeneous environments (see e.g., \cite{blsw,dls,fyz,hosono,Kan97,yz,zr,zr2} and references therein). In addition, many progress has been made in understanding of the long-time dynamics of the corresponding Cauchy problem 
with various types of initial data, covering both monostable and bistable cases (see e.g., \cite{carrere,gl19,llw,pwr,pwz,wxz23} and references therein).

%%%%%%%%%%%%%%%%%%%%%%%%%%%%%%%%%%%%%%%%%%%%%%%%%%%%%%%%%%%%%%%%%%%%%%%%%%%%%%%%%%%%%%%%%%%

\subsection{Periodic traveling wave in rapidly oscillating environments} 

Our first main result is concerned with the existence of periodic traveling wave for system \eqref{eq:main} when $T$ is small, and the convergence of wave speed, as $T\to 0^+$, to that of the following homogenized system
\bqq\label{eq:homogenization limit}
\left\{\begin{array}{ll}
\partial_t u_1=\bar{d}_1 \partial_{x x} u_1+u_1(\bar{r}_1-\bar{a}_1u_1-\bar{k}_1  u_2) & \text { in }\R \times \mathbb{R}, \vspace{5pt}\\
\partial_t u_2=\bar{d}_2 \partial_{x x} u_2+u_2(\bar{r}_2-\bar{a}_2u_2-\bar{k}_2 u_1) & \text { in }\R \times \mathbb{R}. 
\end{array}\right.
\eqq
Hereinafter, for each $i=1,2$, $\bar{d}_i$, $\bar{r}_i$, $\bar{a}_i$ and $\bar{k}_i$ denote, respectively, the arithmetic means of the $1$-periodic functions $d_i(\cdot)$, $r_i(\cdot)$, $a_i(\cdot)$ and $k_i(\cdot)$, that is,
$$
\bar{d}_i=\int_0^1 d_i(t) dt, \quad \bar{r}_i=\int_0^1 r_i(t)dt, \quad \bar{a}_i=\int_0^1 a_i(t)dt \quad\text { and } \quad \bar{k}_i=\int_0^1 k_i(t) dt.
$$
Similarly to the strong competition condition \eqref{homo-strong} to the homogeneous system \eqref{eq:homogeneous}, we assume that the coefficients of system \eqref{eq:homogenization limit} satisfies 
\begin{itemize}
	\item [{\bf (A1)}] $\bar{r}_1 \bar{a}_2<\bar{k}_1\bar{r}_2$ and $\bar{r}_2\bar{a}_1<\bar{k}_2\bar{r}_1$. 
\end{itemize}

It is clear that, under assumption (A1), $(\bar{p}_1,0)$ and $(0,\bar{p}_2)$ are two locally stable equilibria of the kinetic system of \eqref{eq:homogenization limit}, where $\bar{p}_i=\bar{r}_i/\bar{a}_i$, $i=1,2$. It is also easily seen that there exists a unique speed $c_0\in\R$  such that system \eqref{eq:homogenization limit} admits a unique (up to shifts) traveling wave 
$$
\left(u_1, u_2\right)(t,x)=\left(\phi_{1,0}, \phi_{2,0}\right)(x-c_0 t),
$$        
with $0<\phi_{i, 0}<1$ for $i=1,2$, and $\left(\phi_{1,0}, \phi_{2,0}\right)(-\infty)=(\bar{p}_1,0)$, $\left(\phi_{1,0}, \phi_{2,0}\right)(+\infty)=(0,\bar{p}_2)$.
% It is also known that this wave is globally stable and monotone in the sense that  $\phi_{1,0}(\xi)$ (resp. $\phi_{2,0}(\xi)$) is decreasing (resp. increasing) in $\xi \in \mathbb{R}$.

\begin{theo}\label{Theo small T}
 Let  {\rm (A1)} hold. There exists $T_*>0$ such that for any $0<T<T_*$, there is a unique speed $c_T \in \mathbb{R}$ for which system \eqref{eq:main} admits a periodic traveling wave
$\left(u_{1, T}, u_{2, T}\right)(t,x)=\left(\phi_{1, T}, \phi_{2, T}\right)(x-c_T t, t)$
connecting $(0,p_{2,T}(t))$ and $(p_{1,T}(t),0)$. Furthermore, there holds
\bqq\label{small speed cover}
\left|c_T-c_0\right| \leq M_1 T
\eqq
for all small $T>0$, where $M_1>0$ is a constant independent of $T$. 
\end{theo}

Theorem \ref{Theo small T} is parallel to the main results in \cite{dhy}, which concerned traveling waves for Lotka-Volterra competition systems in spatially rapidly oscillating media. The existence proof follows a similar idea: we apply the abstract theory on the existence of bistable traveling waves developed in \cite{fz} by verifying the stability of both semi-trivial periodic states and the instability of all coexistence periodic states. As mentioned earlier, this strategy aligns with previous applications in \cite{dlw, girardin, wo} for spatially periodic problems.

On the other hand, our method to determining the limiting speed and estimating the convergence rate, as stated in \eqref{small speed cover}, is different from that used in \cite{dhy}. Specifically, the proof in \cite{dhy} relies on an application of the implicit function theorem in appropriately chosen Banach spaces, which requires the limiting speed to be non-vanishing.
Their method also provides the convergence of wave profiles as the spatial period tends to $0$. 
However, this technique seems difficult to adapt to our temporally periodic problem.  
Instead, our proof is based on the construction of novel super- and sub-solutions for system \eqref{eq:main} with rapid temporal oscillations, which yields \eqref{small speed cover} regardless of whether the limiting speed $c_0$ vanishes or not,  and does not require any information about the asymptotic behavior of the wave profile. 
This construction involves suitably perturbing the traveling wave of the homogenized system \eqref{eq:homogenization limit}.

%%%%%%%%%%%%%%%%%%%%%%%%%%%%%%%%%%%%%%%%%%%%%%%%%%%%%%%%%%%%%%%%%%%%%%%%%%%%%%%

\subsection{Periodic traveling wave in slowly oscillating environments} 
In the subsection, we present our results on the existence of periodic traveling wave with large $T>0$ and the convergence of wave speed as $T\to +\infty$. To this end, we need to consider the following system with coefficients frozen at $s \in \mathbb{R}$: 
\bqq\label{eq:frozen}
\left\{\begin{array}{ll}
	\partial_t u_1=d_{1}(s) \partial_{xx}u_1+u_1\left(r_{1}(s)-a_{1}(s)u_1-k_{1}(s) u_2\right)   & \text { in }\,\R \times \mathbb{R}, \vspace{5pt}\\
	\partial_t u_2=d_{2}(s) \partial_{xx}u_{2}+u_2\left(r_{2}(s)-a_{2}(s)u_2 -k_{2}(s) u_1\right)  & \text { in }\,\R \times \mathbb{R}. 
\end{array}\right.
\eqq
In \eqref{eq:frozen}, $s$ is viewed as a parameter, since the derivatives only concern the variables $(t,x)$. We assume that for each $s\in\R$, the kinetic system associated with \eqref{eq:frozen} exhibits strong competition dynamics, that is, 
\begin{itemize}
	\item [{\bf (A2)}] 
$r_1(s)a_2(s)<k_1(s)r_2(s)$ and $r_2(s)a_1(s)<k_2(s)r_1(s)$ for all $s\in\R$.  
\end{itemize}
	
Under assumption (A2), it is clear that for each $s\in\R$, system \eqref{eq:frozen} admits a bistable structure with  $(p_1(s),0)$ and $(0,p_2(s))$ being two stable equilibria of the associated kinetic system, where $p_i(s)=r_i(s)/a_i(s)$, $i=1,2$.  Therefore, there exists a traveling wave, denoted by
$$
(u_1, u_2)(t, x; s)=\left(\phi_{1}, \phi_{2}\right)(x-c(s)t;s),
$$
connecting $(0,p_{2}(s))$ and $(p_{1}(s),0)$. Namely, for each $s\in\R$, $(\phi_{1}, \phi_{2})(\xi;s)$ satisfies	
\bqq\label{eq:frozen-wave}
\left\{\begin{array}{l}
d_1(s)  \partial_{\xi\xi}\phi_{1}+c(s) \partial_{\xi} \phi_{1}+\phi_{1}\left(r_{1}(s)-a_{1}(s)\phi_{1} -k_{1}(s) \phi_{2}\right)=0 \quad \text {in }\,\R, \vspace{5pt}\\
 d_2(s) \partial_{\xi\xi}\phi_{2}+c(s)\partial_{\xi} \phi_{2}+\phi_{2}\left(r_{2}(s)-a_{2}(s)\phi_{2} -k_{2}(s)\phi_{1} \right)=0 \quad \text {in }\,\R, \vspace{5pt}\\
 (\phi_{1}, \phi_{2})(-\infty;s)=(p_1(s),0),\quad (\phi_{1}, \phi_{2})(+\infty;s)=(0,p_2(s)).
\end{array}\right.
\eqq
Moreover, the speed $c(s)$ is unique and the  profile $(\phi_{1}, \phi_{2})(\xi;s)$ is unique up to shifts in $\xi$, and both are $1$-periodic in $s$. Denote by $c^*$ the arithmetic mean of the $1$-periodic function $s\mapsto c(s)$, that is, 
\begin{equation}\label{speed-limit-large}
c^*=\int_0^1 c(s) d s. 
\end{equation}

\begin{theo}\label{Theo large T}
Let {\rm (A2)} hold. There exists $T^*>0$ such that for any $T>T^*$, there is a unique $c_T \in \mathbb{R}$ for which system \eqref{eq:main} admits a periodic traveling wave 
$\left(u_{1, T}, u_{2, T}\right)(t,x)=\left(\phi_{1, T}, \phi_{2, T}\right)(x-c_T t, t)$  connecting $(0,p_{2,T}(t))$ and $(p_{1,T}(t),0)$. Furthermore, there holds
\bqq\label{large speed cover}
\left|c_T-c^*\right| \leq \frac{M_2}{T}
\eqq
for all large $T>0$, where $M_2>0$ is a constant independent of $T$. 
\end{theo} 

Several comments are given in order. First of all, Theorem \ref{Theo large T} presents new observations, as to our knowledge, no existing results address the wave propagation  in Lotka-Volterra competition systems under temporally or spatially slow oscillations. Although the proof employs a strategy similar to that of Theorem \ref{Theo small T} for rapidly oscillating environments, the slow oscillation regime introduces greater challenges. For instance, to show the instability of all coexistence periodic states for large $T$, it is necessary to fully characterize the stability of all entire nonnegative solutions (including connecting orbits) of the kinetic system associated with \eqref{eq:frozen} at an arbitrary frozen poin $s$. In contrast, the rapid oscillation case only requires analyzing the stability of nonnegaive equilibria in the kinetic system associated with \eqref{eq:homogenization limit}. A detailed discussion of these methodological contrasts is provided in the opening of Section 4.

Next, we compare our Theorems \ref{Theo small T}-\ref{Theo large T} with the main results in \cite{B,ymo}. 
While conditions (A1) and (A2) are clearly weaker and more natural than \eqref{wang-a2}, our existence results are proved only for either small or large values of $T$. Under assumptions (A1) and (A2), a natural question arises: if $T_*$ and $T^*$ denote, respectively, the supremum and infimum of the values of $T$ for which the conclusions of Theorems \ref{Theo small T}-\ref{Theo large T} hold, is it true that $T_* = T^*$? In other words, does a periodic traveling wave exist for every period $T > 0$? This is a difficult question, due to the complexity of competition dynamics under general periodic oscillations.
We leave it open for future investigation. We also note that, the periodic traveling waves obtained in Theorems \ref{Theo small T}-\ref{Theo large T} are monotone in the spatial variable, unique up to spatial shifts and globally stable (see Propositions \ref{existence-fix} and \ref{stability} below), which aligns with the results in \cite{B,ymo}.

Lastly, the ideas in showing \eqref{small speed cover} and  \eqref{large speed cover} are inspired by a recent work by the first author \cite{D1} on scalar bistable equations of the form 
\begin{equation}\label{scalar-eq}
\partial_t u =\partial_{xx} u+ a(t/T)u(1-u)(u-b(t/T)) \,\,\hbox{ in }\,\, \R^2, 
\end{equation}
where $a(\cdot)$ and $b(\cdot)$ are almost periodic functions satisfying $\inf_{t\in\R}a(t)>0$ and $0<\inf_{t\in\R}b(t)\leq \sup_{t\in\R}b(t)<1$. Although we consider the simpler periodic case here, 
the competition mechanism introduces considerable difficulties that necessitate a more  technically demanding construction of sub- and super-solutions.  Furthermore, while traveling waves in the scalar case \eqref{scalar-eq} connect constant states $0$ and $1$, the waves in our competition model connect periodic states $(p_{1,T},0)$ and $(0,p_{2,T})$. An additional complication arises from the fact that the function $p_{i,T}$ may not simply be $T$-rescalings of a common function. 
This makes the proof of \eqref{large speed cover} in the slow oscillation case particularly more involved. Specifically, we need to establish a uniform convergence of the $T$-rescaled functions of $p_{i,T}(\cdot)$ as $T \to +\infty$, and also estimate the corresponding convergence rate for large $T$.

%%%%%%%%%%%%%%%%%%%%%%%%%%%%%%%%%%%%%%%%%%%%%%%%%%%%%%%%%%%%%%%%%%%%%%%%%%%%%%%%%%%%%%%%%%%
	
\subsection{Sign of wave speeds}

Our last main result is concerned with the sign of wave speeds for time periodic traveling waves of \eqref{eq:main}. 
This question remains largely open even in the homogeneous case, where some interesting partial results have been obtained 
under various constraints on the coefficients.
We mainly investigate the influence of temporal oscillations on the sign of wave speeds by examining the case where both species possess identical diffusion rates and intraspecific competition rates, namely,
\begin{equation}\label{equiv-ad}
d_{1}(t)\equiv d_{2}(t)\quad \hbox{and}\quad a_{1}(t)\equiv a_{2}(t).
\end{equation}
In this special case, the sign of bistable wave speeds for the homogeneous system \eqref{eq:homogeneous} has been completed understood. More precisely, if the coefficients $d_i,r_i,a_i,k_i$ are positive constants satisfying \eqref{homo-strong} and \eqref{equiv-ad}, then it is known from \cite[Theorems 1-2]{gl} (see also \cite{gno}) that
\begin{equation*}
	\left\{\begin{array}{ll}
		c>0 & \displaystyle \hbox{when }\,\, r_1\geq r_2,\,\, k_1\leq k_2, \hbox{ and }\, \frac{r_1}{k_1}\neq \frac{r_2}{k_2}; \vspace{6pt}\\
		c<0 & \displaystyle \hbox{when }\,\, r_1\leq r_2,\,\, k_1\geq k_2, \hbox{ and }\,  \frac{r_1}{k_1}\neq \frac{r_2}{k_2}; \vspace{6pt}\\
		c=0 & \hbox{when }\,\,  r_1= r_2,\,\, k_1= k_2.
	\end{array}\right.
\end{equation*}
In the following theorem, we extend the above criteria for the homogeneous system \eqref{eq:homogeneous} to the time periodic system \eqref{eq:main} with fix period $T>0$. 

\begin{theo}\label{Theorem  sign of speed}
Let $T>0$ be fixed, and \eqref{0-1-stable}, \eqref{equiv-ad} hold. Assume that system \eqref{eq:change of variable} admits a periodic traveling wave $(u_{1, T}, u_{2, T})(t,x)=(\phi_{1,T}, \phi_{2,T})(x-c_Tt, t)$ connecting $(0,p_{2,T}(t))$ and $(p_{1,T}(t),0)$ with speed $c_T\in\R$. Then, the following statements hold true:
\begin{itemize}
\item [{\rm (i)}]
If $k_{1}(t) \leq k_{2}(t)$  and $r_{1}(t) \geq r_{2}(t)$ for all $t \in \mathbb{R}$, then $c_T \geq0$;
\item [{\rm (ii)}]
If $k_{1}(t) \geq k_{2}(t)$  and $r_{1}(t) \leq r_{2}(t)$ for all $t \in \mathbb{R}$, then $c_T \leq0$;
\item [{\rm (iii)}] If the assumptions in {\rm (i)} $($resp. {\rm (ii)}$)$ hold, and $r_1(t)/k_{1}(t) \not\equiv r_2(t)/k_{2}(t)$, then $c_T>0$  $($resp. $c_T <0$$)$. Consequently, if  $k_{1}(t) \equiv k_{2}(t)$ and $r_{1}(t)\equiv r_{2}(t)$ then $c_T=0$.	
\end{itemize}
\end{theo}

We recall that assumption \eqref{0-1-stable} is imposed to guarantee the local stability 
$(p_{1,T},0)$ and $(0,p_{2,T})$ for each $T>0$.
Theorem \ref{Theorem sign of speed} is proved by comparing the periodic traveling wave of \eqref{eq:main} with that of a reflected system, exploiting a symmetry inherent in the competitive structure under assumption \eqref{equiv-ad}.
This approach is inspired by a similar comparison method used for spatially periodic systems in \cite[Proposition 1.2 and Theorem 1.3]{dl}. It is worth noting that, the conditions in statement (iii), which ensure the strict sign of wave speeds, is much weaker than those imposed in \cite{dl}, where it is necessary that the pair of coefficients $k_1(\cdot)$ and $k_2(\cdot)$, as well as $r_1(\cdot)$ and $r_2(\cdot)$, 
be separable by constants. We also mention that \cite{dl} provided an example in which the signs of wave speeds in two opposite directions are different (one is zero and the other one is nonzero), 
while such non-symmetric propagation phenomena never happen in our spatially homogeneous model \eqref{eq:main}. 
On the other hand, for our system, we observe that variations in the temporal period can alter the sign of wave speeds, as illustrated in the following theorem.

\begin{theo}\label{theo-change-sign}
There exist positive $1$-periodic functions $d_i,r_i,a_i,k_i \in C^{2}(\mathbb{R})$ satisfying {\rm (A1)}, {\rm (A2)} and \eqref{equiv-ad}
such that $c_0c^*<0$, where $c_0$ and $c^*$ are, respectively, the limiting speeds obtained in Theorems {\rm \ref{Theo small T}}-{\rm \ref{Theo large T}}. Consequently, there exist $0<T_1<T_2<+\infty$ such that the unique wave speeds $c_{T_1}$ and $c_{T_2}$ of system \eqref{eq:main} with periods $T = T_1$ and $T = T_2$, respectively, have opposite signs. 
\end{theo}

Theorem \eqref{theo-change-sign} is proved by utilizing the exact traveling wave solution of a special homogeneous Lotka-Volterra system provided by Rodrigo and Mimura in \cite{M3}. 
Specifically, by choosing the periodic coefficients $d_i,r_i,a_i,k_i$ such that both systems \eqref{eq:homogenization limit} and \eqref{eq:frozen} reduce to the examples in \cite{M3}, and applying the formulas for the limiting wave speeds derived in Theorems \ref{Theo small T} and \ref{Theo large T}, we are able to compute the limiting speeds explicitly.

Finally, we comment on our results from an ecological perspective. Theorem \ref{Theorem sign of speed} describes an intuitive phenomenon: When both species possess the same dispersal and intraspecific competition abilities, the species with a higher growth rate and a lower interspecific competition rate can successfully invade from right to left and eventually win in the competition. It is important to note that these comparisons of growth rates and interspecific competition rates hold at any time. When this condition fails, more interesting phenomena happen: Different species may become dominant under different periods, as shown in Theorem \ref{theo-change-sign}. This phenomenon highlights the significant influence of temporal variations on the interspecific competition outcomes.

\vskip 10pt

{\bf Outline of this paper.}  
In Section 2, we transform system \eqref{eq:main} into a cooperative diffusion system, and present some results for this system under a fixed period $T>0$, including the comparison principle, as well as the existence, uniqueness, and stability of periodic traveling waves. In the subsequent sections, we work directly with this equivalent cooperative system. Sections 3 and 4 are devoted to the proofs of Theorems \ref{Theo small T} and \ref{Theo large T}, which address the existence of periodic traveling waves and the asymptotic behavior of wave speeds for small and large values of $T$, respectively.  The proofs of Theorems \ref{Theorem sign of speed} and \ref{theo-change-sign} on the sign of wave speeds are given in Section 5. Lastly, Section 6 is an appendix, proving some uniform estimates of the wave profile $(\phi_{1}, \phi_{2})(\xi;s)$ for the homogeneous system \eqref{eq:frozen-wave} with respect to the frozen parameter $s$. These estimates play an important role in showing Theorem \ref{Theo large T}.

%%%%%%%%%%%%%%%%%%%%%%%%%%%%%%%%%%%%%%%%%%%%%%%%%%%%%%%%%%%%%%%%%%%%%%%%%%%%%%%%%%%%%%%%%%%
%%%%%%%%%%%%%%%%%%%%%%%%%%%%%%%%%%%%%%%%%%%%%%%%%%%%%%%%%%%%%%%%%%%%%%%%%%%%%%%%%%%%%%%%%%%

\section{Preliminaries}
In this section, we collect some preliminary materials for later use. Throughout this section, we fix $T>0$ and assume \eqref{0-1-stable} holds.
Let us first set some notations that will be frequently used throughout the remaining parts of this paper. 
Since $d_i,r_i,a_i,k_i \in C^{2}(\mathbb{R})$ are positive 1-periodic functions for $i=1,2$, there exist two constants $\theta_+>\theta_->0$
such that 
\begin{equation}\label{bound-rak}
\theta_-\leq \min\left\{d_i(t),\, r_i(t),\, a_i(t),\,k_i(t)\right\} \leq \max\left\{d_i(t),\, r_i(t),\, a_i(t),\, k_i(t)\right\} \leq \theta_+ \quad \hbox{for all } \, t\in\R.
\end{equation}
Moreover, remember that for each $T>0$, $p_{i,T}$ is the unique positive solution of \eqref{T-periodic equilibria}. A simple comparison argument implies that 
\begin{equation}\label{bound-pT}
	\gamma_-\leq 	\frac{\min_{s\in\R} r_i(s)}{\max_{s\in\R} a_i(s)} \leq p_{i,T}(\cdot) \leq \frac{\max_{s\in\R} r_i(s)}{\min_{s\in\R} a_i(s)} \leq \gamma_+  \quad \hbox{in }\, \R,
\end{equation}
for some constants $\gamma_+>\gamma_->0$ independent of $T$. 

\subsection{Comparison principle for cooperative systems}
We begin by transforming the competitive system \eqref{eq:main} into a cooperative system. 
Let $(v_{1,T},v_{2,T})(t,x)$ be defined as follows:
\bqq\label{change}
v_{1,T}(t, x)=\frac{u_{1,T}(t, x)}{p_{1,T}(t)}, \quad v_{2,T}(t, x)=\frac{p_{2,T}(t)-u_{2,T}(t, x)}{p_{2,T}(t)}.
\eqq
Substituting \eqref{change} into \eqref{eq:main}, we obtain
\bqq\label{eq:change of variable}
\left\{\begin{array}{ll}
\partial_t v_{1,T}=d_{1,T}(t)\partial_{xx} v_{1,T}+g_{1,T}\left(t, v_{1,T}, v_{2,T}\right), & (t,x)\in \R \times \mathbb{R},\vspace{5pt} \\
\partial_t v_{2,T}= d_{2,T}(t)\partial_{xx}v_{2,T}+g_{2,T}\left(t, v_{1,T}, v_{2,T}\right), & (t,x)\in \R \times \mathbb{R},
\end{array}\right.
\eqq
where  
\begin{equation}\label{def-giT}
\left\{\begin{aligned}
	g_{1,T}\left(t, v_{1}, v_{2}\right)&=a_{1,T}(t) p_{1,T}(t)(1-v_{1})v_{1}-k_{1,T}(t)p_{2,T}(t)v_{1}(1-v_{2}),\vspace{5pt} \\
	g_{2,T}\left(t, v_{1}, v_{2}\right)&=-a_{2,T}(t) p_{2,T}(t)(1-v_{2})v_{2}+k_{2,T}(t)p_{1,T}(t)v_{1}(1-v_{2}).
\end{aligned}\right.
\end{equation}
Clearly, the states $(0,p_{2,T})$, $(p_{1,T},0)$ and $(0,0)$ of \eqref{eq:main} become, respectively, 
$$
\mathbf{0}:=(0,0),\quad \mathbf{1}:=(1,1), \quad \mathbf{e}_0:=(0,1). 
$$
Then, finding a periodic traveling wave $(u_1,u_2)(t,x)=(\phi_1,\phi_2)(x-c_Tt,t)$ of \eqref{eq:main} connecting $(0, p_{2,T})$ and $(p_{1,T}, 0)$ amounts to finding a periodic traveling wave  $(v_1,v_2)(t,x)=(\varphi_1,\varphi_2)(x-c_Tt,t)$ of \eqref{eq:change of variable} 
connecting $\mathbf{0}$ and $\mathbf{1}$.  
It is readily seen that, $(\varphi_{1,T},\varphi_{2,T})\in C(\R^2;[0,1]^2)$ satisfies 
\begin{equation}\label{eq:change-tw}
\left\{\begin{array}{l}
	\partial_t \varphi_{1,T}=d_{1,T}(t)\partial_{\xi \xi} \varphi_{1,T}+c_T\partial_\xi \varphi_{1,T}+g_{1,T}\left(t, \varphi_{1,T}, \varphi_{2,T}\right) \,\, \text { in }\,  \R\times \mathbb{R}, \vspace{5pt} \\
	\partial_t \varphi_{2,T}= d_{2,T}(t)\partial_{\xi \xi}\varphi_{2,T}+c_T\partial_\xi \varphi_{2,T}+g_{2,T}\left(t, \varphi_{1,T}, \varphi_{2,T}\right) \,\,  \text { in }\,  \R\times \mathbb{R}, \vspace{5pt} \\
	(\varphi_1,\varphi_2)(\xi,t) \hbox{ is $T$-periodic in } t\in\R, \vspace{5pt} \\
	(\varphi_1,\varphi_2)(-\infty,t)=\mathbf{1}\quad \hbox{and}\quad (\varphi_1,\varphi_2)(+\infty,t)=\mathbf{0} \,\, \hbox{ uniformly in } \,t\in\R.
\end{array}\right.
\end{equation}

Denote by $\mathcal{C}$ the Banach space of all bounded and uniformly continuous functions from $\R$ to $\R^2$ with the supremum norm. We recall some order relations in $\mathcal{C}$. For any $\mathbf{u}=\left(u_1, u_2\right)$, $\mathbf{v}=\left(v_1, v_2\right)$ in $\mathcal{C}$, we say that $\mathbf{u} \geq \mathbf{v}$ if $u_1(x) \geq v_1(x)$ and $u_2(x) \geq v_2(x)$ for all $x\in\R$; $\mathbf{u}>\mathbf{v}$ if $\mathbf{u} \geq \mathbf{v}$ but $\mathbf{u}$ is not identically equal to $\mathbf{v}$; $\mathbf{u} \gg \mathbf{v}$ if $u_1(x)>v_1(x)$ and $u_2(x)>v_2(x)$ for all $x\in\R$. 

In the sequel, we denote
$\mathcal{C}_{\mathbf{1}}=\{ \mathbf{u} \in \mathcal{C}: \mathbf{1} \geq  \mathbf{u} \geq \mathbf{0} \}$.
By definition, any periodic traveling wave $(v_1,v_2)(t,x)$ of \eqref{eq:change of variable} connecting  $\mathbf{0}$ and $\mathbf{1}$ satisfies $(v_1,v_2)(t,\cdot)\in \mathcal{C}_{\mathbf{1}}$ for all $t\in\R$. 
Furthermore, the strong maximum principle of parabolic equations implies that 
$\mathbf{1} \gg (v_1,v_2)(t,\cdot) \gg\mathbf{0}$ for all $t\in\R$.  

\begin{defi}\label{defi sup sub}
	A pair of functions $(v_1, v_2) \in C^{1,2}\left((0, \infty) \times \mathbb{R}; \mathbb{R}^2\right)$ is said to be a super-solution $($resp. sub-solution$)$ of \eqref{eq:change of variable} if
	$$
	\left\{\begin{array}{l}
		\partial_t v_1-d_{1,T}(t)\partial_{x x} v_1-g_{1,T}\left(t, v_1, v_2\right) \geq 0 \,\, (\text{resp.} \,\leq 0) \quad \text { in } \, (0, \infty) \times \mathbb{R}, \vspace{5pt} \\
		\partial_t v_2-d_{2,T}(t)\partial_{x x} v_2-g_{2,T}\left(t, v_1, v_2\right) \geq 0\,\, (\text{resp.}\, \leq 0) \quad\text { in }\, (0, \infty) \times \mathbb{R} .
	\end{array}\right.
	$$	
Furthermore, a super-solution $($resp. sub-solution$)$ is called strict if it is not a solution. 
\end{defi}

The comparison principle possessed by system \eqref{eq:change of variable} can be presented as follows. 

\begin{lem}\label{comparison principle}
	Fix a period $T>0$. Let $(v_1^{+}, v_2^{+})$ $($resp. $(v_1^{-}, v_2^{-})$$)$ be a super-solution $($resp. sub-solution$)$ of \eqref{eq:change of variable} with initial data $(v_{1,0}^{+}, v_{2,0}^{+}) \in \mathcal{C}$ $($resp. $(v_{1,0}^{-}, v_{2,0}^{-}) \in \mathcal{C}$$)$. Suppose that either $(v_1^{+}, v_2^{+})(t,\cdot) \in \mathcal{C}_{\mathbf{1}}$ or $ (v_1^{-}, v_2^{-})(t,\cdot) \in \mathcal{C}_{\mathbf{1}}$ holds for all $t \geq 0$. Then, the following statements hold true:
	\begin{enumerate}
		\item [{\rm (i)}] If $(v_{1,0}^{+}, v_{2,0}^{+}) \geq (v_{1,0}^{-}, v_{2,0}^{-})$, then $(v_1^{+}, v_2^{+})(t, \cdot) \geq(v_1^{-}, v_2^{-})(t, \cdot)$ for any $t>0$;
		\item [{\rm (ii)}] If $\mathbf{1} \gg (v_{1,0}^{+}, v_{2,0}^{+})>(v_{1,0}^{-}, v_{2,0}^{-}) \gg \mathbf{0}$, then $(v_1^{+}, v_2^{+})(t,\cdot) \gg (v_1^{-}, v_2^{-})(t, \cdot)$ for any $t>0$.
	\end{enumerate}
\end{lem} 

The proof of Lemma \ref{comparison principle} follows from that of \cite[Lemma 2.2]{dl} with a minor modification.

%%%%%%%%%%%%%%%%%%%%%%%%%%%%%%%%%%%%%%%%%%%%%%%%%%%%%%%%%%%%%%%%%%%%%%%%%%%%%%%%%%%%%%%%%%%%%%%%%%%%%%%%%%%%%%%%%%%%

\subsection{Existence of periodic traveling waves for fixed $T>0$}
In this subsection, we employ the abstract theory in \cite{fz} to prove the existence of periodic traveling waves for the cooperative system \eqref{eq:change of variable} with fixed $T>0$. As mentioned earlier, this approach has been adopted in \cite{B,ymo} under the additional assumption \eqref{wang-a2}, which ensures that the kinetic system of \eqref{eq:change of variable}
\bqq\label{eq:variable-kinetic}
\left\{\begin{array}{ll}
	\displaystyle \frac{dv_{1}}{dt}=g_{1,T}\left(t, v_{1}, v_{2}\right), & t \in \R,\vspace{5pt} \\
	\displaystyle \frac{dv_{2}}{dt}=g_{2,T}\left(t, v_{1}, v_{2}\right), & t \in \R,
\end{array}\right.
\eqq
has a unique periodic solution strictly between $\mathbf{0}$ and $\mathbf{1}$. 
Without such an assumption, we will show that if all the periodic solutions lying strictly between $\mathbf{0}$ and $\mathbf{1}$ are linearly unstable, then the existence results in \cite{B,ymo} remain valid.

To apply the existence result in \cite{fz}, let us introduce some notations. For any $t\geq 0$, define $Q_t:\mathcal{C}_{\mathbf{1}}\to \mathcal{C}_{\mathbf{1}}$ by 
$$Q_t[(v_{1,0},v_{2,0})](\cdot)=(v_1,v_2)(t,\cdot;(v_{1,0},v_{2,0})), $$
where $(v_1,v_2)(t,\cdot;(v_{1,0},v_{2,0}))$ denotes the solution of the Cauchy problem of \eqref{eq:change of variable} with initial data $(v_{1,0},v_{2,0})$. Clearly, the family of mappings $\{Q_t\}_{t\in \R^+}$ forms a $T$-periodic semiflow on the metric space $\mathcal{C}$ in the sense that 
$Q_{0}[\boldsymbol\varphi]=\boldsymbol\varphi$ for all $\boldsymbol\varphi\in\mathcal{C}_{\mathbf{1}}$;
 $Q_{t}\circ Q_{T}[\boldsymbol\varphi]= Q_{t+T}[\boldsymbol\varphi]$ for all $\boldsymbol\varphi\in\mathcal{C}_{\mathbf{1}}$, $t\geq 0$;
 $Q_{t}[\boldsymbol\varphi]$ is jointly continuous in $(t,\boldsymbol\varphi)$ on $[0, +\infty)\times \mathcal{C}_{\mathbf{1}}$.
It is well known that the map $Q_{T}$ is termed as the Poincar\'e  map associated with this periodic semiflow. Denote by $E\subset [0,1]^2$ the set of all spatially homogeneous fixed pints of $Q_T$. Clearly, $\mathbf{0}$,  $\mathbf{1}$, $\mathbf{e}_0\in E$.  

We recall here some notions on the stability of the fixed points of $Q_T$.
\begin{defi}\label{defi-stability}
For the map $Q_T: [0,1]^2\to [0,1]^2$, a fixed point $\boldsymbol{\alpha} \in E$ is said to be strongly stable $($resp. unstable$)$ from below if there exist a number $\delta>0$ and a vector $\mathbf{e} \in (0,+\infty)\times (0,+\infty)$ such that for any $\epsilon \in (0,\delta]$, 
$$Q_T[\boldsymbol{\alpha}-\epsilon \mathbf{e}]\gg \boldsymbol{\alpha}-\epsilon \mathbf{e} \quad (resp. \quad Q_T[\boldsymbol{\alpha}-\epsilon \mathbf{e}]\ll \boldsymbol{\alpha}-\epsilon \mathbf{e}).$$
Strong stability $($resp. unstability$)$ from above can be defined in a similar way. Moreover, $\boldsymbol{\alpha} \in E$ is said to be strongly stable $($resp. unstable$)$ uniformly with respective to $T>0$ if $\delta$ and $\mathbf{e}$ can be chosen independent of $T$.
\end{defi}

We first consider the stability of the fixed points  $\mathbf{0}$ and $\mathbf{1}$. 
Let $\bar{\mathbf{v}}=\left(\bar{v}_1, \bar{v}_2\right)$ be an arbitrary $T$-periodic of \eqref{eq:variable-kinetic}. 
The corresponding linearized system of \eqref{eq:variable-kinetic}  at $\bar{\mathbf{v}}$ is 
\bqq\label{linearized-0}
\frac{d \mathbf{w}}{dt}= A_T(\bar{\mathbf{v}}) \mathbf{w} \,\,\hbox{ for }\, t\in\R,
\eqq
where $\mathbf{w}=(w_1,w_2)\in C^1(\R;\R^2)$, and 
\begin{equation}\label{matrxi-A}
\begin{array}{l}
	A_T(\overline{\mathbf{v}}) = \vspace{7pt} \\
	\left(\begin{array}{cc}
		a_{1,T}p_{1,T}-2a_{1,T}p_{1,T}\bar{v}_1-k_{1,T}p_{2,T}\left(1-\bar{v}_2\right) & k_{1,T}p_{2,T}\bar{v}_1 \vspace{5pt} \\
		k_{2,T}p_{1,T}\left(1-\bar{v}_2\right) & -a_{2,T}p_{2,T}+2a_{2,T}p_{2,T}\bar{v}_2-k_{2,T}p_{1,T}\bar{v}_1
	\end{array}\right). \end{array}
\end{equation}
In particular, one has
$$
A_T(\mathbf{0})=\left(\begin{array}{cc}
	a_{1,T}p_{1,T}-k_{1,T}p_{2,T} & 0 \vspace{5pt} \\
	k_{2,T}p_{1,T} & -a_{2,T}p_{2,T}
\end{array}\right),
$$
and 
$$
A_T(\mathbf{1})=\left(\begin{array}{cc}
	-a_{1,T}p_{1,T} &  k_{1,T}p_{2,T} \vspace{5pt} \\
	0 & a_{2,T}p_{2,T}-k_{2,T}p_{1,T}
\end{array}\right).
$$
Note that at $\bar{\mathbf{v}}=\mathbf{0}$ and $\bar{\mathbf{v}}=\mathbf{1}$, the linear system \eqref{linearized-0} is not irreducible, and hence, the corresponding linear periodic system may not admit a principal Floquet multiplier. However, thanks to the special triangular form of $A_T(\mathbf{0})$, by using the first inequality of \eqref{0-1-stable} and the fact that $-\int_0^Ta_{2,T}(t)p_{2,T}(t)dt=-1/T \int_0^1r_2(t)dt<0$, one can conclude from the proof of \cite[Theorem 1.2]{ymo} that, for the map $Q_T: [0,1]^2\to [0,1]^2$,  the fixed point $\mathbf{0}$ is strongly stable from above. Similarly, the second inequality of \eqref{0-1-stable} implies that  $\mathbf{1}$ is  strongly stable from below. As a consequence, we have the following lemma. 

\begin{lem}\label{lem-stable-0-1}
Let \eqref{0-1-stable} hold. Then, the fixed points $\mathbf{0}$ and $\mathbf{1}$ are strongly stable from above and below, respectively, for the map $Q_T: [0,1]^2\to[0,1]^2$.  
\end{lem}

We will show in Section {\rm 3.1} (resp. Section {\rm 4.1}) that, 
under assumption {\rm (A1)} (resp. {\rm (A2)}), 
\eqref{0-1-stable} holds for all small $T>0$ (resp. all large $T>0$), and the strong stability of $\mathbf{0}$ and $\mathbf{1}$ is actually uniform in these regimes.

Now, we turn to consider the intermediate fix point $\boldsymbol{\alpha} \in E$ with $\mathbf{1} \gg \boldsymbol{\alpha} \gg\mathbf{0}$. Let $\bar{\mathbf{v}}_{\boldsymbol{\alpha}}(t)$ be the $T$-periodic solution of \eqref{eq:variable-kinetic} with $\bar{\mathbf{v}}_{\boldsymbol{\alpha}}(0)=\boldsymbol{\alpha}$. It is easily seen that  $\mathbf{1}\gg \bar{\mathbf{v}}_{\boldsymbol{\alpha}}(t)\gg \mathbf{0}$ for all $t\in\R$. 
Since the off-diagonal entries of $A_T(\bar{\mathbf{v}}_{\boldsymbol{\alpha}})$ are positive for any $t \in \mathbb{R}$, it is then known from the Krein-Rutmann theory that 
the periodic eigenvalue problem 
\bqq\label{inter-eigen}
\left\{\begin{array}{l}
\displaystyle \frac{d \boldsymbol \psi}{dt}- A_T(\bar{\mathbf{v}}_{\boldsymbol{\alpha}}) \boldsymbol \psi= \lambda \boldsymbol \psi \,\,\hbox{ for }\, t\in\R,\vspace{5pt} \\
\boldsymbol \psi (t+T)=\boldsymbol \psi(t)\,\,\hbox{ for }\, t\in\R,
\end{array}\right.
\eqq
admits a principal Floquet multiplier $\lambda_T(\bar{\mathbf{v}}_{\boldsymbol{\alpha}})$ with a unique (up to multiplication) eigenfunction $\boldsymbol \psi_{\boldsymbol{\alpha}} \gg(0,0)$.

\begin{lem}\label{inter-unstable}
Let $\boldsymbol{\alpha} \in \tilde{E}:=\{\boldsymbol{\alpha} \in E: \mathbf{1} \gg \boldsymbol{\alpha} \gg\mathbf{0}\}$. If $\boldsymbol{\alpha}$ is linearly unstable, that is, $\lambda_T(\bar{\mathbf{v}}_{\boldsymbol{\alpha}})<0$, then 
$\boldsymbol \alpha$ is strongly unstable from above and below.
\end{lem}

\begin{proof}
Since $\lambda_T(\bar{\mathbf{v}}_{\boldsymbol{\alpha}})<0$, by standard comparison arguments, one can find some small $\delta>0$ such that for any $\epsilon \in (0,\delta]$, 
$\bar{\mathbf{v}}_{\boldsymbol{\alpha}}(t)-\epsilon{\boldsymbol \psi_{\boldsymbol{\alpha}}}(t)$ and $\bar{\mathbf{v}}_{\boldsymbol{\alpha}}(t)+\epsilon{\boldsymbol \psi_{\boldsymbol{\alpha}}}(t)$ are, respectively, super-solution and sub-solution for the problem \eqref{eq:variable-kinetic}. 
It further follows from Lemma \ref{comparison principle} and the periodicity of $\bar{\mathbf{v}}_{\boldsymbol{\alpha}}(t)$ and $\boldsymbol \psi_{\boldsymbol{\alpha}}(t)$ that
$$Q_T[\boldsymbol{\alpha}-\epsilon\boldsymbol \psi_{\boldsymbol{\alpha}}(0)] \ll \boldsymbol{\alpha}-\epsilon\boldsymbol \psi_{\boldsymbol{\alpha}}(0) \quad \hbox{and}\quad 
Q_T[\boldsymbol{\alpha}+\epsilon\boldsymbol \psi_{\boldsymbol{\alpha}}(0)] \gg \boldsymbol{\alpha}+\epsilon\boldsymbol\psi_{\boldsymbol{\alpha}}(0).
$$
This ends the proof of Lemma \ref{inter-unstable}. 
\end{proof}

Based on the above preparations, we are ready to present the main result of this subsection. 
\begin{prop}\label{existence-fix}
Let $T>0$ be fixed and \eqref{0-1-stable} hold. Assume that $\lambda_T(\bar{\mathbf{v}}_{\boldsymbol{\alpha}})<0$ for all $\boldsymbol{\alpha} \in \tilde{E}$. 
Then, system \eqref{eq:change of variable} admits a periodic traveling wave
$\left(v_{1}, v_{2}\right)(t,x)=\left(\varphi_{1,T}, \varphi_{2,T}\right)(x-c_Tt, t)$
connecting  $\mathbf{0}$ and $\mathbf{1}$. Furthermore, $\varphi_{i,T}(\xi,t)$ is decreasing in $\xi\in\R$ for each $i=1,2$.  
\end{prop}

\begin{proof}
According to \cite[Theorem 3.3]{fz}, to obtain the existence of a periodic traveling wave connecting $\mathbf{0}$ and $\mathbf{1}$, it suffices to verify that the map $Q_T: \mathcal{C}_{\mathbf{1}}\to \mathcal{C}_{\mathbf{1}}$ satisfies the following hypotheses. 
\begin{itemize}
\item [(H1)] (Translation invariance) $T_y \circ Q_T[\boldsymbol\varphi]= Q_T\circ T_y[\boldsymbol\varphi]$ for all $\boldsymbol\varphi \in \mathcal{C}_{\mathbf{1}}$ and $y\in\R$, where $T_y$ is the translation operator defined by $T_y[\boldsymbol\varphi](x):=\boldsymbol\varphi(x-y)$.
\item [(H2)] (Continuity) $Q_T: \mathcal{C}_{\mathbf{1}}\to \mathcal{C}_{\mathbf{1}}$  is continuous with respect to the compact open topology. 
\item [(H3)] (Monotonicity) $Q_T$ is order-preserving in the sense that $Q_T[\boldsymbol\varphi]\geq Q_T[\boldsymbol\phi]$ whenever $\boldsymbol\varphi \geq \boldsymbol\phi$ in $\mathcal{C}_{\mathbf{1}}$. 
\item [(H4)] (Compactness) $Q_T: \mathcal{C}_{\mathbf{1}}\to \mathcal{C}_{\mathbf{1}}$ is compact with respect to the compact open topology.
\item [(H5)] (Bistability) The fixed points $\mathbf{0}$ and $\mathbf{1}$ are strongly stable from above and below, respectively, and the set $E\setminus \{\mathbf{0}, \mathbf{1}\}$ is totally unordered in the sense that no two elements are ordered.  
\item [(H6)] (Counter-propagation) For each $\boldsymbol{\alpha} \in E\setminus \{\mathbf{0}, \mathbf{1}\}$, there holds $c^*_-(\boldsymbol{\alpha},\mathbf{1})+c^*_+(\mathbf{0},\boldsymbol{\alpha})>0$, where $c^*_-(\boldsymbol{\alpha},\mathbf{1})$ (resp. $c^*_+(\mathbf{0},\boldsymbol{\alpha})$) is the leftward (resp. rightward) spreading speed of the map $Q_T$ restricted on the subspace $\{ \mathbf{u} \in \mathcal{C}: \mathbf{1} \geq  \mathbf{u} \geq \boldsymbol{\alpha} \}$ (resp. $\{ \mathbf{u} \in \mathcal{C}: \boldsymbol{\alpha} \geq  \mathbf{u} \geq \mathbf{0} \}$). We refer to \cite{B,fz} for the detailed definitions of $c^*_-(\boldsymbol{\alpha},\mathbf{1})$ and $c^*_+(\mathbf{0},\boldsymbol{\alpha})$.
\end{itemize}
It is standard to check that the hypotheses (H1)-(H4) are indeed satisfied. Next, we verify (H5). It follows directly from Lemma \ref{lem-stable-0-1} that $\mathbf{0}$ and $\mathbf{1}$ are strongly stable from above and below, respectively.  Moreover, Lemma \ref{inter-unstable}, together with 
Lemma \ref{comparison principle} (ii) and the proof of \cite[Proposition 2.1]{fz} implies that the set $\tilde{E}$ is totally unordered. Due to $\tilde{E}\cup\{\mathbf{e}_0\} = E\setminus \{\mathbf{0}, \mathbf{1}\}$, this immediately gives that the set $ E\setminus \{\mathbf{0}, \mathbf{1}\}$ is totally unordered. Thus, (H5) is satisfied. Since $\lambda_T(\bar{\mathbf{v}}_{\boldsymbol{\alpha}})<0$ for all $\alpha \in \tilde{E}$, the verification of (H6) follows from essentially the same arguments as those used in the proof of \cite[Theorem 3.2]{B}. Thus, 
a periodic traveling wave $\left(v_{1}, v_{2}\right)(t,x)=\left(\varphi_{1,T}, \varphi_{2,T}\right)(x-ct,t)$ is obtained. 
Finally, it also follows from \cite[Theorem 3.3]{fz} that $\varphi_{1,T}(\xi,t)$ and $\varphi_{2,T}(\xi,t)$ are nonincreasing in $\xi \in\R$. Applying the strong maximum of parabolic equations, one can conclude that 
 $\partial_{\xi}\varphi_{i,T}(\xi,t)<0$ in $\R^2$ for each $i=1,2$. This ends the proof of Proposition \ref{existence-fix}.
\end{proof}

%%%%%%%%%%%%%%%%%%%%%%%%%%%%%%%%%%%%%%%%%%%%%%%%%%%%%%%%%%%%%%%%%%%%%%%%%%%%%%%%%%%%%%%%%%%%%%%%%%%%%%%%%%%%%%%%%%

\subsection{Uniqueness and stability of periodic traveling waves for fixed $T>0$}

In this subsection, we present the uniqueness and global stability of periodic traveling wave connecting  $\mathbf{0}$ and $\mathbf{1}$ of system \eqref{eq:change of variable} with fixed $T>0$, which will play an important role in determining the limits of wave speed as $T\to 0^+$ and $T\to +\infty$.  

\begin{prop}\label{stability}
	Let $T>0$ be fixed and \eqref{0-1-stable} hold. Assume that system \eqref{eq:change of variable} admits a periodic traveling wave $\left(\varphi_{1,T}, \varphi_{2,T}\right)(x-c_Tt, t)$
	connecting  $\mathbf{0}$ and $\mathbf{1}$, and that $\varphi_{i,T}(\xi,t)$ is decreasing in $\xi\in\R$ for $i=1,2$.  
	Then, the following statements hold true.
	\begin{itemize}
	\item [{\rm (i)}] $($Global stability$)$ There exists a small constant $\delta_0\in(0,1/2)$ such that for any pair of functions $(v_{1,0},v_{2,0}) \in C(\R;\R^2)$ satisfying
	$$
	\varphi_{i,T}\left( x-\xi^1,0\right)-\delta_0 \leq v_{i,0}(x) \leq \varphi_{i,T}\left( x-\xi^2,0\right)+\delta_0\quad \text {in }\, x \in \mathbb{R},
	$$
	for $i=1,2$ and some constants $\xi^1 \leq \xi^2$ in $\mathbb{R}$, there holds
	$$
	\sup _{x \in \mathbb{R},\,i=1,2}\left|v_{i}(t,x;(v_{1,0},v_{2,0}))-\varphi_{i,T}\left(x-c_Tt-\xi_T^*,t\right)\right| \rightarrow 0 \quad \text {as} \quad t \rightarrow+\infty,
	$$
	where $(v_1,v_2)(t,x,(v_{1,0},v_{2,0}))$ is the solution of the Cauchy problem of \eqref{eq:change of variable} with initial data $(v_{1,0},v_{2,0})$, and $\xi_T^* \in \mathbb{R}$ is a constant depending on $T$ and $(v_{1,0},v_{2,0})$.  
		
	\item [{\rm (ii)}] $($Uniqueness$)$ If $\left(\tilde{\varphi}_{1,T}, \tilde{\varphi_{2,T}}\right)(x-\tilde{c}_Tt, t)$ is another periodic traveling wave connecting  $\mathbf{0}$ and $\mathbf{1}$, then there exits $\xi_0\in\R$ such that 
	$$  \left(\varphi_{1,T}, \varphi_{2,T}\right) (\cdot,\cdot)\equiv \left(\tilde{\varphi}_{1,T}, \tilde{\varphi}_{2,T}\right)(\cdot+\xi_0,\cdot)\, \hbox{ in }\,\R^2\quad  \hbox{and}\quad c_T=\tilde{c}_T.$$
	\end{itemize}
\end{prop}

Statement (ii) is a direct consequence from (i). The proof of (i) adapts the key arguments from \cite[Theorem 4.4]{B} with necessary adjustments (a similar claim was made in \cite{ymo} without proof). Specifically, \cite[Theorem 4.4]{B} relies on the technical assumption \eqref{wang-a3}, which is used to construct super- and sub-solutions of \eqref{eq:change of variable} via the periodic traveling wave in \cite[Lemma 3.4]{B}. 
In Lemma \ref{fix T sup sub} below, we refine \cite[Lemma 3.4]{B} by removing this technical assumption. 
We point out that, the super- and sub-solutions constructed in Lemma \ref{fix T sup sub} will also play a crucial role in proving Theorem \ref{Theorem  sign of speed} in Section 5. Moreover, the method of proof will serve as an archetype to get more involved comparisons in Sections 3-4.    

To present our super- and sub-solutions, we need to introduce a few more notations which will also be frequently used in subsequent sections.  For each $i=1,2$, we denote
\begin{equation}\label{seperate-lambda}
\lambda_{i,T}^{-}:=- \frac{1}{T} \int_0^T \partial_{v_i}g_{i,T}(t,0,0)dt  \quad\hbox{and}\quad
\lambda_{i,T}^{+}:=- \frac{1}{T} \int_0^T \partial_{v_i}g_{i,T}(t,1,1)dt. 
\end{equation}
It is clear that  
$$\lambda_{1,T}^{+}=\frac{1}{T} \int_0^T a_{1,T}(t)p_{1,T}(t)dt>0,\quad \lambda_{2,T}^{-}=\frac{1}{T} \int_0^T a_{2,T}(t)p_{2,T}(t)dt>0, $$
and that \eqref{0-1-stable} implies 
$$\lambda_{1,T}^{-}=-\frac{1}{T} \int_0^T a_{1,T}(t)p_{1,T}(t)-k_{1,T}(t)p_{2,T}(t)dt>0,$$
and
$$\lambda_{2,T}^{+}=- \frac{1}{T} \int_0^T  a_{2,T}(t)p_{2,T}(t)-k_{2,T}(t)p_{1,T}(t) dt>0.$$
Furthermore, by \eqref{bound-pT}, one finds a positive constant $\mu$ independent of $T>0$ such that 
\bqq\label{set-mu}
0<\mu\leq \min _{i=1,2}\left\{\frac{\lambda_{i,T}^{+}}{2},\frac{\lambda_{i,T}^{-}}{2}\right\}.
\eqq 
It is also easily seen that, for each $i=1,2$, the solutions of the following linear  equations
\begin{equation}\label{eigenfunction-0}
\psi'- \partial_{v_i}g_{i,T}(t,0,0) \psi=\lambda_{i,T}^{-}\psi  \text { in } \mathbb{R}, \quad \psi>0 \text { in } \mathbb{R}, \quad \psi \text { is } T \text {-periodic},
\end{equation}	
and
\begin{equation}\label{eigenfunction-1}
	\psi'- \partial_{v_i}g_{i,T}(t,1,1) \psi=\lambda_{i,T}^{+}\psi  \text { in } \mathbb{R}, \quad \psi>0 \text { in } \mathbb{R}, \quad \psi \text { is } T \text {-periodic},
\end{equation}	
exist, denoted by  $\psi_{i,T}^-(t)$ and $\psi_{i,T}^+(t)$, respectively. As a matter of fact, there holds 
\bqq\label{explicit-eigenfunction}
\left\{\begin{array}{l}
	\displaystyle\psi_{i,T}^{-}(t)=\psi_{i,T}^{-}(0){\rm exp}\left\{T\int_{0}^{t/T}\left(\lambda^-_{i,T}+\partial_{v_i}g_{i,T}(sT,0,0)\right)  ds\right\}\vspace{5pt} \\
	\displaystyle \psi_{i,T}^{+}(t)=\psi_{i,T}^{+}(0){\rm exp}\left\{T\int_{0}^{t/T}\left(\lambda^+_{i,T}+\partial_{v_i}g_{i,T}(sT,1,1)\right)  ds\right\}
\end{array}\right.  \quad \hbox{ for }\, t\in\R.
\eqq
It is clear that $\psi_{i,T}^-(t)$ and $\psi_{i,T}^+(t)$ are unique up to multiplication by positive constants. In several places below, we will normalize them as follows 
\bqq\label{normalize-supsub}
\left\|\psi_{1,T}^{+}\right\|=1, \quad\left\|\psi_{2,T}^{-}\right\|=1, \quad\left\|\frac{k_{1}p_{2,T} \psi_{2,T}^{+}}{\psi_{1,T}^{+}}\right\|=\frac{\mu}{4} \quad \text { and } \quad\left\|\frac{k_{2}p_{1,T} \psi_{1,T}^{-}}{\psi_{2,T}^{-}}\right\|=\frac{\mu}{4},
\eqq
where $\|\cdot\|$ denotes the  $L^{\infty}$-norm in $C(\mathbb{R})$.

Let $\rho \in C^2(\mathbb{R})$ be a nonnegative function satisfying
\begin{equation}\label{function-rho}
\rho(\xi)=0 \hbox{ in } [2, +\infty),\quad   \rho(\xi)=1 \hbox{ in } (-\infty, 0], \quad -1 \leq \rho^{\prime}(\xi) \leq 0 \quad\hbox{and}\quad \left|\rho^{\prime \prime}(\xi)\right| \leq 1  \hbox{ in } \mathbb{R}.	
\end{equation}

\begin{lem}\label{fix T sup sub} 
Let $T>0$ be fixed and all the assumptions in Proposition {\rm \ref{stability}} hold. 
Let $\mu>0$ be a positive constant satisfying \eqref{set-mu}, and let $\psi_{i,T}^-(t)$ and $\psi_{i,T}^+(t)$ be, respectively, the unique solutions of \eqref{eigenfunction-0} and \eqref{eigenfunction-1} satisfying \eqref{normalize-supsub}.
Then, there exist $\varepsilon_0>0$ and $K>0$ such that for any $\xi_0\in\R$ and  $\varepsilon \in (0, \varepsilon_0]$, the pair of functions $\mathbf{v}^+(t,x)=(v_1^+,v_2^+)(t,x)$ $($resp. $\mathbf{v}^-(t,x)=(v_1^-,v_2^-)(t,x)$$)$ defined by 
\begin{equation}\label{define-subsuper-solu}
	\left\{\begin{array}{l}
		v_1^{\pm}(t, x)=\varphi_{1,T}\left( \xi_{\pm}(t,x),t \right) \pm \varepsilon \mathrm{e}^{-\mu t}\left(\rho(\xi_\pm(t,x)) \psi_{1,T}^{+}(t)+(1-\rho(\xi_\pm(t,x))) \psi_{1,T}^{-}(t)\right), \vspace{5pt} \\
		v_2^{\pm}(t, x)=\varphi_{2,T}\left(\xi_\pm(t,x),t\right) \pm \varepsilon \mathrm{e}^{-\mu t}\left(\rho(\xi_\pm(t,x)) \psi_{2,T}^{+}(t)+(1-\rho(\xi_\pm(t,x))) \psi_{2,T}^{-}(t)\right),
	\end{array}\right.
\end{equation}
 is a super-solution $($resp. sub-solution$)$ of \eqref{eq:change of variable} for $(t,x)\in (0,+\infty)\times \R$, where $\xi_{\pm}(t,x)=x-c_Tt+\xi_0 \pm \varepsilon K(\mathrm{e}^{-\mu t}-1)$, and $\rho \in C^2(\mathbb{R})$ is a nonnegative function satisfying \eqref{function-rho}. 
\end{lem}

\begin{proof}
	We only show that $\mathbf{v}^+(t,x)=(v_1^+,v_2^+)(t,x)$ is a super-solution, as the verification of the sub-solution is analogous. Since $T>0$ is fixed, to simplify the notations,  instead of $d_{i, T}(\cdot)$, $r_{i, T}(\cdot)$, $a_{i, T}(\cdot)$, $k_{i, T}(\cdot)$, $p_{i,T}(\cdot)$,  $\varphi_{i,T}(\cdot,\cdot)$, $\lambda_{i,T}^{\pm}$,  $\psi_{i,T}^{\pm}(\cdot)$ and
	$c_T$, we write $d_i(\cdot)$, $r_i(\cdot)$, $a_i(\cdot)$, $k_i(\cdot)$, $p_i(\cdot)$, $\varphi_{i}(\cdot,\cdot)$, $\lambda_{i}^{\pm}$,  $\psi_{i}^{\pm}(\cdot)$ and $c$, respectively. 
	
	Let $\sigma_0 \in(0,1 / 4)$ be a small constant such that
	\bqq\label{sigma_0}
	\sigma_0 \leq \min \left\{\frac{\mu}{4\theta_+\gamma_+},\, \frac{\mu}{4}\left\|\frac{k_1p_2 \psi_{2}^{-}}{\psi_{1}^{-}}\right\|^{-1},\, \frac{\mu}{4}\left\|\frac{k_2p_1 \psi_{1}^{+}}{\psi_{2}^{+}}\right\|^{-1}\right\},
	\eqq
	where $\theta_+$ and $\gamma_+$ are the positive constants provided by \eqref{bound-rak} and \eqref{bound-pT}, respectively. 
	By the definition of periodic traveling waves, there exists $M>0$ sufficiently large such that, for $i=1,2$, 	
	\bqq\label{wave fix T}
	\begin{aligned}
		&\begin{cases}0<\varphi_{i}(\xi,t)<\sigma_0, & \text { for all }\,  \xi \geq M,\,t\in\R, \vspace{5pt} \\
		1-\sigma_0<\varphi_{i}(\xi,t)<1, & \text { for all }\, \xi \leq-M, \,t\in\R, \vspace{5pt} \\
		 \sigma_0 / 2<\varphi_{i}(\xi,t)<1-\sigma_0 / 2 & \text { for all }\, -M<\xi<M,\,t\in\R.\end{cases}
	\end{aligned}
	\eqq
	Since $\varphi_{i}(\xi,t)$ is decreasing in $\xi$ and periodic in $t$, it is clear that 
	\begin{equation}\label{choose-beta0}
		\beta_0:=-\max _{-M\leq\xi\leq M,\,t\in\R,\,i=1,2}\varphi_{i}'(\xi,t) \in (0,+\infty).
	\end{equation} 
	We will show that $K>0$ and  $\varepsilon_0\in (0,1)$ satisfying 	
		\bqq\label{choose-K}
	K=\frac{2\left( \left(\mu+|c|+\theta_++4\theta_+\gamma_+ \right) \displaystyle\max_{i=1,2}   \left(\left\|\psi_{i}^{+}+\psi_{i}^{-}\right\|\right)+\displaystyle\max_{i=1,2} \left( \left\| \left(\psi_{i}^{+}\right)^{\prime}\right\| +\left\|\left(\psi_{i}^{-}\right)^{\prime} \right\| \right)\right) }{\mu\beta_0},
	\eqq
	and 
	\bqq\label{choose-epsilon0}
	0<\varepsilon_0 \leq \min \left\{\frac{\sigma_0}{2 \displaystyle\max_{i=1,2}\left\|\psi_{i}^{-}+\psi_{i}^{+}\right\|},\, \frac{\beta_0}{2\displaystyle\max_{i=1,2}\left\|\psi_{i}^{-}-\psi_{i}^{+}\right\|}\right\},
	\eqq
	are the desired constants of the present lemma. 
    To do so, according to Definition \ref{defi sup sub}, it suffices to show that, for any $\epsilon \in (0,\epsilon_0)$ and $\xi_0\in\R$, 
    $$
    \left\{\begin{array}{l}
    	\mathcal{L}_1 (\mathbf{v}^+):=\partial_t v^{+}_{1}(t,x)-d_{1}(t)\partial_{x x}v^{+}_{1}(t,x)-g_{1}\left(t, \mathbf{v}^+(t,x)\right) \geq 0, \vspace{5pt} \\
    	\mathcal{L}_2 (\mathbf{v}^+):=\partial_t v^{+}_{2}(t,x)-d_{2}(t)\partial_{x x} v^{+}_{2}(t,x)-g_{2}\left(t, \mathbf{v}^+(t,x)\right) \geq 0, 
    \end{array}\right.
    $$
    for $(t, x) \in(0, +\infty) \times \mathbb{R}$. We will only verity that $\mathcal{L}_1 (\mathbf{v}^+)\geq 0$, as the proof of $\mathcal{L}_2 (\mathbf{v}^+) \geq 0$ is similar.
    
    Let us denote 	
	$$ q(t)=\varepsilon e^{-\mu t}\quad \hbox{and}\quad \eta(t)=\varepsilon K(\mathrm{e}^{-\mu t}-1)\quad\hbox{for }\, t>0.$$
	Since $\boldsymbol{\varphi}(\xi,t):=(\varphi_{1},\varphi_{2})(\xi,t)$ is an entire solution of \eqref{eq:change-tw}, direct computations give that, for $(t, x) \in(0, +\infty) \times \mathbb{R}$
	\bqq\label{concrete expression fix T}
	\begin{aligned}
		\mathcal{L}_1 (\mathbf{v}^+)=\,\, & \eta^{\prime}\partial_{\xi}\varphi_{1}  +q^{\prime}\left(\rho \psi_{1}^{+}+(1-\rho) \psi_{1}^{-}\right) -\left(g_{1}\left(t, \mathbf{v}^+\right)-g_{1}\left(t, \boldsymbol{\varphi}\right)\right) \vspace{5pt} \\
		&\quad +q\left[(\eta^{\prime}-c)\rho'\left(\psi_{1}^{+}-\psi_{1}^{-}\right)-d_{1} \rho''\left(\psi_{1}^{+}-\psi_{1}^{-}\right)+\rho\left(\psi_{1}^{+}\right)^{ \prime}+(1-\rho)\left(\psi_{1}^{-}\right)^{\prime }\right],
	\end{aligned}
	\eqq
	where $\rho(\cdot)$, $\rho^{\prime}(\cdot)$, $\rho^{\prime \prime}(\cdot)$ are evaluated at $\xi_+(t,x)=x-ct+\xi_0 +\eta(t)$; and
	$\partial_{\xi}\varphi_{1}(\cdot,\cdot)$, $\boldsymbol{\varphi}(\cdot,\cdot)$ are evaluated at $(\xi_+,t)$. 
	 We will complete the proof of $\mathcal{L}_1 (\mathbf{v}^+) \geq 0$ by considering three cases: (a) $\xi_+(t,x)\geq M$; (b) $\xi_+(t,x) \leq -M$; $(c)  -M<\xi_+(t,x)<M$.
	
In the case where $(t,x)\in (0,+\infty)\times \R$ such that $\xi_+(t,x)\geq M$, 
 with our choice of $\rho$ satisfying \eqref{function-rho}, we have $\rho=\rho^{\prime}=\rho^{\prime \prime}=0$, and 
 $$g_{1}\left(t, \mathbf{v}^+\right)-g_{1}\left(t, \boldsymbol{\varphi}\right)
 	=a_{1}p_{1}q\psi_{1}^{-}  \left(1-(v^{+}_{1}+\varphi_{1})\right)-k_{1}p_{2} \left(q\psi_{1}^{-} +\varphi_{1}\varphi_{2}-v^{+}_{1}v^{+}_{2}\right).$$
Since $\eta(\cdot)$ is decreasing, and  $\varphi_{1}$ is decreasing in its first variable, it follows that 
	$$
		\mathcal{L}_1 (\mathbf{v}^+) \geq  q^{\prime} \psi_{1}^{-}+ q(\psi_{1}^{-}) ' -a_{1}p_{1}q\psi_{1}^{-}  \left(1-(v^{+}_{1}+\varphi_{1})\right)+k_{1}p_{2} \left(q\psi_{1}^{-}+\varphi_{1}\varphi_{2}-v^{+}_{1}v^{+}_{2}\right).
	$$
By the first line of \eqref{wave fix T} and \eqref{choose-epsilon0}, it is easily seen that  $q\psi_{1}^{-}\in(0,\sigma_0/2)$,  and  $\varphi_{1}+q \psi_{1}^{-}\in(0,2\sigma_0)$. Noticing that $\psi_{1}^{-}$ is a solution of \eqref{eigenfunction-0} with $i=1$,  we have
	$$
	\begin{aligned}
		\mathcal{L}_1 (\mathbf{v}^+)&\geq q^{\prime} \psi_{1}^{-}+ q(\psi_{1}^{-}) '-a_{1}p_{1}q\psi_{1}^{-}+k_{1}p_{2}q\psi_{1}^{-}+k_{1}p_{2} \left(\varphi_{1}\varphi_{2}-v^{+}_{1}v^{+}_{2}\right)
	\vspace{5pt} \\
		& \displaystyle = q^{\prime} \psi_{1}^{-}+\lambda_{1}^{-}q \psi_{1}^{-}  -k_{1}p_{2}\varphi_{2}q \psi_1^{-}-q\left(\varphi_{1}+q \psi_{1}^{-}\right) \frac{k_{1}p_{2} \psi_{2}^{-}}{\psi_{1}^{-}} \psi_{1}^{-} \vspace{5pt} \\
			& \displaystyle \geq  q^{\prime} \psi_{1}^{-}+\lambda_{1}^{-}q \psi_{1}^{-}  -\sigma_0k_{1}p_{2}q \psi_1^{-}-2\sigma_0  \frac{k_{1}p_{2} \psi_{2}^{-}}{\psi_{1}^{-}} q \psi_{1}^{-}. 
	\end{aligned}
	$$
	 Furthermore, it follows from \eqref{sigma_0} that $0<\sigma_0k_{1}p_{2} \leq \sigma_0 \theta_+\gamma_+ \leq \mu/4$ and  
	 $0<2\sigma_0  k_{1}p_{2} \psi_{2}^{-}/\psi_{1}^{-} \leq \mu/2 $. This together with the fact that $q(t)=\varepsilon e^{-\mu t}$ implies that
\begin{equation*}%\label{geqM}
		\mathcal{L}_1 (\mathbf{v}^+) \geq q^{\prime} \psi_{1}^{-}+\lambda_{1}^{-}q \psi_{1}^{-} -\frac{\mu}{4}q \psi_{1}^{-}-\frac{\mu}{2}q \psi_{1}^{-}= \left(\lambda_{1}^{-} -\frac{7}{4}\mu\right)q \psi_{1}^{-}.
\end{equation*}
	Finally, thanks to the restriction of $\mu$ in \eqref{set-mu}, we obtain $\mathcal{L}_1 (\mathbf{v}^+) \geq 0$ when $\xi_+(t,x)\geq M$.

	Next, we consider the case where $(t,x)\in (0,+\infty)\times \R$ such that $\xi_+ \leq -M$. Similarly as above,  we have $\rho=1, \rho^{\prime}=\rho^{\prime \prime}=0$, 
%	$$g_{1}\left(t, \mathbf{v}^+\right)-g_{1}\left(t, \boldsymbol{\varphi}\right)=a_{1}p_{1}q\psi_{1}^{+}  \left(1-(v^{+}_{1}+\varphi_{1})\right)-k_{1}p_{2} \left(q\psi_{1}^{+}+\varphi_{1}\varphi_{2}-v^{+}_{1}v^{+}_{2}\right),$$
	and 
    $$
		\mathcal{L}_1 (\mathbf{v}^+)
		\geq q^{\prime} \psi_{1}^{+}+ q(\psi_{1}^{+}) ' -a_{1}p_{1}q\psi_{1}^{+}  \left(1-(v^{+}_{1}+\varphi_{1})\right)+k_{1}p_{2} \left(q\psi_{1}^{+}+\varphi_{1}\varphi_{2}-v^{+}_{1}v^{+}_{2}\right).
	$$ 
	By the second line of  \eqref{wave fix T} and \eqref{choose-epsilon0},  there holds $q\psi_{1}^{+}\in(0,\sigma_0/2)$ and $v^{+}_{1}+\varphi_{1}=2\varphi_{1}+q\psi_{1}^{+}\in(2-2\sigma_0,2+\sigma_0)$. Since $\psi_{1}^{+}(\cdot)$ is a solution of  \eqref{eigenfunction-1} with $i=1$, and since $0<2\sigma_0a_1p_1\leq 2\sigma_0\theta_+\gamma_+\leq \mu/2$ (see \eqref{sigma_0}), it follows that 
	$$
	\begin{aligned}
		\mathcal{L}_1(\mathbf{v}^+) 
		&\geq  q^{\prime} \psi_{1}^{+}+ q(\psi_{1}^{+}) '+a_{1}p_{1}q\psi_{1}^{+}(1-2\sigma_0) +k_{1}p_{2}q\left(\psi_{1}^{+}-\varphi_{1}\psi_{2}^{+}-\varphi_{2}\psi_{1}^{+}-q\psi_{1}^{+}\psi_{2}^{+}\right) 
		\vspace{5pt} \\
		&\geq q^{\prime} \psi_{1}^{+}+\lambda_{1}^{+}q \psi_{1}^{+}-\frac{\mu}{2} q \psi_{1}^{+}  +k_{1}p_{2}(1-\varphi_{2})q \psi_{1}^{+}-\left(\varphi_{1}+ q\psi_{1}^{+}\right) \frac{k_{1}p_{2} \psi_{2}^{+}}{\psi_{1}^{+}} q\psi_{1}^{+} \vspace{5pt} \\
		&\geq q^{\prime} \psi_{1}^{+}+\lambda_{1}^{+}q \psi_{1}^{+}-\frac{\mu}{2} q \psi_{1}^{+}  -\left(\varphi_{1}+q \psi_{1}^{+}\right) \frac{k_{1}p_{2} \psi_{2}^{+}}{\psi_{1}^{+}} q\psi_{1}^{+}.
	\end{aligned}
	$$
	Notice that $1-\sigma_0\leq \varphi_{1}+q \psi_{1}^{+} \leq 1+\sigma_0 \leq 3/2$. Then, by using the third term of the normalization \eqref{normalize-supsub}, we have
	\begin{equation*}%\label{<-M fix T}
		\mathcal{L}_1 (\mathbf{v}^+) \geq q^{\prime} \psi_{1}^{+}+\lambda_{1}^{+}q \psi_{1}^{+} -\frac{\mu}{2}q \psi_{1}^{+}-\frac{3\mu}{8}q \psi_{1}^{+}= \left(\lambda_{1}^{+}-\frac{15}{8}\mu\right) q \psi_{1}^{+}.
	\end{equation*}
	Therefore, since $0<\mu \leq \lambda_{1}^{+}/2$ by \eqref{set-mu}, we obtain $\mathcal{L}_1 (\mathbf{v}^+) \geq 0$ when $\xi_+\leq -M$.
	
	It remains to consider the values $(t,x)\in (0,+\infty)\times \R$ such that $-M<\xi_+(t,x)<M$. In this case, we have $\varphi_{1},\varphi_{2} \in\left[\sigma_0 / 2,1-\sigma_0 / 2\right]$, whence	
	by \eqref{choose-epsilon0}, $ v^{+}_{1}, v^{+}_{2} \in\left[\sigma_0/2 ,1 \right]$. 
	Since the function $\rho \in C^2(\mathbb{R})$ satisfying \eqref{function-rho}, 
	and since $|q(\psi_{1}^{+}-\psi_{1}^{-})|\leq \beta_0/2$ (due to \eqref{choose-epsilon0}), it is straightforward to check that
	$$
	\begin{aligned}
		&q\rho'\left(\psi_{1}^{+}-\psi_{1}^{-}\right)\eta'+q\left[-(d_{1} \rho''+c)\left(\psi_{1}^{+}-\psi_{1}^{-}\right)+\left(\rho\left(\psi_{1}^{+}\right)^{ \prime}+(1-\rho)\left(\psi_{1}^{-}\right)^{\prime }\right)\right]\vspace{5pt} \\
		\geq \,\, & \frac{\beta_0}{2}\eta'-\left( \left(\theta_++|c|\right) \left\|\psi_{1}^{+}+\psi_{1}^{-}\right\|+ \left\| \left(\psi_{i}^{+}\right)^{\prime}\right\| +\left\|\left(\psi_{i}^{-}\right)^{\prime} \right\| \right)q,
	\end{aligned}
	$$
	and that
$$\left.\begin{array}{ll}
		&g_{1}\left(t, \boldsymbol{\varphi}\right)- g_{1}\left(t, \mathbf{v}^+\right)  \vspace{5pt} \\
		=& a_{1}p_{1}\left( \varphi_{1}-v^{+}_{1}\right)  \left(1-(v^{+}_{1}+\varphi_{1})\right)-k_{1}p_{2} \left(\varphi_{1}(1-\varphi_{2})-v^{+}_{1}(1-v^{+}_{2})\right)  \vspace{5pt} \\
		\geq & -a_{1}p_{1} \left\|1-(v^{+}_{1}+\varphi_{1})\right\|\left\|\psi_{1}^{+}+\psi_{1}^{-}\right\|q   -k_{1}p_{2}v^{+} \left\|\psi_{2}^{+}+\psi_{2}^{-}\right\|q   \vspace{5pt} \\
		\geq& - 3\theta_+\gamma_+  \left\|\psi_{1}^{+}+\psi_{1}^{-}\right\|q-\theta_+\gamma_+  \left\|\psi_{2}^{+}+\psi_{2}^{-}\right\|q   \vspace{5pt} \\
		\geq& -4\theta_+\gamma_+ \max_{i=1,2}\left\|\psi_{i}^{+}+\psi_{i}^{-}\right\| q.
\end{array}\right.
$$
Moreover, by \eqref{choose-beta0} and the monotonicity of $\eta(\cdot)$, we have $\eta^{\prime}\partial_{\xi}\varphi_{1} \geq -\beta_0\eta'$. 
It then follows from \eqref{concrete expression fix T} and our choice of the positive constant $K$  that
$$\left.\begin{array}{ll}
	 \mathcal{L}_1 (\mathbf{v}^+)\!\! & \displaystyle \geq -\beta_0 \eta^{\prime} +\frac{\beta_0}{2} \eta^{\prime} +\left( \left\|\psi_{1}^{+}+\psi_{1}^{-}\right\|\right)  q^{\prime} \vspace{5pt} \\
	& \qquad\qquad  \displaystyle   -\left(\left(\theta_++|c| +4\theta_+\gamma_+\right)\max_{i=1,2}\left\|\psi_{i}^{+}+\psi_{i}^{-}\right\| + \max_{i=1,2} \left( \left\| \left(\psi_{i}^{+}\right)^{\prime}\right\| +\left\|\left(\psi_{i}^{-}\right)^{\prime} \right\| \right)  \right) q \vspace{5pt} \\
	& \geq  \displaystyle -\frac{\beta_0}{2} \eta^{\prime} - \frac{\beta_0}{2}\mu K q.
\end{array}\right.
$$
Since $\eta(t)=\varepsilon K(\mathrm{e}^{-\mu t}-1)$, we thus obtain that $\mathcal{L}_1 (\mathbf{v}^+)\geq 0$ when $-M<\xi_+(t,x)<M$. 

Finally, combining the above, we can conclude that $\mathcal{L}_1 (\mathbf{v}^+)\geq 0$ for $(t, x) \in(0, +\infty) \times \mathbb{R}$. Proceeding similarly as above, there holds  $\mathcal{L}_2 (\mathbf{v}^+)\geq 0$ for $(t, x) \in(0, +\infty) \times \mathbb{R}$. Therefore, $\left(v^{+}_{1}, v^{+}_{2}\right)$ is a supsolution of \eqref{eq:change of variable}. This completes the proof of Lemma~\ref{fix T sup sub}.
\end{proof}

Having in hand Lemma \ref{fix T sup sub}, the proof of Proposition \ref{stability} follows the same lines as those used in \cite[Theorem 4.4]{B}; therefore, we omit the details. 

%%%%%%%%%%%%%%%%%%%%%%%%%%%%%%%%%%%%%%%%%%%%%%%%%%%%%%%%%%%%%%%%%%%%%%%%%%%%%%%%%%%%%%%%%%%%%%%
%%%%%%%%%%%%%%%%%%%%%%%%%%%%%%%%%%%%%%%%%%%%%%%%%%%%%%%%%%%%%%%%%%%%%%%%%%%%%%%%%%%%%%%%%%%%%%%%

\section{Proof of Theorem \ref{Theo small T}}
In this section, we first establish the existence of periodic traveling waves for the transformed system \eqref{eq:change of variable} with small $T$, and then characterize the homogenization speed limit as $T\to 0^+$.

\subsection{Existence of periodic traveling waves when $T$ is small}

This subsection is devoted to proving the following existence result by applying Proposition \ref{existence-fix}.

\begin{theo}\label{exitstence wave small}
	Let  {\rm (A1)} hold. There exists $T_*>0$ such that, for any $0<T<T_*$, system \eqref{eq:change of variable} admits a periodic traveling wave $\left(v_{1,T}, v_{2,T}\right)(t,x)=\left(\varphi_{1,T}, \varphi_{2,T}\right)(x-c_Tt, t)$ connecting  $\mathbf{0}$ and $\mathbf{1}$.
\end{theo} 

The proof of Theorem \eqref{exitstence wave small} is adapted from \cite{dhy} concerning spatially periodic traveling waves. 
For clarity, we divide the proof of Theorem \ref{exitstence wave small} into several lemmas. Recall that $p_{i,T}(\cdot)$ is the unique positive solution of equation \eqref{T-periodic equilibria}, and that $\bar{p}_i=\bar{r}_i/\bar{a}_i$. We first show that the constant  
$\bar{p}_i$ is the homogenization limit of $p_{i,T}(T\cdot)$.  

\begin{lem}\label{converge-pi-0}
For each $i=1,2$, $p_{i,T}(Ts)$ converges to $\bar{p}_i$ as $T\to 0^+$ uniformly in $s\in\R$.
\end{lem}

\begin{proof}
For each $T>0$, define $q_{i,T}(s)=p_{i,T}(Ts)$ for $s\in\R$. It is clear that 
$q_{i,T}(\cdot)$ is positive, $1$-periodic, and satisfies
\bqq \label{p_2,n}
\frac{d}{ds}q_{i,T}(s)=Tq_{i,T}(s)\left(r_i(s)-a_i(s)q_{i,T}(s)\right)\quad \hbox{for }\, s\in\R.
\eqq
Since $q_{i,T}(\cdot)$ is bounded in $\R$ uniformly with respect to $T>0$ by \eqref{bound-pT}, it follows that $q_{i,T}'(s)\to 0$ as $T\to 0^+$ uniformly in $s\in\R$. 
Then by the periodicity, one finds some constant $q_{i,0}>0$ such that 
$q_{i,T}(s) \to q_{i,0}$ as $T\to 0^+$ uniformly in $s\in\R$. 
On the other hand, for each $T>0$, multiplying both sides of \eqref{p_2,n} by $T^{-1}q_{i,T}^{-1}$ and integrating over $[0, 1]$, we have
$0=\int_0^1 \left(r_i(s)-a_i(s)q_{i,T}(s)\right)ds$. 
Passing to the limit as $T \rightarrow 0^+$ gives $q_{i,0}=\bar{r}_i/\bar{a}_i=\bar{p}_i$.  
\end{proof}

The following lemma provides the strong stability of the fixed points $\mathbf{0}$ and $\mathbf{1}$ when $T$ is small.  

\begin{lem}\label{linearly stable claim}
Let  {\rm (A1)} hold. There exists $T_1>0$ such that \eqref{0-1-stable} holds for all $0<T < T_1$, and that for the map $Q_T: [0,1]^2\to[0,1]^2$, the fixed points $\mathbf{0}$ and $\mathbf{1}$ are, respectively, strongly stable from above and below, uniformly with respect to $T\in (0,T_1)$ in the sense of Definition {\rm \ref{defi-stability}}.  
\end{lem}

\begin{rmk}\label{rem-unif-stable}
The uniform stability of $\mathbf{0}$ and $\mathbf{1}$ in $T\in (0,T_1)$ implies that, there exist a small constant $\epsilon_0>0$ and two vectors $\mathbf{e}_{\pm} \in (0,+\infty)\times (0,+\infty)$ such that for any $T\in (0,T_1)$, the map $Q_T: [0,1]^2\to[0,1]^2$ admits no fixed points in the sets $\{\mathbf{x}\in\R^2: \epsilon_0 \mathbf{e}_{-} \geq  \mathbf{x} > \mathbf{0} \}$ and $\{\mathbf{x}\in\R^2: \mathbf{1} >  \mathbf{x} \geq \mathbf{1}- \epsilon_0 \mathbf{e}_{+} \}$. This property will be useful in showing the instability of intermediate fixed points in Lemma {\rm \ref{linearly unstable claim}} below.
\end{rmk}

\begin{proof}[Proof of Lemma {\rm \ref{linearly stable claim}}]
Recall that $\lambda_{i,T}^{-}$ and $\lambda_{i,T}^{+}$ are the real numbers defined in \eqref{seperate-lambda}.
Since
\begin{equation*}
	\lambda^-_{1, T}=-\int_0^1 \partial_{v_1}g_{1,T}(sT,0,0) ds=-\int_0^1 \left( a_{1}(s)p_{1,T}(sT)-k_{1}(s)p_{2,T}(sT)\right)ds, 
\end{equation*}
by using Lemma \ref{converge-pi-0}, we obtain 
$\lambda^-_{1, T} \to \bar{k}_1\bar{p}_2-\bar{r}_1$ as $T\to 0^+$. 
In a similar way, one can conclude that as $T\to 0^+$,  
$\lambda^+_{1, T} \to  \bar{r}_1$, $\lambda^-_{2, T}\to \bar{r}_2$ and $\lambda^+_{2, T} \to \bar{k}_2\bar{p}_1-\bar{r}_2$. 
By assumption (A1), all the limits are positive constants. Then, one finds some small $T_1>0$ such that
\begin{equation}\label{choose-sigma0}
	\min_{i=1,2} \left\{\lambda^-_{i,T}, \lambda^+_{i,T}\right\}> \sigma_0:=\frac{1}{2} \min \left\{\bar{r}_1, \bar{r}_2, \bar{k}_1\bar{p}_2-\bar{r}_2, \bar{k}_2\bar{p}_1-\bar{r}_1 \right\}\quad\hbox{for all }\,0<T<T_1. 
\end{equation} 
This in particular implies that \eqref{0-1-stable} holds for all $0<T<T_1$, whence by Lemma \ref{lem-stable-0-1}, the fixed points  $\mathbf{0}$ and $\mathbf{1}$ are, respectively, strongly stable from above and below. 

It remains to show that the strong stability holds uniformly in $T \in (0,T_1)$. We only show that $\mathbf{0}$ is strongly stable from above uniformly in $T$, as the analysis for the uniform stability of 
$\mathbf{1}$ follows analogously.
For each $T>0$, let $\psi_{1,T}^{-}$ and $\psi_{2,T}^{-}$ be the unique solutions of the linear problem \eqref{eigenfunction-0} under the following normalization conditions 
\bqq\label{normalization small T stable}
\max _{t \in \mathbb{R}}\psi_{1,T}^{-}(t)=1 \quad \text { and } \quad \min _{t \in \mathbb{R}} \psi_{2,T}^{-}(t)=\frac{2 \gamma_+ \theta_+ }{\sigma_0 },
\eqq
where $\sigma_0$, $\theta_+$ and $\gamma_+$  are, respectively, the positive constants given in \eqref{choose-sigma0}, \eqref{bound-rak} and \eqref{bound-pT}.  
Since $p_{2,T}$ and $\lambda^-_{2,T}$ are bounded in $T>0$ (due to \eqref{bound-pT}), it follows from \eqref{explicit-eigenfunction} that there exists $\nu_0>0$ such that 
\begin{equation}\label{bound-psi2-small}
\psi_{2,T}^{-}(t)\leq \nu_0 \quad \hbox{for all }\,  t\in\R,\, T \in (0,T_1).
\end{equation}
Now, we choose a small constant $\epsilon_0$ (independent of $T$) such that
	\bqq \label{epsilon_0}
	0<\epsilon_0 \leq   \frac{\sigma_0 }{2\nu_0 \theta_+ \gamma_+}.
	\eqq
Letting $\epsilon \in (0, \epsilon_0]$ and $T \in (0,T_1)$ be arbitrary, we verify that  $(\epsilon\psi_{1,T}^{-},\epsilon\psi_{2,T}^{-})(\cdot)$ is a strict super-solution of \eqref{eq:change of variable} in the sense of Definition \ref{defi sup sub}. Indeed, for any $t\geq  0$, by using \eqref{bound-psi2-small} and \eqref{epsilon_0}, one computes that 
$$
\begin{aligned}
	& \left(\epsilon\psi_{1,T}^{-}\right)' - g_{1,T}\left( t,\epsilon\psi_{1,T}^{-},\epsilon\psi_{2,T}^{-}\right)  \vspace{5pt} \\
	> \, &   \left(\epsilon\psi_{1,T}^{-}\right)'-\left(a_{1,T}(t)p_{1,T}(t)-k_{1,T}(t)p_{2,T}(t)\right)\left(\epsilon\psi_{1,T}^{-}\right)-\left(k_{1,T}(t)p_{2,T}(t)\left(\epsilon\psi_{2,T}^{-}\right)\right)\left(\epsilon\psi_{1,T}^{-}\right)\vspace{5pt} \\
	\geq\,&  \left(\epsilon\psi_{1,T}^{-}\right)'- \partial_{v_1}g_{1,T}(t,0,0)\left(\epsilon\psi_{1,T}^{-}\right)-\sigma_0\left(\epsilon\psi_{1,T}^{-}\right)\vspace{5pt} \\
	=\,& \left(\lambda_{1,T}^{-}-\sigma_0\right)\left(\epsilon\psi_{1,T}^{-}\right).
\end{aligned}
$$
On the other hand, it follows from \eqref{normalization small T stable} and \eqref{epsilon_0} that for $t \geq 0$,
$$
\begin{aligned}
&\left(\epsilon\psi_{2,T}^{-}\right)'- g_{2,T}\left( t,\epsilon\psi_{1,T}^{-},\epsilon\psi_{2,T}^{-}\right)\vspace{5pt} \\
>\,	& \left(\epsilon\psi_{2,T}^{-}\right)' +a_{2,T}(t)p_{2,T}(t)\left(1-\left(\epsilon\psi_{2,T}^{-}\right)\right)\left(\epsilon\psi_{2,T}^{-}\right)-\left(k_{2,T}(t)p_{1,T}(t)\frac{\psi_{1,T}^{-}}{\psi_{2,T}^{-}}\right)\left(\epsilon\psi_{2,T}^{-}\right)\vspace{5pt} \\
\geq\,	& \left(\epsilon\psi_{2,T}^{-}\right)'
+a_{2,T}(t)p_{2,T}(t)\left(\epsilon\psi_{2,T}^{-}\right)-\frac{\sigma_0}{2}\left(\epsilon\psi_{2,T}^{-}\right)-\frac{\sigma_0}{2}\left(\epsilon\psi_{2,T}^{-}\right)\vspace{5pt} \\
=\,& \left(\lambda_{2,T}^{-}-\sigma_0\right)\epsilon\psi_{2,T}^{-}(t).
\end{aligned}
$$
Combining the above, one sees from the choice of $\sigma_0$ that for $t\geq 0$, 
$$\left(\epsilon\psi_{1,T}^{-}\right)' - g_{1,T}\left( t,\epsilon\psi_{1,T}^{-},\epsilon\psi_{2,T}^{-}\right) > 0 \quad\hbox{and} \quad \left(\epsilon\psi_{2,T}^{-}\right)'- g_{2,T}\left( t,\epsilon\psi_{1,T}^{-},\epsilon\psi_{2,T}^{-}\right) > 0.$$
Therefore, $(\epsilon\psi_{1,T}^{-},\epsilon\psi_{2,T}^{-})(\cdot)$ is a strict super-solution of \eqref{eq:change of variable}.
Finally, by using the comparison principle (see Lemma \ref{comparison principle}), one can conclude  that
$$ \left(\epsilon\psi_{1,T}^{-},\epsilon\psi_{2,T}^{-}\right)(t) \gg   (v_1,v_2)(t) \,\,\hbox{ for all }\, t>0, $$
where $(v_1,v_2)(\cdot)$ is the spatially homogeneous solution of \eqref{eq:change of variable} with initial data $(\epsilon\psi_{1,T}^{-},\epsilon\psi_{2,T}^{-})(0)$. 
Due to the periodicity of  $(\epsilon\psi_{1,T}^{-},\epsilon\psi_{2,T}^{-})(\cdot)$, this in particular implies that 
$(\epsilon\psi_{1,T}^{-},\epsilon\psi_{2,T}^{-})(0) \gg  (v_1,v_2)(T)$. 
Setting $\mathbf{e}=(\psi_{1,T}^{-},\psi_{2,T}^{-})(0)$, one then reaches that 
$ \epsilon \mathbf{e} \gg Q_T[\epsilon \mathbf{e}]$ for all $\epsilon \in (0,\epsilon_0]$ and $T\in (0,T_1)$. Namely, $\mathbf{0}$ is strongly stable from above uniformly in  $T\in (0,T_1)$. This ends the proof of Lemma \ref{linearly stable claim}.
\end{proof}

Next, we verify that all the intermediate fixed points strictly between $\mathbf{0}$ and $\mathbf{1}$ are linearly unstable when $T$ is small.

\begin{lem}\label{linearly unstable claim}
There exists $T_2>0$ such that $\lambda_1(\bar{\mathbf{v}}, T)<0$ for every $0<T<T_2$ and for every $T$-periodic solution $\bar{\mathbf{v}}$ of \eqref{eq:variable-kinetic} with $(0,0) \ll \bar{\mathbf{v}} \ll(1,1)$, where $\lambda_1(\bar{\mathbf{v}}, T)$ is the principal Floquet multiplier for problem \eqref{inter-eigen}. 
\end{lem}

\begin{proof}
Assume by contradiction that there are some sequences $\left(T_n\right)_{n \in \mathbb{N}}$ in $(0,+\infty)$, $\left(\bar{\mathbf{v}}_n\right)_{n \in \mathbb{N}}$ and $\left(\boldsymbol \psi_n\right)_{n \in \mathbb{N}}$ in $C^1(\R;\R^2)$ such that $T_n \rightarrow 0^{+}$ as $n \rightarrow+\infty$ and, for each $n \in \mathbb{N}$, $\bar{\mathbf{v}}_n=\left(\bar{v}_{1, n}, \bar{v}_{2, n}\right)$ satisfies
\bqq\label{unstable  equilibria}
\left\{\begin{array}{ll}
\displaystyle \frac{d}{dt} \bar{v}_{1, n}=g_{1,T_n}\left(t, \bar{\mathbf{v}}_n\right) & \text {in } \mathbb{R}, \vspace{5pt} \\ 
\displaystyle \frac{d}{dt}\bar{v}_{2, n}= g_{2,T_n}\left(t, \bar{\mathbf{v}}_n\right) & \text {in } \mathbb{R}, \vspace{5pt} \\
  (0,0) \ll \bar{\mathbf{v}}_n(\cdot) \ll(1,1) & \text {in } \mathbb{R},\vspace{5pt} \\
  \bar{\mathbf{v}}_n(\cdot)= \bar{\mathbf{v}}_n(T_n+\cdot)  & \text {in } \mathbb{R},
\end{array}\right.
\eqq
and $\boldsymbol \psi_n=\left(\psi_{1, n}, \psi_{2, n}\right)$ satisfies
\bqq\label{inter-eigen-n}
\left\{\begin{array}{ll}
	\displaystyle \frac{d \boldsymbol \psi_n}{dt}- A_{T_n}(\bar{\mathbf{v}}_n) \boldsymbol \psi_n= \lambda_1(\bar{\mathbf{v}}_n, T_n) \boldsymbol \psi_n & \text {in } \mathbb{R}, \vspace{5pt} \\
	\boldsymbol \psi_n (\cdot+T_n)=\boldsymbol \psi_n(\cdot) & \text {in } \mathbb{R}, \vspace{5pt} \\
	\boldsymbol \psi_n(\cdot)\gg \mathbf{0} & \text {in } \mathbb{R},
\end{array}\right.
\eqq
with $\lambda_1(\bar{\mathbf{v}}_n, T_n) \geq 0$, where $A_{T_n}(\bar{\mathbf{v}}_n)$ is the 
matrix function defined as in \eqref{matrxi-A}. Since $\boldsymbol \psi_n$ is unique up to multiplication, we normalize it as follows 
\bqq\label{normalize unstable eigenfunction}
\max _{t \in \mathbb{R}} \left( \psi_{1, n}(t)+ \psi_{2, n}(t)\right)=1.
\eqq
Notice that, thanks to \eqref{bound-pT}, all the entries in $A_{T_n}(\bar{\mathbf{v}}_n)$ are bounded functions in $\R$ uniformly with respect to $n\in\N$. One infers that the sequence $\left(\lambda_1\left(\bar{\mathbf{v}}_n, T_n\right)\right)_{n \in \mathbb{N}}$ is bounded. Then, up to extraction of some subsequence, there is a real number $\lambda_\infty \geq 0$ such that 
$\lambda_1\left(\bar{\mathbf{v}}_n, T_n\right) \rightarrow \lambda_\infty$  as $n \rightarrow+\infty$.  

For each $n\in\N$, denote  $p_{i,n}(s):=p_{i,T_n}(T_ns)$,
$$
\boldsymbol{\alpha}_n(s)=(\alpha_{1,n},\alpha_{2,n})(s):=\bar{\mathbf{v}}_n(T_n s) \quad \text {and} \quad \boldsymbol{\beta}_n(s)=(\beta_{1,n},\beta_{2,n})(s):=\boldsymbol \psi_n\left(T_n s\right)
$$
for $s\in\R$. Clearly, $p_{i,n}(s)$, $\boldsymbol{\alpha}_n(s)$ and $\boldsymbol{\beta}_n(s)$ are 1-periodic functions. 
Proceeding similarly as in the proof of Lemma \ref{converge-pi-0}, one can conclude that as $n\to+\infty$, the following convergences hold in $C^1(\R)$: 
$p_{i,n}(\cdot) \to \bar{p}_i$, $\boldsymbol{\alpha}_n(\cdot) \rightarrow \boldsymbol{\alpha}_{\infty}$ and $\boldsymbol{\beta}_n(\cdot) \rightarrow \boldsymbol{\beta}_{\infty}$,
where $\boldsymbol{\alpha}_{\infty}=(\alpha_{1, \infty},\alpha_{2, \infty}) \in[0,1]^2$ and $\boldsymbol{\beta}_{\infty}=(\beta_{1, \infty},\beta_{2, \infty}) \in \mathbb{R}_{+}^2$ are two constant vectors.
It is easily seen from \eqref{normalize unstable eigenfunction} that $\beta_{1, \infty}+\beta_{2, \infty}=1$. Furthermore, integrating the equations of $\bar{v}_{1,2}$ and $\bar{v}_{2,n}$ in \eqref{unstable  equilibria} over  $[0,T_n]$, and passing to the limits as $n \rightarrow \infty$, we obtain
\begin{equation*}
\left\{\begin{array}{l}
\bar{r}_1(1-\alpha_{1, \infty})\alpha_{1, \infty}-\bar{k}_1\bar{p}_2\left(1-\alpha_{2, \infty}\right) \alpha_{1, \infty}=0, \vspace{5pt} \\
\bar{r}_2(1-\alpha_{2, \infty})\alpha_{2, \infty}-\bar{k}_2\bar{p}_1 \left(1-\alpha_{2, \infty}\right)\alpha_{1, \infty}=0.
\end{array}\right.
\end{equation*}
Then, according to condition (A1), $\boldsymbol{\alpha}_{\infty}$ must be one of the constant vectors $\mathbf{0}$, $\mathbf{1}$, $\mathbf{e}_0=(0,1)$ and  
$$ \mathbf{e}_*:=(v_1^*,v_2^*)= \left(\frac{\bar{a}_1\bar{a}_2-\bar{a}_1\bar{k}_1\bar{r}_2/\bar{r}_1}{\bar{a}_1\bar{a}_2-\bar{k}_1\bar{k}_2},\, \frac{\bar{a}_2\bar{k}_2\bar{r}_1/\bar{r}_2-\bar{k}_1\bar{k}_2}{\bar{a}_1\bar{a}_2-\bar{k}_1\bar{k}_2}\right) \in(0,1)^2.$$

By Remark \ref{rem-unif-stable} and the fact that $\bar{\mathbf{v}}_n(T_n s)$ converges uniformly to $\boldsymbol{\alpha}_{\infty}$ in $s\in\R$ as $n\to+\infty$, one observes that $\boldsymbol{\alpha}_{\infty} \neq \mathbf{0}$ and $\boldsymbol{\alpha}_{\infty} \neq \mathbf{1}$. Next, we exclude the case that $\boldsymbol{\alpha}_{\infty} = \mathbf{e}_0$. Assume by contradiction that $\boldsymbol{\alpha}_{\infty}=\mathbf{e}_0$. Since $\boldsymbol{\beta}_{\infty}=\left(\beta_{1, \infty},\beta_{2, \infty}\right)$ is nonnegative and satisfies $\beta_{1, \infty}+\beta_{2, \infty}=1$, it is easily seen that either $\beta_{1, \infty}>0$ or $\beta_{2, \infty}>0$. If $\beta_{1, \infty}>0$, then integrating the equation of $\beta_{1, n}$ over $[0,1]$ gives
$$
\begin{aligned}
	& -\int_0^1\left[a_1(s)p_{1,n}(s)(1-2\alpha_{1, n}(s))-k_1(s) p_{2, n}(s)\left(1-\alpha_{2, n}(s)\right)\right]\beta_{1, n} ds \vspace{5pt} \\
	& \qquad\qquad\qquad-\int_0^1 k_1(s)p_{2, n}(s) \alpha_{1, n}(s) \beta_{2, n}(s) ds=\lambda_1\left(\bar{\mathbf{v}}_n, T_n\right) \int_0^1 \beta_{1, n}(s) ds.
\end{aligned}
$$
Passing to the limits as $n \rightarrow \infty$, we obtain $\lambda_\infty=-\int_0^1 a_1(s)\bar{p}_1 ds=-\bar{r}_1<0$,
which is a contradiction with $\lambda_\infty \geq 0$. Similarly, if $\beta_{2, \infty}>0$, one can find a contradiction by integrating the equation of $\beta_{2, n}$ over $[0,1]$. Thus, $\boldsymbol{\alpha}_{\infty}$ cannot be $\mathbf{e}_0$.

As a consequence, we have $\boldsymbol{\alpha}_{\infty}=\mathbf{e}_*$. In this situation, there must hold $\beta_{1, \infty}>0$ and $\beta_{2, \infty}>0$. Otherwise, we may assume without loss of generality that $\beta_{2, \infty}=0$ and $\beta_{1, \infty}=1$. Integrating the equation of $\beta_{2, n}$ over $[ 0,1]$ gives
$$
\begin{aligned}
	& \lambda_1\left(\bar{\mathbf{v}}_n, T_n\right) \int_0^1 \beta_{2, n}(s) ds=-\int_0^1 k_2(s)p_{1, n}(s)\left(1-\alpha_{2, n}(s)\right) \beta_{1, n}(s) ds \vspace{5pt} \\
	&\qquad \qquad \qquad +\int_0^1\left[a_2(s)p_{2,n}(s)(1-2\alpha_{2, n}(s))+k_2(s)p_{1, n}(s) \alpha_{1, n}(s)\right] \beta_{2, n}(s) ds.
\end{aligned}
$$
Sending to the limits as $n \rightarrow \infty$ yields $\bar{p}_1\left(1-\alpha_{2, \infty}\right) \bar{k}_2=0$, which is impossible. Integrating the equation of $\beta_{1, n}$ over $[0,1]$, one can derive a similar contradiction in the case where $\beta_{1, \infty}=0$ and $\beta_{2, \infty}=1$. Therefore, both $\beta_{1, \infty}$ and $\beta_{2, \infty}$ must be positive.

Now, multiplying both sides of the equation of $\beta_{1, n}$ by $\beta_{1, n}^{-1}$ and integrating over $[0,1]$ gives
$$
\begin{aligned}
	& -\int_0^1 \left[  a_1(s)p_{1,n}(s)(1-2\alpha_{1, n}(s))-k_1(s) p_{2, n}(s)\left(1-\alpha_{2, n}(s)\right)\right] ds \vspace{5pt} \\
	& \qquad\qquad\qquad-\int_0^1 k_1(s)p_{2, n}(s) \alpha_{1, n}(s)\frac{\beta_{2, n}(s)}{\beta_{1, n}(s)} ds=\lambda_1\left(\bar{\mathbf{v}}_n, T_n\right).
\end{aligned}
$$
Taking the limits as $n \rightarrow +\infty$, we have $\lambda_\infty =\bar{r}_1v_1^*-\bar{k}_1 \bar{p}_2v_1^*\beta_{2, \infty}\beta_{1, \infty}^{-1}$.
Similarly, multiplying both sides of the equation of $\beta_{2, n}$ by $\beta_{2, n}^{-1}$ and passing to the limits as $n \rightarrow +\infty$ yields that 
$ \lambda_\infty=\bar{r}_2(1-v_2^*)-\bar{k}_2 \bar{p}_1(1-v_2^*)\beta_{1, \infty}\beta_{2,\infty}^{-1}$.
Combining the above, we obtain that $B_0 \boldsymbol{\beta}_{\infty} = -\lambda_\infty \boldsymbol{\beta}_{\infty}$, where $B_0$ is a constant matrix given by
$$
B_0=\left(\begin{array}{cc}
	-\bar{r}_1v_1^* & \bar{k}_1 \bar{p}_2v_1^* \vspace{5pt} \\
	\bar{k}_2 \bar{p}_1(1-v_2^*) & -\bar{r}_2(1-v_2^*)
\end{array}\right).
$$
This means that $-\lambda_\infty$ is an eigenvalue of $B_0$, and $\boldsymbol{\beta}_{\infty}\gg \mathbf{0}$ is the corresponding eigenvector. 
Notice that the off-diagonal entries of $B_0$ are negative, and that by (A1),  $$\operatorname{det}\left(B_0\right)=\frac{\left(\bar{r}_1 \bar{a}_2-\bar{r}_2 \bar{k}_1\right)\left(\bar{r}_2 \bar{a}_1-\bar{r}_1 \bar{k}_2\right)}{\bar{a}_1 \bar{a}_2-\bar{k}_1 \bar{k}_2}<0.$$
It is then known from the Perron-Frobenius theorem that the dominant eigenvalue of $B_0$, denoted by $\lambda^+(B_0)$,  is positive, and it is the unique eigenvalue with an associated strictly positive eigenvector. This immediately gives $\lambda^+(B_0)=-\lambda_\infty >0$, which is a contradiction with $\lambda_\infty \geq 0$ again. The proof of Lemma \ref{linearly unstable claim} is thus complete.
\end{proof}

Finally, we take $T_*=\min\{T_1,T_2\}$. Based on Lemmas \ref{linearly stable claim}-\ref{linearly unstable claim}, 
Theorem \ref{exitstence wave small} follows directly from Proposition \ref{existence-fix}.

%%%%%%%%%%%%%%%%%%%%%%%%%%%%%%%%%%%%%%%%%%%%%%%%%%%%%%%%%%%%%%%%%%%%%%%%%%%%%%%%%%%%%%%%%%%%%%%%%%%%%%%%%
\subsection{Convergence of $c_T$ as $T\rightarrow 0^+$}

Let $T_*>0$ be the constant provided by Theorem \ref{exitstence wave small}, and for each $0<T<T_*$, let $(\varphi_{1,T},\varphi_{2,T})( x-c_Tt,t)$ be a periodic traveling wave of \eqref{eq:change of variable} connecting $\mathbf{0}$ and $\mathbf{1}$. It follows directly from Proposition \ref{stability} that $c_T \in \mathbb{R}$ is the unique wave speed and that $(\varphi_{1,T},\varphi_{2,T})$ is unique up to shift in the first variable. For definiteness and later use, we normalize this wave as follows 
\bqq\label{normalize small T wave}
\varphi_{1,T}( 0,0)=\frac{1}{2}\quad\hbox{for each } \,0<T<T_*.
\eqq
In this subsection, we prove the convergence of $c_T$ as $T\to 0^+$, and estimate the convergence rate as well. For clarity, we restate the result as follows.

\begin{theo}\label{converge-samllT}
	Let  {\rm (A1)} hold. Then the estimate \eqref{small speed cover} holds true.
\end{theo}

Theorem \ref{converge-samllT} is proved by a sub- and super-solution method, which relies on suitably perturbing the traveling wave for the homogenized equation \eqref{eq:homogenization limit} connecting $(0,\bar{p}_2)$ and $(\bar{p}_1,0)$. Recall that $(\phi_{1,0}, \phi_{2,0})(x-c_0 t)$ is such a wave. Taking the change of variable 
$\boldsymbol{\varphi}_0(\cdot)=(\varphi_{1,0},\varphi_{2,0})(\cdot):=\left(\phi_{1,0}(\cdot)/\bar{p}_1, 1-\phi_{2,0}(\cdot)/\bar{p}_2 \right)$,
one readily sees that $(\varphi_{1,0},\varphi_{2,0})(\cdot)$ satisfies  
\bqq\label{eq: homo-limit-change}
\left\{\begin{array}{l}
	\bar{d}_1 \varphi_{1,0}''+c_0\varphi_{1,0}'+\bar{g}_{1}(\varphi_{1,0},\varphi_{2,0})=0  \quad \hbox{in }\,\, \R, \vspace{5pt}\\  
	\bar{d}_2  \varphi_{2,0}''+c_0\varphi_{2,0}'+\bar{g}_{2}(\varphi_{1,0},\varphi_{2,0})=0 \quad \hbox{in }\,\, \R, \vspace{5pt}\\  
	(\varphi_{1,0},\varphi_{2,0})(-\infty)=\mathbf{1} \quad\hbox{and}\quad (\varphi_{1,0},\varphi_{2,0})(+\infty)=\mathbf{0}, 
\end{array}\right.
\eqq
and that for each $i=1,2$, $\varphi_{i,0}(\cdot) \in (0,1)$ is decreasing in $\R$, where 
\begin{equation*}
\left\{\begin{array}{l}
	\bar{g}_{1}(v_1,v_2)=\bar{r}_1v_{1} (1-v_{1})-\bar{k}_1\bar{p}_2 v_1 (1-v_2)\vspace{5pt}\\  
	\bar{g}_{2}(v_1,v_2)=-\bar{r}_2v_{2} (1-v_{2})+\bar{k}_2\bar{p}_1 v_1 (1-v_2)\end{array}\right. \quad \hbox{ for }\, (v_1,v_2)\in\R^2.
\end{equation*}

For each $T\in (0,T_*)$ and $i=1,2$, let $D_{i,T}: \mathbb{R}  \rightarrow \mathbb{R}$ be defined by
\bqq\label{defi D_i}
D_{i,T}(t)=\int_0^{t/T} \left(d_{i}(\tau)-\bar{d}_i\right) d \tau \quad \text {for }\, t \in \mathbb{R}, 
\eqq
and $G_{i,T}: \mathbb{R}^3  \rightarrow \mathbb{R}$ be defined by
\bqq\label{defi G_i}
G_{i,T}(t,v_1,v_2)=\int_0^{t/T}\left(g_{i,T}(T\tau, v_1,v_2)-\bar{g}_i(v_1,v_2)\right) d \tau \quad \text {for }\,  (t,v_1,v_2) \in \mathbb{R}^3,
\eqq
where $(g_{1,T}, g_{2,T})$ is given in \eqref{def-giT}. 
It is clear that $D_{i,T}\in C^2(\R)$ is $T$-periodic, and that $G_{i,T}\in C^2(\R^3)$ is $T$-periodic in the first variable. 
Recall that $\mu>0$ is a constant (independent of $T$) satisfying \eqref{set-mu}. 
Let $\psi_{i,T}^{\pm}(\cdot)$ be the unique solutions of the linear problems \eqref{eigenfunction-0} and \eqref{eigenfunction-1} under the normalization \eqref{normalize-supsub}.
Then, since $p_{i,T}(\cdot)$ and $\lambda^{\pm}_{i,T}$ are bounded in $T>0$ (due to \eqref{bound-rak} and \eqref{bound-pT}), 
it follows from \eqref{explicit-eigenfunction} that there exist two constants $\nu_+>\nu_->0$ (independent of $T$) such that 
\bqq\label{bound-eigenfunction}
\nu_-\leq \psi_{i,T}^{\pm}(\cdot)\leq  \nu_+  \, \hbox{ in }\, \R \,  \hbox{ for all } \, 0<T<T_*.
\eqq

Recall further that $\rho(\cdot)\in C^2(\R)$ is a nonnegative function satisfying \eqref{function-rho}.
Our super-solution for \eqref{eq:change of variable} when $T>0$ is small is presented in the following lemma.  

\begin{lem}\label{small T supsub}
	There exists $\varepsilon_0>0$ such that for any $\varepsilon\in (0,\varepsilon_0)$, there exists $T_3=T_3(\varepsilon)>0$ sufficiently small such that for any $T \in (0, T_3)$, the pair of functions $\mathbf{v}_T^+(t,x)=(v_{1,T}^+,v_{2,T}^+)(t,x)$ defined by  
	$$
\begin{aligned}
	v_{i,T}^{+}(t, x)=\, &\varphi_{i,0}\left(\xi_T\right)+q_T(t)\left( \rho(\xi_T)\psi_{i,T}^{ +}(t)+\left( 1-\rho(\xi_T)\right) \psi_{i,T}^{-}(t)\right) \vspace{5pt}\\  
	&\qquad \qquad+T \left( D_{i,T}(t)\varphi_{i,0}''(\xi_T)+G_{i,T}\left(t,\varphi_{1,0}(\xi_T), \varphi_{2,0}(\xi_T)\right)\right)
\end{aligned}$$
for $i=1,2$,  is a super-solution of \eqref{eq:change of variable} for $(t,x) \in (0,+\infty) \times \R$,
where $\xi_T=x-c_0 t+\eta_T(t)$, and $q_T(\cdot)$, $\eta_T(\cdot)$ are $C^1([0,+\infty))$ functions satisfying
	\bqq\label{defi q T}
	q_T(0)=\varepsilon, \quad -A_1 \leq q_T^{\prime}(t)<0<q_T(t) \quad \text {for all }\, t > 0,
	\eqq
	and
	\bqq\label{defi eta T}
	\eta_T(0)=0, \quad-A_2 \leq \eta_T^{\prime}(t)<0 \quad \text {for all }\, t > 0,
	\eqq
	for some positive constants $A_1$, $A_2$ independent of $T$.
\end{lem}

\begin{proof}
	 Let us first choose a small constant  
	\bqq\label{delta_0}
	\sigma_0 =  \min\left\{ \frac{\mu\nu_-}{4\nu_+\theta_+\gamma_+},\, \frac{1}{4} \right\},
	\eqq
	where $\mu$, $\nu_{\pm}$, $\theta_+$, $\gamma_+$ are given by \eqref{set-mu}, \eqref{bound-eigenfunction}, \eqref{bound-rak} and \eqref{bound-pT}, respectively.
	By the last line of \eqref{eq: homo-limit-change} and the fact that $\varphi_{i,0}(\cdot)$ is decreasing in $\R$, 
    there exist $M>0$ and $\beta_0$>0 such that, for each $i=1,2$, 
	\bqq\label{wave small T}
	\left\{\begin{array}{ll}
		0<\varphi_{i,0}(\xi)<\sigma_0 & \text {for all } \, \xi \geq M, \vspace{5pt}\\ 
		1-\sigma_0<\varphi_{i,0}(\xi)<1 & \text {for all }\,  \xi \leq-M, \vspace{5pt}\\ 
			\sigma_0 / 2<\varphi_{i,0}(\xi)<1-\sigma_0 / 2 & \text {for all } \, -M<\xi<M,
	\end{array}\right.
	\eqq
	and
\begin{equation}\label{choose-beta0-samll}
\beta_0\leq -\max _{-M\leq \xi\leq M}\varphi_{i,0}'(\xi).
\end{equation} 

Now, we choose  
\bqq\label{delta_1} 
\varepsilon_0 = \min \left\{\frac{\sigma_0}{2\nu_+},\,\frac{\beta_0}{2(\nu_+-\nu_-)}\right\}, 
\eqq
and for any $\varepsilon \in (0, \varepsilon_0)$, we define
	$$
	\left\{\begin{array}{l}
		\mathcal{L}_1 (\textbf{v}^{+}_{T}):=\partial_t v^{+}_{1,T}-d_{1,T}(t)\partial_{x x} v^{+}_{1,T}-g_{1,T}\left(t, \textbf{v}^{+}_{T}\right), \vspace{5pt}\\ 
		\mathcal{L}_2 (\textbf{v}^{+}_{T}):=\partial_t v^{+}_{2,T}-d_{2,T}(t)\partial_{x x} v^{+}_{2,T}-g_{2,T}\left(t, \textbf{v}^{+}_{T}\right), 
	\end{array}\right.$$
     for $(t, x) \in(0, +\infty) \times \mathbb{R}$.
To show Lemma \ref{small T supsub}, it suffices to find suitable  $C^1([0,+\infty))$ functions $q_T$ and $\eta_T$ satisfying \eqref{defi q T} and \eqref{defi eta T} such that, if $T\in (0,T_*)$ is sufficiently small, there holds $\mathcal{L}_1 (\textbf{v}^{+}_{T}) \geq 0$ and $\mathcal{L}_2 (\textbf{v}^{+}_{T}) \geq 0$ in $(0, +\infty) \times \mathbb{R}$.  We will only verify $\mathcal{L}_1 (\textbf{v}^{+}_{T})\geq 0$, as the verification of $\mathcal{L}_2 (\textbf{v}^{+}_{T}) \geq 0$ is analogous.
	Since $(\varphi_{1,0}, \varphi_{2,0})(\cdot)$ is an entire solution of \eqref{eq: homo-limit-change}, direct computations give that, for $(t, x) \in(0, +\infty) \times \mathbb{R}$
	\begin{equation*}
	\left.\begin{array}{ll}
		\mathcal{L}_1 (\textbf{v}^{+}_{T})=\!\!\! & \eta^{\prime}_T  \varphi'_{1,0}+q^{\prime}_T\left(\rho \psi_{1,T}^{+}+(1-\rho) \psi_{1,T}^{-}\right)- q_Td_{1,T} \rho''\left(\psi_{1,T}^{+}-\psi_{1,T}^{-}\right) \vspace{5pt}\\ 
		&\quad +q_T\left[(\eta^{\prime}_T-c_0)\rho'\left(\psi_{1,T}^{+}-\psi_{1,T}^{-}\right)+\rho\left(\psi_{1,T}^{+}\right)^{ \prime}+(1-\rho)\left(\psi_{1,T}^{-}\right)^{\prime }\right]  \vspace{5pt}\\ 
		&\quad +T(\eta^{\prime}_T-c_0)\left[D_{1,T}(t)\varphi'''_{1,0}+\partial_{v_1}G_{1,T}(t,\boldsymbol{\varphi}_0)\varphi'_{1,0}+\partial_{v_2}G_{1,T}(t,\boldsymbol{\varphi}_0)\varphi'_{2,0}\right]\vspace{5pt}\\ 
		&\quad -Td_{1,T}\left[D_{1,T}(t)\varphi^{(4)}_{1,0}+\partial_{v_1v_1}G_{1,T}(t,\boldsymbol{\varphi}_0)(\varphi'_{1,0})^2+\partial_{v_2v_2}G_{1,T}(t,\boldsymbol{\varphi}_0)(\varphi'_{2,0})^2\right.\vspace{5pt}\\ 
		&\quad \left.+2\partial_{v_1v_2}G_{1,T}(t,\boldsymbol{\varphi}_0)\varphi'_{1,0}\varphi'_{2,0}+\partial_{v_1}G_{1,T}(t,\boldsymbol{\varphi}_0)\varphi''_{1,0}+\partial_{v_2}G_{1,T}(t,\boldsymbol{\varphi}_0)\varphi''_{2,0}\right] \vspace{5pt}\\ 
		&\quad -\left[g_{1,T}\left(t, \textbf{v}^{+}_{T}\right)-g_{1,T}\left(t, \boldsymbol{\varphi}_0\right) \right],
	\end{array}\right.
	\end{equation*}
	where $\rho(\cdot)$, $\varphi_{i,0}(\cdot)$ and their derivatives are evaluated at $\xi_T$, and  $q_T(\cdot)$, $\eta_T(\cdot)$, $\psi_{1,T}^{\pm}(\cdot)$ and their derivatives are evaluated at $t$.  
	
	Let us first derive some rough estimates for $\mathcal{L}_1 (\textbf{v}^{+}_{T})$.
	By the definitions of $D_{1,T}$ and $G_{1,T}$ (see \eqref{defi D_i} and \eqref{defi G_i}), one can find a constant $K_1>0$ independent of $T$ such that for each $i,j=1,2$,
	\begin{equation}\label{DGT-bound}
	\left|D_{1,T}\left(t \right)\right|,\,\left| G_{1,T}\left(t, \boldsymbol{\varphi}_0\right)\right|,\,\left| \partial_{v_i}G_{1,T}\left(t, \boldsymbol{\varphi}_0\right)\right|,\,\left| \partial_{v_iv_j}G_{1,T}\left(t, \boldsymbol{\varphi}_0\right)\right| \leq K_1 
	\end{equation}
for $(t,x)\in (0,+\infty)\times\R$.
    Moreover, since $\boldsymbol{\varphi}_0\in C^4(\R;\R^2)$ connects $\mathbf{0}$ and $\mathbf{1}$, there exits $K_2>0$ such that for each $i=1,2$,
	\begin{equation}\label{wave bound}
	\left|\varphi_{i,0}(\xi)\right|,\,\left|\varphi'_{i,0}(\xi)\right|,\,\left|\varphi''_{i,0}(\xi)\right|,\,\left|\varphi'''_{i,0}(\xi)\right|,\, \left|\varphi^{(4)}_{i,0}(\xi)\right|\leq K_2  \quad \hbox{for }\,\, \xi\in \R.
    \end{equation}
    It is then easily checked that 
	\bqq\label{concrete expression}
\left.\begin{array}{ll}
	\mathcal{L}_1 (\textbf{v}^{+}_{T})\geq\!\!\! & \eta^{\prime}_T  \varphi'_{1,0}+q^{\prime}_T\left(\rho \psi_{1,T}^{+}+(1-\rho) \psi_{1,T}^{-}\right)- q_Td_{1,T} \rho''\left(\psi_{1,T}^{+}-\psi_{1,T}^{-}\right) \vspace{5pt}\\ 
	&\quad +q_T\left[(\eta^{\prime}_T-c_0)\rho'\left(\psi_{1,T}^{+}-\psi_{1,T}^{-}\right)+\rho\left(\psi_{1,T}^{+}\right)^{ \prime}+(1-\rho)\left(\psi_{1,T}^{-}\right)^{\prime }\right]  \vspace{5pt}\\ 
	&\quad -\left[g_{1,T}\left(t, \textbf{v}^{+}_{T}\right)-g_{1,T}\left(t, \boldsymbol{\varphi}_0\right) \right]-L_1T,
\end{array}\right.
\eqq
for $(t,x)\in (0,+\infty)\times \R$, where $L_1=3(|c_0|+A_2)K_1K_2+\theta_+(3K_1K_2+4K_1K_2^2)$
(recall that $\|d_1\|\leq \theta_+$ by \eqref{bound-rak}). Here, we use the requirement that $\eta_T$ satisfies \eqref{defi eta T}.
Furthermore, by \eqref{def-giT}, it is straightforward to compute that 
\begin{equation*}
	\left.\begin{array}{ll}
		&g_{1,T}\left(t, \textbf{v}^{+}_{T}\right)-g_{1,T}\left(t, \boldsymbol{\varphi}_0\right) \vspace{5pt}\\ 
		=\!\!\!&a_{1,T}p_{1,T}\left( v^{+}_{1,T}-\varphi_{1,0}\right)  \left(1-(v^{+}_{1,T}+\varphi_{1,0})\right)-k_{1,T}p_{2,T} \left(v^{+}_{1,T}-\varphi_{1,0}+\varphi_{1,0}\varphi_{2,0}-v^{+}_{1,T}v^{+}_{2,T}\right).
	\end{array}\right.
\end{equation*}
Denote by 
$$H_{i,T}(t,x):=\rho(\xi_T)\psi_{i,T}^{ +}(t)+\left( 1-\rho(\xi_T)\right) \psi_{i,T}^{-}(t)\quad \hbox{for } \,\, (t,x)\in (0,+\infty)\times \R.$$ 
One can find some constant $L_2>0$ independent of $T$ such that 
\begin{equation}\label{estimate-gT}
	\left.\begin{array}{ll}
		g_{1,T}\left(t, \textbf{v}^{+}_{T}\right)-g_{1,T}\left(t, \boldsymbol{\varphi}_0\right) 
		&\!\!\!\leq k_{1,T}p_{2,T}q_T\left(- H_{1,T} +\varphi_{1,0}H_{2,T}  +\varphi_{2,0}H_{1,T}+H_{1,T}H_{2,T}q_T\right)  \vspace{5pt}\\ 
		& \qquad\qquad+a_{1,T}p_{1,T}H_{1,T}q_T \left(1-(v^{+}_{1,T}+\varphi_{1,0})\right) +L_2T
	\end{array}\right.
\end{equation}
for $(t,x)\in (0,+\infty)\times \R$.
With $q_T$ required to satisfy \eqref{defi q T}, it follows from \eqref{delta_1} and \eqref{bound-eigenfunction} that
$0<H_{i,T}\leq \nu_+$ and
\begin{equation}\label{esti=qHT}
 0<q_T H_{i,T} \leq \varepsilon \nu_+ \leq \varepsilon_0\nu_+ \leq \sigma_0.
\end{equation}
By using \eqref{bound-rak} and \eqref{bound-pT}, we have
$-a_{1,T}p_{1,T}H_{1,T}q_T ( D_{i,T}\varphi_{1,0}''+G_{1,T})T \leq L_3T$, where $L_3=\theta_+\gamma_+\sigma_0K_1(K_2+1)$,
and hence,  
$$a_{1,T}p_{1,T}H_{1,T}q_T \left(1-(v^{+}_{1,T}+\varphi_{1,0})\right) \leq a_{1,T}p_{1,T}H_{1,T}q_T\left( 1-(2\varphi_{1,0}+H_{1,T}q_T) \right)+L_3T. $$
It then follows directly from \eqref{estimate-gT} that 
\begin{equation}\label{estimate-gT-new}
	\left.\begin{array}{ll}
		g_{1,T}\left(t, \textbf{v}^{+}_{T}\right)-g_{1,T}\left(t, \boldsymbol{\varphi}_0\right) 
		&\!\!\!\leq k_{1,T}p_{2,T}q_T\left(- H_{1,T} +\varphi_{1,0}H_{2,T}  +\varphi_{2,0}H_{1,T}+H_{1,T}H_{2,T}q_T\right)  \vspace{5pt}\\ 
		& \qquad+a_{1,T}p_{1,T}H_{1,T}q_T \left( 1-(2\varphi_{1,0}+H_{1,T}q_T) \right) +(L_2+L_3)T.
	\end{array}\right.
\end{equation}

Below, we complete the proof by considering three cases.
	
	{\bf Case (a)}: $(t,x)\in (0,+\infty)\times \R$ such that $\xi_T\geq M$.
	
	In this case, with our choice of $\rho$ satisfying \eqref{function-rho}, we have $\rho=\rho^{\prime}=\rho^{\prime \prime}=0$, and $H_{i,T}\equiv \psi_{i,T}^{-}$. Since $\eta_T(t)$ is required to be decreasing in $t$, and since $\varphi_{1,0}(\xi)$ is decreasing in $\xi$, we have $\eta^{\prime}_T  \varphi'_{1,0}\geq 0$. Then, it follows immediately from \eqref{concrete expression}  that
	\begin{equation}\label{L1-largep-small}
		\mathcal{L}_1 (\textbf{v}^{+}_{T})
		\geq  q^{\prime}_T \psi_{1,T}^{-}+ q_T(\psi_{1,T}^{-}) ' -\left[g_{1,T}\left(t, \textbf{v}^{+}_{T}\right)-g_{1,T}\left(t, \boldsymbol{\varphi}_0\right)\right]-L_1T. 
	\end{equation}
Furthermore, since $2\varphi_{1,0}+H_{1,T}q_T>0$,  one readily checks from \eqref{estimate-gT-new}  that
	 \begin{equation*}
	 	\left.\begin{array}{ll}
	 		\mathcal{L}_1 (\textbf{v}^{+}_{T})\!\!\!&\geq  q^{\prime}_T \psi_{1,T}^{-}+ q_T(\psi_{1,T}^{-}) '-a_{1,T}p_{1,T}q_T\psi_{1,T}^{-}+k_{1,T}p_{2,T}q_T\psi_{1,T}^{-} \vspace{5pt}\\ 
	 		&\qquad -k_{1,T}p_{2,T}q_T\left(\varphi_{1,0}\psi_{2,T}^{-}+\varphi_{2,0}\psi_{1,T}^{-}+q_T\psi_{1,T}^{-}\psi_{2,T}^{-}\right)
	 		-(L_1+L_2+L_3)T.
	 	\end{array}\right.
	 \end{equation*}
 Noticing that $\psi_{1,T}^{-}(\cdot)$ is a solution of \eqref{eigenfunction-0} with $i=1$,  we obtain that
\begin{equation*}
\left.\begin{array}{ll}
\displaystyle \mathcal{L}_1 (\textbf{v}^{+}_{T}) \geq q^{\prime}_T \psi_{1,T}^{-}\!\!\!\! & \displaystyle +\lambda_{1,T}^{-}q_T \psi_{1,T}^{-} 
		 -k_{1,T}p_{2,T}\varphi_{2,0}q_T \psi_{1,T}^{-} \vspace{5pt}\\ 
		& \displaystyle -q_T\left(\varphi_{1,0}+q_T \psi_{1,T}^{-}\right) \frac{k_{1,T}p_{2,T} \psi_{2,T}^{-}}{\psi_{1,T}^{-}} \psi_{1,T}^{-}-(L_1+L_2+L_3)T.
\end{array}\right.
\end{equation*} 
Since $\varphi_{2,0}\in(0,\sigma_0)$ and $\varphi_{1,0}+q_T \psi_{1,T}^{-}\in(0,2\sigma_0)$, by the setting of $\sigma_0$ and $\mu$ in \eqref{delta_0} and \eqref{set-mu} respectively, it follows that
$$ \mathcal{L}_1 (\textbf{v}^{+}_{T}) \geq q^{\prime}_T \psi_{1,T}^{-}+2\mu q_T \psi_{1,T}^{-} -\frac{\mu}{4}q_T \psi_{1,T}^{-}-\frac{\mu}{2}q_T \psi_{1,T}^{-}-(L_1+L_2+L_3)T.$$ 
Finally, thanks to \eqref{bound-eigenfunction}, there holds 
 \begin{equation}\label{geqM-samll}
	\mathcal{L}_1 (\textbf{v}^{+}_{T})	\geq \nu_+ q^{\prime}_T +\mu \nu_- q_T-(L_1+L_2+L_3)T.
\end{equation}
Therefore, we can conclude that, with $q_T$ and $\eta_T$ satisfying \eqref{defi q T} and \eqref{defi eta T}, for any small $T>0$,  \eqref{geqM-samll} holds for all $(t,x)\in (0,+\infty)\times \R$ such that $\xi_T\geq M$. 	
	
	{\bf Case (b)}: $(t,x)\in (0,+\infty)\times \R$ such that $\xi_T\leq -M$.
	
   In this case,  there holds $\rho=1$, $\rho^{\prime}=\rho^{\prime \prime}=0$ and $H_{i,T}\equiv \psi_{i,T}^{+}$. Similarly as above, we have
	$$\mathcal{L}_1 (\textbf{v}^{+}_{T})
		\geq  q^{\prime}_T \psi_{1,T}^{+}+ q_T(\psi_{1,T}^{+}) ' -\left[g_{1,T}\left(t, \textbf{v}^{+}_{T}\right)-g_{1,T}\left(t, \boldsymbol{\varphi}_0\right)\right]-L_1T.$$
 Since $q_T\psi_{1,T}^{+}\in(0,\sigma_0)$, by using the second line of \eqref{wave small T}, we have
$1-(2\varphi_{1,0}+q_T\psi_{1,T}^{+}) \leq -1+2\sigma_0$.		
It then follows from \eqref{estimate-gT-new} that		
 \begin{equation*}
 	\left.\begin{array}{ll}
 		g_{1,T}\left(t, \textbf{v}^{+}_{T}\right)-g_{1,T}\left(t, \boldsymbol{\varphi}_0\right) 
 		\leq \!\!\!& -k_{1,T}p_{2,T}q_T\psi_{1,T}^{+}+k_{1,T}p_{2,T}q_T\left(\varphi_{1,0}\psi_{2,T}^{+}+\varphi_{2,0}\psi_{1,T}^{+}+q_T\psi_{1,T}^{+}\psi_{2,T}^{+}\right)    \vspace{5pt}\\ 
 		&\quad +a_{1,T}p_{1,T}q_T\psi_{1,T}^{+}(-1+2\sigma_0)+(L_2+L_3)T.
 	\end{array}\right.
 \end{equation*}
Furthermore, noticing that $\psi_{1,T}^{+}$ is a solution of \eqref{eigenfunction-1} with $i=1$, we obtain that 
\begin{equation*}
\left.\begin{array}{ll}
		\mathcal{L}_1(\textbf{v}^{+}_{T}) \!\!\! & \geq q^{\prime}_T \psi_{1,T}^{+}+\lambda_{1,T}^{+}q_T \psi_{1,T}^{+}-2\sigma_0a_{1,T}p_{1,T}q_T\psi_{1,T}^{+}   +k_{1,T}p_{2,T}(1-\varphi_{2,0})q_T \psi_{1,T}^{+} \vspace{5pt}\\
	&\qquad \displaystyle	-q_T\left(\varphi_{1,0}+q_T \psi_{1,T}^{+}\right) \frac{k_{1,T}p_{2,T} \psi_{2,T}^{+}}{\psi_{1,T}^{+}} \psi_{1,T}^{+}-(L_1+L_2+L_3)T.
\end{array}\right.
\end{equation*}
Since $1-\varphi_{2,0}\geq 0$ and $\varphi_{1,0}+q_T \psi_{1,T}^{+} \in (1-\sigma_0, 1+\sigma_0)$,  it follows from \eqref{delta_0}, the third normalization in \eqref{normalize-supsub} and \eqref{set-mu} that
\begin{equation*}
	\left.\begin{array}{ll}
\mathcal{L}_1 (\textbf{v}^{+}_{T}) \!\!\! &  \displaystyle \geq q^{\prime}_T \psi_{1,T}^{+}+\lambda_{1,T}^{+}q_T \psi_{1,T}^{+} -\frac{\mu}{2}q_T \psi_{1,T}^{+}-\frac{\mu}{2}q_T \psi_{1,T}^{+}-(L_1+L_2+L_3)T \vspace{5pt}\\
& \geq q^{\prime}_T \psi_{1,T}^{+}+\mu q_T \psi_{1,T}^{+}-(L_1+L_2+L_3)T.
\end{array}\right.
\end{equation*}
By using  \eqref{bound-eigenfunction} again, we reach the same estimate as that in \eqref{geqM-samll}, but now $(t,x)\in (0,+\infty)\times \R$ is restricted to satisfy $\xi_T\leq -M$. 

Combining the above, we set $L_0=L_1+L_2+L_3$, and choose
	\bqq\label{chosse T_3}
	T_3=\min\left\{\frac{\varepsilon \mu \nu_-}{2L_0},\,\,T_*\right\}>0. 
	\eqq
For any $0<T< T_3$, let the function $q_T(t)$ satisfy 
	$$
	\nu_+q^{\prime}_T +\mu\nu_- q_T -L_0T=0 \, \text { for }\, t>0,  \quad \text {and} \quad q_T(0)=\varepsilon,
	$$
namely, 
	\bqq\label{function-qT}
	q_T(t)=\frac{L_0}{ \mu \nu_-}T+\left(\varepsilon-\frac{L_0}{ \mu \nu_-}T \right) {\rm exp} \left( -\frac{\mu\nu_-}{\nu_+}t  \right)\quad \text {for } \,  t \geq 0.
	\eqq
	Clearly, for any $0<T<T_3$, $q_T(t)$ satisfies \eqref{defi q T} with $A_1=\varepsilon \mu\nu_-/\nu_+$ independent of $T$, and $\mathcal{L}_1 (\textbf{v}^{+}_{T}) \geq 0$ for $(t, x) \in (0, +\infty) \times \mathbb{R}$ such that $\xi_T \leq -M$   or $\xi_T \geq M$, as soon as $\eta_T$ satisfies \eqref{defi eta T}. 
	
	It remains to find suitable $\eta_T$ satisfying \eqref{defi eta T} such that $\mathcal{L}_1 (\textbf{v}^{+}_{T}) \geq 0$ in the following case. 
	
	{\bf Case (c)}: $(t,x)\in (0,+\infty)\times \R$ such that $-M<\xi_T< M$.
	
	In this case, since $\eta_T(\cdot)$ is required to be decreasing, it follows from \eqref{choose-beta0-samll} that $\eta'_T  \varphi'_{1,0}\geq -\beta_0 \eta'_T$. Moreover, since $\rho(\cdot)$ satisfies \eqref{function-rho},  by using \eqref{delta_1} and \eqref{bound-eigenfunction}, one readily checks that 
	$$ q_T\rho'\left(\psi_{1,T}^{+}-\psi_{1,T}^{-}\right)\eta_T'  \geq \left|q_T\rho'\left(\psi_{1,T}^{+}-\psi_{1,T}^{-}\right) \right |\eta_T' \geq  \varepsilon_0 (\nu_+-\nu_-)   \eta_T' \geq \frac{\beta_0}{2}\eta_T', $$
	and that
	$$-(d_{1,T} \rho''+c_0)\left(\psi_{1,T}^{+}-\psi_{1,T}^{-}\right)q_T \geq  -\left(\|d_1\|+|c_0|\right)\left\|\psi_{1,T}^{+}-\psi_{1,T}^{-}\right\| q_T \geq -(\nu_+-\nu_-)(\theta_++|c_0|) q_T. $$
	Since $q_T(\cdot)$ is decreasing, there holds 
	$$q^{\prime}_T\left(\rho \psi_{1,T}^{+}+(1-\rho) \psi_{1,T}^{-}\right)\geq  \nu_+ q^{\prime}_T. $$
Remember that $\psi_{1,T}^{\pm}$ are  solutions of \eqref{eigenfunction-0} and \eqref{eigenfunction-1} with $i=1$. By using \eqref{bound-rak}, \eqref{bound-pT} and \eqref{seperate-lambda}, one computes that 	
\begin{equation*}
\left.\begin{array}{ll}
		&\left(\rho\left(\psi_{1,T}^{+}\right)^{ \prime}+(1-\rho)\left(\psi_{1,T}^{-}\right)^{\prime }\right)q_T \vspace{5pt}\\
		\geq  & \!\!\! -\left(\left\|\left(\lambda_{1,T}^{+}+\partial_{v_1}g_{1,T}(\cdot,1,1)\right)\psi_{1,T}^{+}(\cdot)  \right\|+\left\|\left(\lambda_{1,T}^{-}+\partial_{v_1}g_{1,T}(\cdot,0,0)\right)\psi_{1,T}^{-}(\cdot)\right\| \right)q_T    \vspace{5pt}\\
		\geq & \!\!\! -6\theta_+\gamma_+ \nu_+q_T.
		\end{array}\right.
\end{equation*}
Let us now turn to estimate $g_{1,T}(t, \textbf{v}^{+}_{T})-g_{1,T}(t, \boldsymbol{\varphi}_0)$. 
 Recall that $0<H_{i,T}\leq \nu_+$ and that $0<q_T H_{i,T} \leq  \sigma_0$ (see \eqref{esti=qHT}). It follows from \eqref{estimate-gT-new} that
$$  g_{1,T}\left(t, \textbf{v}^{+}_{T}\right)-g_{1,T}\left(t, \boldsymbol{\varphi}_0\right)\leq \theta_+\gamma_+ (2\nu_++\nu_+\sigma_0)q_T+(L_2+L_3)T \leq 3\theta_+\gamma_+ \nu_+q_T+(L_2+L_3)T. $$

Combining the above and setting 
$B_0=(\nu_+-\nu_-)(\theta_++|c_0|)+9\theta_+\gamma_+ \nu_+$, 
one sees from  \eqref{concrete expression} that 
$$	\mathcal{L}_1 (\textbf{v}^{+}_{T})  \geq -\frac{\beta_0}{2}\eta_T'+\nu_+ q'_T-B_0q_T-L_0T. $$
Then, with $q_T(t)$ given by \eqref{function-qT}, by choosing $\eta_T(t)$ such that
$$-\frac{\beta_0}{2}\eta_T'+\nu_+ q'_T-B_0q_T-L_0T=0\,\hbox{ for }\, t>0 \quad \hbox{and}\quad \eta_T(0)=0,  $$
we have $\mathcal{L}_1 (\textbf{v}^{+}_{T})  \geq 0$ for all $(t,x)\in (0,+\infty)\times \R$ such that $-M<\xi_T< M$. It is easily computed that 
\begin{equation}\label{function-etaT}
	\eta_T(t)=-\frac{2(B_0+\mu\nu_-)}{\beta_0}\left( \frac{L_0}{ \mu \nu_-}Tt + \left(\varepsilon-\frac{L_0}{ \mu \nu_-}T \right) \frac{\nu_+}{\mu\nu_-} \left(1-{\rm exp} \left( -\frac{\mu\nu_-}{\nu_+}t  \right)  \right) \right)
\end{equation}
for $t\geq 0$. 
It is also straightforward to check that for any $0<T<T_3$ (remember that $T_3$ is given by \eqref{function-qT}), there holds
	$$
	-\frac{3\varepsilon(B_0+\mu\nu_-)}{\beta_0} \leq \eta_T^{\prime}(t) \leq -\frac{\varepsilon(B_0+\mu\nu_-)}{\beta_0}  <0 
	\quad \text{for all }\,  t > 0.
	$$
	As a consequence, the function $\eta_T(t)$ given by \eqref{function-etaT} satisfies \eqref{defi eta T} with $A_2=3\varepsilon (B_0+\mu\nu_-)/\beta_0$. Now, we can conclude that for any $\varepsilon\in (0,\varepsilon_0)$ and any $T \in (0, T_3)$, there holds 
	 $\mathcal{L}_1 (\textbf{v}^{+}_{T})  \geq 0$ for all $(t,x)\in (0,+\infty)\times \R$. 
	 
Proceeding similarly as above, with $q_T(\cdot)$ and $\eta_T(\cdot)$ given by \eqref{function-qT} and \eqref{function-etaT}, one can show that  $\mathcal{L}_2 (\textbf{v}^{+}_{T})  \geq 0$ for all $(t,x)\in (0,+\infty)\times \R$.   
Therefore, $(v^{+}_{1,T}, v^{+}_{2,T})$ is a super-solution of \eqref{eq:change of variable} for $(t,x)\in (0,+\infty)\times \R$.	This completes the proof of Lemma \ref{small T supsub}.
\end{proof}

The following lemma gives the sub-solution of \eqref{eq:change of variable} when $T>0$ is small.

\begin{lem}\label{small T supsub2}
Let $\varepsilon_0>0$ be given by Lemma {\rm \ref{small T supsub}} and, for any $\varepsilon \in (0, \varepsilon_0)$, let $T_3>0$ be given by \eqref{chosse T_3} as in Lemma {\rm \ref{small T supsub}}. Then for any $0<T<T_3$, the pair of functions  $\mathbf{v}_T^-(t,x)=(v_{1,T}^-,v_{2,T}^-)(t,x)$ defined by  
$$
\begin{aligned}
	v_{i,T}^{-}(t, x)=\, &\varphi_{i,0}(\tilde{\xi}_T)-q_T(t)\left(\rho(\tilde{\xi}_T)\psi_{i,T}^{ +}(t)+\left( 1-\rho(\tilde{\xi}_T)\right) \psi_{i,T}^{-}(t)\right) \vspace{5pt}\\  
	&\qquad \qquad-T \left( D_{i,T}(t)\varphi_{i,0}''(\tilde{\xi}_T)+G_{i,T}\left(t, \varphi_{1,0}(\tilde{\xi}_T),\varphi_{2,0}(\tilde{\xi}_T)\right)\right)
\end{aligned}$$
for $i=1,2$,  is a sub-solution of \eqref{eq:change of variable} for $(t,x) \in (0,+\infty) \times \R$,
where $\tilde{\xi}_T=x-c_0 t-\eta_T(t)$, and
$q_T(\cdot)$, $\eta_T(\cdot)$ are $C^1([0,+\infty))$ functions satisfying \eqref{defi q T} and \eqref{defi eta T}, given by \eqref{function-qT} and \eqref{function-etaT}, respectively.
\end{lem}

The proof is analogous to that of Lemma \ref{small T supsub}; therefore, we omit the details. We are now ready to prove the convergence of $c_T$ as $T\to 0^+$.

\begin{proof}[Proof of Theorem {\rm \ref{converge-samllT}}]
	Fix $\varepsilon \in (0, \varepsilon_0)$, and for any $T\in (0,T_3)$, let $\mathbf{v}_T^+(t,x)=(v_{1,T}^+,v_{2,T}^+)(t,x)$ and $\mathbf{v}_T^-(t,x)=(v_{1,T}^-,v_{2,T}^-)(t,x)$ be, respectively, the super-solution and sub-solution of \eqref{eq:change of variable} provided by Lemmas \ref{small T supsub}-\ref{small T supsub2}.  
	Since $q_T(0)=\varepsilon$, $\eta_T(0)=0$, $D_{i,T}(0)=0$ and $G_{i,T}(0,\cdot,\cdot)\equiv 0$, it is easily  seen that
	$$
	v_{i,T}^{-}(0, x) \leq \varphi_{i,0}(x) \leq v_{i,T}^{+}(0, x) \quad \text {for }\,  x \in \mathbb{R},\, i=1,2.
	$$
	Denote by $\mathbf{v}_T(t,x)=(v_{1,T},v_{2,T})(t,x)$ the solution of the Cauchy problem of \eqref{eq:change of variable} with initial data $\mathbf{v}_T(0,x)= (\varphi_{1,0},\varphi_{2,0})(x)$. Since $\mathbf{0} \ll (\varphi_{1,0},\varphi_{2,0})(\cdot) \ll \mathbf{1}$ in $\R$, the strong maximum principle of parabolic equations implies that $\mathbf{0} \ll \mathbf{v}_T(t,x) \ll \mathbf{1}$ for $(t,x)\in (0,+\infty)\times \R$. Furthermore, applying the comparison principle for cooperative systems (see Lemma \ref{comparison principle}), one obtains that  
	\bqq\label{speed cover1}
	v_{i,T}^{-}(t, x) \leq v_{i,T}(t, x) \leq v_{i,T}^{+}(t, x) \quad \text{for }\,  (t,x)\in (0,+\infty)\times \R,\, i=1,2.
	\eqq

	On the other hand, recall that $(\varphi_{1,T},\varphi_{2,T})(x-c_Tt,t)$ is the unique periodic traveling wave of \eqref{eq:change of variable} satisfying the normalization  \eqref{normalize small T wave}. It is known from Proposition \ref{stability} that this wave is globally stable in the sense that for each $i=1,2$, 
	\bqq\label{speed cover2}
	\sup _{x \in \mathbb{R}}\left|v_{i,T}(t, x)-\varphi_{i,T}\left(x-\xi_T^*-c_Tt,t\right)\right| \rightarrow 0  \quad \text {as }\,  t \rightarrow+\infty
	\eqq
	for some constant $\xi_T^* \in \mathbb{R}$.
	
	Combining \eqref{speed cover1} and \eqref{speed cover2}, one can find some $t_0=t_0(T)>0$ sufficiently large such that for each $i=1,2$, 
	$$
	v_{i,T}^{-}(t, x)-\frac{\sigma_0}{4} \leq v_{i,T}(t, x)-\frac{\sigma_0}{4} \leq 
	\varphi_{i,T}\left(x-c_Tt-\xi_T^*,t\right) \leq v_{i,T}(t, x)+\frac{\sigma_0}{4} \leq v_{i,T}^{+}(t, x)+\frac{\sigma_0}{4}
	$$
	for all $t \geq t_0$ and $x \in \mathbb{R}$, where $\sigma_0>0$ is the constant provided by \eqref{delta_0}. Recall from the proof of Lemma \ref{small T supsub} that 
	$0<q_T(\rho\psi_{i,T}^{+}+( 1-\rho) \psi_{i,T}^{-}) \leq \sigma_0$ (see \eqref{esti=qHT}). Furthermore, thanks to \eqref{DGT-bound} and \eqref{wave bound}, even if means redefining $T_3$, one has $|T( D_{i,T}\varphi_{i,0}''+G_{i,T})| \leq \sigma_0/4$. Therefore, by the definitions of $v_{i,T}^{\pm}(t, x)$, there holds
	$$\varphi_{i,0}\left(x-c_0 t-\eta_T(t)\right)-\frac{3}{2}\sigma_0 \leq\varphi_{i,T}\left(x-c_Tt-\xi_T^*,t\right) \leq \varphi_{i,0}\left(x-c_0 t+\eta_T(t)\right)+\frac{3}{2}\sigma_0  $$ 
	for all $t \geq t_0$ and $x \in \mathbb{R}$. Now, choosing $x=c_Tt+\xi_T^*$, $t=nT$ and $i=1$ in the above inequalities, by  \eqref{normalize small T wave} and  the $T$-periodicity of $\varphi_{i,T}$ in its second variable, one obtains that
	\begin{equation*}
	\varphi_{1,0}\left((c_T-c_0) nT+\xi_T^*-\eta_T(nT)\right)-\frac{3}{2}\sigma_0  \leq \frac{1}{2} \leq \varphi_{1,0}\left((c_T-c_0) nT+\xi_T^*+\eta_T(nT)\right)+\frac{3}{2}\sigma_0
	\end{equation*}
	for all $n\in \N$ sufficiently large, whence by the fact that $\sigma_0 \in (0,1/4)$,  there holds
   $$\varphi_{1,0}\left((c_T-c_0) nT+\xi_T^*-\eta_T(nT)\right)  \leq \frac{7}{8} \quad \hbox{and} \quad  \varphi_{1,0}\left((c_T-c_0) nT+\xi_T^*+\eta_T(nT)\right)\geq \frac{1}{8}.$$
  This together with the definition of traveling wave implies that there exists a constant $M_0>0$ independent of $T$ such that
	$$
	\left|(c_T-c_0) nT+\xi_T^*\right| \leq\left|\eta_T(nT)\right|+M_0 \quad  \text{for all large }\,  n\in \N.
	$$
	Remember that the function $\eta_T(t)$ is given by \eqref{function-etaT}. Multiplying both sides of the above inequality by $(nT)^{-1}$ yields that
	$$
	\left|  c_T -c_0\right|
		\leq  \frac{\left|\xi_T^*\right|+M_0 }{nT}+\frac{2(B_0+\mu\nu_-)L_0}{\beta_0\mu \nu_-} T+\frac{1}{nT} \frac{2(B_0+\mu\nu_-)\varepsilon\nu_+}{\beta_0\mu\nu_-}.$$
Finally, passing to the limits as $n\rightarrow+\infty$ gives that $	\left|c_T-c_0\right| \leq M_1 T$, where 
 $M_1=2L_0(B_0+\mu\nu_-)/(\beta_0\mu \nu_-)$. It is clear that this constant is independent of $T$, and the proof of Theorem \ref{converge-samllT} is thus complete.
\end{proof}

%%%%%%%%%%%%%%%%%%%%%%%%%%%%%%%%%%%%%%%%%%%%%%%%%%%%%%%%%%%%%%%%%%%%%%%%%%%%%%%
%%%%%%%%%%%%%%%%%%%%%%%%%%%%%%%%%%%%%%%%%%%%%%%%%%%%%%%%%%%%%%%%%%%%%%%%%%%%%%

\section{Proof of Theorem \ref{Theo large T}}

In this section, under assumption (A2), we prove the existence of periodic traveling wave for sytem \eqref{eq:main} and the convergence of wave speed in slowly oscillating environments. While the proof strategy parallels that of the rapidly oscillating case, the slow oscillation case introduces greater difficulties. Firstly, the (normalized) eigenfunctions to the linear problems \eqref{eigenfunction-0}-\eqref{eigenfunction-1} may fail to be uniformly bounded for large $T$, in contrast to the uniform boundedness crucial for proving Theorem \ref{Theo large T}. To overcome this difficulty, we use various solutions of system \eqref{eq:frozen} with frozen coefficients to construct sub- and super-solutions (see Lemmas \ref{linearly stable claim large T} and \ref{large T supsub} below). The construction requires establishing the convergence of the (rescaled) semi-trivial states  $(p_{1,T},0)$ and $(0,p_{2,T})$ as $T\to+\infty$, estimating the corresponding convergence rate (see Lemma \ref{uniform cover s}), and deriving certain uniform estimates of the traveling wave for \eqref{eq:frozen} with respect to the frozen coefficient $s$ (see Proposition \ref{fix s continue}). Secondly, when $T\to+\infty$, a coexistece periodic state of \eqref{eq:main} may approach an arbitrary entire solution of the associated kinetic system for \eqref{eq:frozen} (including equilibria and connecting orbits). To show the instability of such states, we need to classify all entire solutions of the underlying autonomous ODE system, and characterize their stability properties (see Lemma \ref{linearly unstable claim large T}).

\subsection{Existence of bistable traveling waves when $T$ is large}

In this subsection, we apply Proposition \ref{existence-fix} to prove the existence of periodic traveling wave when $L$ is large. The main result is restated as follows. 

\begin{theo}\label{exitstence wave large}
	Let  {\rm (A2)} hold. There exists $T^*>0$ such that, for any $T>T^*$, system \eqref{eq:change of variable} admits a periodic traveling wave $\left(v_{1,T}, v_{2,T}\right)(t,x)=\left(\varphi_{1,T}, \varphi_{2,T}\right)(x-c_Tt, t)$ connecting $\mathbf{0}$ and $\mathbf{1}$.
\end{theo} 

The proof of Theorem \ref{exitstence wave large} is divided into several lemmas that are parallel to those in Section 3.1. Recall that $p_i(\cdot)=r_i(\cdot)/a_i(\cdot)$. The first lemma concerns the convergence of $p_{i,T}(Ts)$ to $p_i(\cdot)$ as $T\to+\infty$.  We point out that this result holds without requiring (A2).

\begin{lem}\label{uniform cover s}
There exists $\tilde{T}_1>0$ such that for each $i=1,2$,  
\begin{equation}\label{pit-pi-large}
\left|p_{i,T}(sT)-p_i(s)\right|\leq \frac{C_0}{T} 
\end{equation}
for all $s\in\R$ and $T\geq \tilde{T}_1$, where $C_0>0$ is a constant independent of $T$ and $s$. Consequently, $p_{i,T}(sT)$ converges to $p_i(s)$ as $T\rightarrow+\infty$ uniformly in $s\in\R$.		
\end{lem}

\begin{proof}
We only consider the case that $i=1$, since the proof in the case that $i=2$ is identical. It is clear that $p_{1,T}(sT)$ is $1$-periodic in $s$, and a direct calculation gives that
\begin{equation*}
p_{1,T}(sT)=\frac{p_{1,T}^0 {\rm exp}\left\{T \displaystyle\int_0^{s} r_1(\tau) d \tau\right\}}{1+p_{1,T}^0 T \displaystyle\int_0^{s}  {\rm exp} \left\{T\int_0^\tau r_1(y) d y \right\} a_1(\tau) d \tau}, 	
\end{equation*}
where
$$	p_{1,T}^0=\frac{ \displaystyle {\rm exp}\left\{T\int_0^1 r_1(\tau) d\tau\right\}-1} { \displaystyle T\int_0^1 {\rm exp}\left\{T\int_0^\tau r_1(y) d y\right\} a_1(\tau) d \tau}. $$
By the $1$-periodicity of $p_{1,T}(sT)$ and $p_1(s)$, it suffices to show that \eqref{pit-pi-large} holds for all $s\in [1,2]$ and all large $T>0$. To this end, for each $T>0$, let us denote 
\begin{equation*}
R_T(s):=T\int_0^{s}  {\rm exp} \left\{T\int_0^\tau r_1(y) d y \right\} a_1(\tau) d \tau \quad\hbox{for }\, s\in [1,2]. 
\end{equation*}
Setting $g(\cdot):=a_1(\cdot)/r_1(\cdot)$ in $\R$ and integrating by parts gives that 
$$R_T(s)={\rm exp}\left\{T \displaystyle\int_0^{s} r_1(\tau) d \tau\right\} g(s)-g(0)- \int_0^{s}  {\rm exp} \left\{T\int_0^\tau r_1(y) d y \right\} g'(\tau) d \tau. $$
Integrating by parts again yields that
\begin{equation}\label{def-QTs} 
\left.\begin{array}{ll}
&R_T(s)={\rm exp}\left\{T \displaystyle\int_0^{s} r_1(\tau) d \tau\right\} g(s)-g(0) - \displaystyle\frac{1}{T} \left\{ {\rm exp}\left\{T \displaystyle\int_0^{s} r_1(\tau) d \tau\right\} \frac{g'(s)}{r_1(s)}-\frac{g'(0)}{r_1(0)} \right\}  \vspace{7pt}\\
&\qquad\qquad\qquad\qquad \qquad\qquad \displaystyle +\, \frac{1}{T} \int_0^{s}  {\rm exp} \left\{T\int_0^\tau r_1(y) d y \right\} \left(\frac{g'}{r_1}\right)'(\tau) d \tau. \end{array}\right.
\end{equation} 
Remember that the functions $a_1(\cdot)$ and $r_1(\cdot)$ are positive, $1$-periodic and of class $C^2(\R)$. Then, one can find some constant $C_1>0$ such that 
$$ |g(s)|+\left|\frac{g'}{r_1}(s)\right|+ \left|\left(\frac{g'}{r_1}\right)'(s)\right| \leq C_1 \quad\hbox{for all }\, s\in\R. $$
This implies that 
$$ \left|\frac{1}{T} \int_0^{s}  {\rm exp} \left\{T\int_0^\tau r_1(y) d y \right\} \left(\frac{g'}{r_1}\right)'(\tau) d \tau\right|\leq \frac{2C_1}{T} {\rm exp}\left\{T \displaystyle\int_0^{s} r_1(\tau) d \tau\right\} \quad\hbox{for all } \,s\in [1,2].$$
Let $\tilde{T}_1>1$ be sufficiently large such that 
\begin{equation}\label{choose-tildeT}
{\rm exp}\left\{T \displaystyle\int_0^{1} r_1(\tau) d \tau\right\} \geq T \,\quad \hbox{ for all }\, T\geq \tilde{T}_1. 
\end{equation}
It then follows from \eqref{def-QTs} that exists exists a constant $C_2>C_1$ (independent of $T$) such that
\begin{equation}\label{QTs-exp}
	\left| R_T(s){\rm exp}\left\{-T \displaystyle\int_0^{s} r_1(\tau) d \tau\right\}  -g(s) \right| \leq \frac{C_2}{T} \quad \hbox{for all }\, T\geq \tilde{T}_1,\,s\in [1,2]. 
\end{equation}
Since $p_1(\cdot)=g^{-1}(\cdot)$, this is equivalent to the following 
\begin{equation}\label{QTs-exp-1}
	\left| R_T^{-1}(s){\rm exp}\left\{T \displaystyle\int_0^{s} r_1(\tau) d \tau\right\}  -p_1(s) \right| \leq \frac{C_3}{T} \quad \hbox{for all }\, T\geq \tilde{T}_1,\,s\in [1,2],
\end{equation}
for some constant $C_3>0$ independent of $T$.

We are now ready to complete the proof. Let $T\geq \tilde{T}_1$ and $s\in [1,2]$ be arbitrary.
By the expression of $p_{1,T}(sT)$, we have 
\begin{equation}\label{ex-p1t}
p_{1,T}(sT)= \left( \left(p_{1,T}^0\right)^{-1} {\rm exp}\left\{-T \displaystyle\int_0^{s} r_1(\tau) d \tau\right\} + R_T(s){\rm exp}\left\{-T \displaystyle\int_0^{s} r_1(\tau) d \tau\right\} \right)^{-1}. 
\end{equation} 
On the one hand, since $p_{1,T}^0$ is positive, it follows directly from \eqref{ex-p1t} and \eqref{QTs-exp-1} that 
$$  p_{1,T}(sT) \leq R_T^{-1}(s){\rm exp}\left\{T \displaystyle\int_0^{s} r_1(\tau) d \tau\right\} \leq p_1(s)+\frac{C_3}{T}.$$
On the other hand, by the expression of $p_{1,T}^0$, we have  
$$p_{1,T}^0{\rm exp}\left\{T \displaystyle\int_0^{s} r_1(\tau) d \tau\right\} \geq 
p_{1,T}^0{\rm exp}\left\{T \displaystyle\int_0^{1} r_1(\tau) d \tau\right\}
=\frac{{\rm exp} \left\{T \displaystyle\int_0^1 r_1(\tau) d \tau\right\}-1}{Q_T(1){\rm exp}\left\{-T \displaystyle\int_0^{1} r_1(\tau) d \tau\right\} }. $$
And hence, by using \eqref{choose-tildeT} and \eqref{QTs-exp}, we obtain
$$p_{1,T}^0{\rm exp}\left\{T \displaystyle\int_0^{s} r_1(\tau) d \tau\right\}  \geq  \frac{T-1}{g(1)+C_2T^{-1}}.$$
It then further follows from \eqref{ex-p1t} that, by using \eqref{QTs-exp-1} again, 
$$ p_{1,T}(sT) \geq \left( \frac{g(1)+C_2T^{-1}}{T-1} + R_T(s){\rm exp}\left\{-T \displaystyle\int_0^{s} r_1(\tau) d \tau\right\} \right)^{-1} \geq p_1(s)-\frac{C_4}{T}$$
for some constant $C_4>0$ independent of $T$. Now, taking $C_0=\max\{C_3,C_4\}$, we can conclude that \eqref{pit-pi-large} holds for all $s\in [1,2]$ and $T\geq \tilde{T}_1$. Due to the periodicity of $p_{1,T}(sT)$ and $p_1(s)$, this ends the proof of  Lemma \ref{uniform cover s}.
\end{proof}

Next, we verify the conditions in Proposition \ref{existence-fix} when $T$ is large. Let us first set a notation which will be frequently used in the sequel. By (A2) and the $1$-periodicity of the functions $r_i(\cdot)$, $k_i(\cdot)$ and $p_i(\cdot)$, one finds some constant $\gamma_0>0$ such that 
\bqq\label{gamma_0}
\left\{\begin{array}{l}
r_1(s)-k_1(s)p_2(s)\leq-\gamma_0, \vspace{5pt} \\
r_2(s)-k_2(s)p_1(s)\leq-\gamma_0,\vspace{5pt} \\
r_i(s)\geq\gamma_0,  
\end{array}\right. \quad \hbox{for all }\,\,s\in\R,\,i=1,2.
\eqq

\begin{lem}\label{linearly stable claim large T}
Let  {\rm (A2)} hold. There exists $\tilde{T}_2>0$ such that \eqref{0-1-stable} holds for all $T> \tilde{T}_2$, and that for the map $Q_T: [0,1]^2\to[0,1]^2$, the fixed points $\mathbf{0}$ and $\mathbf{1}$ are, respectively, strongly stable from above and below, uniformly with respect to $T>\tilde{T}_2$ in the sense of Definition {\rm \ref{defi-stability}}.  
\end{lem}

\begin{proof}
Let us first verify that \eqref{0-1-stable} holds for all large $T>0$. Indeed, it follows directly from Lemma \ref{uniform cover s} that, as $T\to+\infty$,  
$$\frac{1}{T}\int_0^T \left[a_{1,T}(t) p_{1,T}(t)-k_{1,T}(t)p_{2,T}(t)\right]dt \to \int_0^1 \left[a_1(t) p_1(t)-k_1(t)p_{2}(t)\right]dt,  $$
and
$$\frac{1}{T}\int_0^T \left[a_{2,T}(t)p_{2,T}(t)-k_{2,T}(t)p_{1,T}(t)\right]dt \to 	\int_0^1 \left[a_{2}(t)p_{2}(t)-k_{2}(t)p_{1}(t)\right]dt.$$
Due to condition (A2), both limits are negative constants. This immediately gives \eqref{0-1-stable} when $T>0$ is sufficiently large, whence by Lemma \ref{lem-stable-0-1}, the fixed points  $\mathbf{0}$ and $\mathbf{1}$ are, respectively, strongly stable from above and below. 

Next, we show that the strong stability holds uniformly in large $T$. We only prove that $\mathbf{0}$ is strongly stable from above uniformly in $T$, as the proof for the uniform stability of $\mathbf{1}$ is analogous. 
By (A2) and Lemma \ref{uniform cover s}, one finds some $\tilde{T}_2>0$ sufficiently large such that  
\bqq\label{prop cover}
\left\{\begin{array}{rl}
\displaystyle	-a_{1}(s) p_{1,T}(sT)+k_{1}(s)p_{2,T}(sT) \!\!\!& \displaystyle \geq-r_{1}(s) +k_{1}(s)p_{2}(s)-\frac{\gamma_0}{4}\quad \hbox{for  }\,\,s\in \R,   \vspace{5pt} \\
\displaystyle	a_{2}(s)p_{2,T}(sT) \!\!\! & \displaystyle \geq r_2(s)-\frac{\gamma_0}{8}\quad \hbox{for  }\,\,s\in \R,
\end{array}\right.
\eqq
for all $T\geq\tilde{T}_2$, where $\gamma_0$ is the positive constant provided by \eqref{gamma_0}. Let $\tilde{\epsilon}_0>0$ be a small constant satisfying 
 \bqq\label{tilde epsilon_1}
0<\tilde{\epsilon}_0\leq \min \left\{ 1, \,  \frac{ 8\theta_+ \gamma_+}{7\gamma_0}\right\},
\eqq
where $\theta_+$ and $\gamma_+$ are the positive constants given by \eqref{bound-rak} and \eqref{bound-pT}, respectively. 
In our arguments below, we let $\epsilon \in (0,\tilde{\epsilon}_0)$ and $T\geq \tilde{T}_2$ be arbitrary, and verify that 
the pair of functions
\begin{equation}\label{super-large-stable}
(v_{1,\epsilon},v_{2,\epsilon})(t):=\left(\epsilon C_1e^{-\frac{\gamma_0}{2}t},\, \epsilon C_2e^{-\frac{\gamma_0}{2}t}\right)\quad \hbox{for }\,\,t\geq 0 
\end{equation}
is a strict super-solution of \eqref{eq:change of variable} in the sense of Definition \ref{defi sup sub}, where 
$$ C_2=\frac{\gamma_0}{8 \theta_+\gamma_+}\quad \hbox{and}\quad C_1=\frac{\gamma_0C_2}{4\theta_+\gamma_+}.$$		
It is straightforward to compute that for $t\geq 0$, 		
$$
\begin{aligned}
	&v'_{1,\epsilon}(t)- g_{1,T}\left( t,v_{1,\epsilon},v_{2,\epsilon}\right)  \vspace{5pt} \\
	> \, & \displaystyle v'_{1,\epsilon}(t) +\left(-a_{1,T}\left(t \right) p_{1,T}(t)+k_{1,T}(t)p_{2,T}(t)\right)v_{1,\epsilon}(t)-\left(k_{1,T}(t)p_{2,T}(t)v_{2,\epsilon}(t)\right)v_{1,\epsilon}(t) \vspace{5pt} \\
	= \, & \displaystyle -\frac{\gamma_0}{2}v_{1,\epsilon}(t) +\left(-a_{1}\left(t/T \right) p_{1,T}\left(t\right)+k_{1}\left(t/T \right)p_{2,T}(t)\right)v_{1,\epsilon}(t)-\left(k_{1}\left(t/T \right)p_{2,T}\left(t\right)v_{2,\epsilon}(t)\right)v_{1,\epsilon}(t) \vspace{5pt} \\
	\geq  \, & \displaystyle   -\frac{\gamma_0}{2}v_{1,\epsilon}(t)+ \left(-r_{1}(t/T) +k_{1}(t/T)p_{2}(t/T)-\frac{\gamma_0}{4}\right)v_{1,\epsilon}(t)-\frac{\gamma_0}{8}v_{1,\epsilon}(t), 
\end{aligned}$$ 
where the last inequality follows from the first line of \eqref{prop cover}, the choice of $C_2$ and the fact that $\epsilon<\tilde{\epsilon}_0\leq 1$. Then, by using the first inequality of \eqref{gamma_0}, we obtain 
$$v'_{1,\epsilon}(t)- g_{1,T}\left( t,v_{1,\epsilon},v_{2,\epsilon}\right) >  \left( -\frac{\gamma_0}{2} +\gamma_0-\frac{\gamma_0}{4}-\frac{\gamma_0}{8}\right)v_{1,\epsilon}(t) \geq 0. $$
On the other hand, it is easily checked that for $t \geq0$,
$$
\begin{aligned}
&v'_{2,\epsilon}(t)- g_{2,T}\left( t,v_{1,\epsilon},v_{2,\epsilon}\right)  \vspace{5pt} \\
> \, & \displaystyle v'_{2,\epsilon}(t)+a_{2,T}(t)p_{2,T}(t)(1-v_{2,\epsilon}(t))v_{2,\epsilon}(t)-k_{2,T}(t)p_{1,T}(t)v_{1,\epsilon}(t)  \vspace{5pt} \\
= \, & \displaystyle -\frac{\gamma_0}{2}v_{2,\epsilon}(t) +a_{2}(t/T)p_{2,T}(t)(1-v_{2,\epsilon}(t))v_{2,\epsilon}(t)-\frac{C_1}{C_2}k_{2}(t/T)p_{1,T}(t)v_{2,\epsilon}(t). 
\end{aligned}$$ 
By the second line of \eqref{prop cover} and the choice of $C_1$ and $C_2$, we have
$$v'_{2,\epsilon}(t)- g_{2,T}\left( t,v_{1,\epsilon},v_{2,\epsilon}\right) > -\frac{\gamma_0}{2}v_{2,\epsilon}(t)+ \left( r_2(t/T)-\frac{\gamma_0}{8}\right) (1-\tilde{\epsilon}_0C_2)v_{2,\epsilon}(t) -\frac{\gamma_0}{4}v_{2,\epsilon}(t).$$
It then follows from the third inequality of \eqref{gamma_0} and \eqref{tilde epsilon_1} that 
$$v'_{2,\epsilon}(t)- g_{2,T}\left( t,v_{1,\epsilon},v_{2,\epsilon}\right) > \left(-\frac{\gamma_0}{2}+\frac{6}{7}\left(\gamma_0-\gamma_0/8\right)-\frac{\gamma_0}{4} \right)v_{2,\epsilon}(t)\geq 0.$$
Combining the above, we see that for any $\epsilon \in (0,\tilde{\epsilon}_0)$ and $T\geq \tilde{T}_2$, $(v_{1,\epsilon},v_{2,\epsilon})(t)$ is a strict spatially homogeneous super-solution of \eqref{eq:change of variable} in $t\geq 0$. Finally, following the same lines as those used at the end of the proof of Lemma \ref{linearly stable claim}, one can derive that 
$ \epsilon \mathbf{e} \gg Q_T[\epsilon \mathbf{e}]$ with $\mathbf{e}=(\epsilon C_1, \epsilon C_2)$. Therefore, $\mathbf{0}$ is strongly stable from above uniformly in  $T\geq \tilde{T}_2$.  This completes the proof of Lemma \ref{linearly stable claim large T}.
\end{proof}

By the proof of Lemma \ref{linearly stable claim large T}, one can conclude that for any $T>\tilde{T}_2$ and $\tau\in\R$, the pair of functions  $(v_{1,\epsilon},v_{2,\epsilon})(t)$ provide by \eqref{super-large-stable} is also a strict super-solution of the following system 
\begin{equation*}
\left\{\begin{array}{ll}
	\displaystyle \frac{dv_{1}}{dt}=g_{1,T}\left(t+\tau, v_{1}, v_{2}\right), & t >0,\vspace{5pt} \\
	\displaystyle \frac{dv_{2}}{dt}=g_{2,T}\left(t+\tau, v_{1}, v_{2}\right), & t >0.
\end{array}\right.
\end{equation*}
One then infers that there exists  $\tilde{\epsilon}_1>0$ independent of $T$ such that any $T$-periodic solution $\bar{\mathbf{v}}$ of \eqref{eq:variable-kinetic} with $(0,0) \ll \bar{\mathbf{v}}(\cdot) \ll(1,1)$ must satisfy $\bar{\mathbf{v}}(t)\not\in \{(x_1,x_2)\in [0,1]^2: 0\leq x_i\leq \tilde{\epsilon}_1, i=1,2\}$ for all $t\in\R$. By constructing an analogous sub-solution, there holds $\bar{\mathbf{v}}(t)\not\in \{(x_1,x_2)\in [0,1]^2: 1-\tilde{\epsilon}_1 \leq x_i\leq 1, i=1,2\}$ for all $t\in\R$. This observation will be useful in verifying the linear instability of $\bar{\mathbf{v}}$ when $T$ is large, as stated in the following lemma.

\begin{lem}\label{linearly unstable claim large T}
	There exists $\tilde{T}_3>0$ such that $\lambda_1(\bar{\mathbf{v}}, T)<0$ for every $T>\tilde{T}_3$ and for every $T$-periodic solution $\bar{\mathbf{v}}$ of \eqref{eq:variable-kinetic} with $(0,0) \ll \bar{\mathbf{v}} \ll(1,1)$, where $\lambda_1(\bar{\mathbf{v}}, T)$ is the principal Floquet multiplier for problem \eqref{inter-eigen}. 	
\end{lem}
	
To facilitate the presentation of our proof, we need to classify the asymptotic dynamics of the following kinetic system with coefficients frozen at an arbitrary $s\in\R$: 
\bqq\label{eq:fix-s-Lv} 
\left\{\begin{array}{l}
	\displaystyle  \frac{dv_{1}}{dt}=g_1(s,v_{1},v_{2}),\vspace{5pt} \\
	\displaystyle \frac{dv_{2}}{dt} =g_2(s,v_{1},v_{2}),
\end{array}\right.
\eqq
for $t>0$, where
\begin{equation}\label{g-frozen}
	\left\{\begin{array}{ll}
		g_1(s,v_{1},v_{2})=r_1(s)v_{1} (1-v_{1}) -k_1(s)p_2(s) v_{1} (1-v_{2}),\vspace{5pt} \\
		g_2(s,v_{1},v_{2})=-r_2(s)v_{2} (1-v_{2})+k_2(s)p_1(s) v_{1} (1-v_{2}).
	\end{array}\right.
\end{equation}
Under the condition (A2), it is clear that for each $s\in\R$, the autonomous system \eqref{eq:fix-s-Lv} has only four equilibria in $[0,1]^2$:   
\begin{equation}\label{4s-equilibrium}
	\mathbf{0}=(0,0),\quad \mathbf{1}=(1,1), \quad \mathbf{e}_0=(0,1) \ \text{ and }\  \textbf{e}^*_{s}=(v_{1,s}^*,v_{2,s}^*)\in (0,1)^2,
\end{equation}
where $\textbf{e}^*_{s}$ is a saddle equilibrium.  The following lemma provides a separatrix curve governing the asymptotic behavior of solutions of \eqref{eq:fix-s-Lv}.

\begin{lem}\label{separatrix}{\rm(\cite{M4} Proposition 2.1)}
	Let {\rm (A2)} hold, and for any $s\in\R$, let $(v_1,v_2)(t;s)$ be the solution of \eqref{eq:fix-s-Lv} with initial values $(v_1,v_2)(0;s) \in [0,1]^2$. Then, there exists a decreasing function $h(x)$ $($depending on $s$$)$ defined on $(0, x_{\infty})$ with $x_{\infty} \in (v_{1,s}^*, +\infty]$ such that
	$$\begin{array}{l}
		\text { if } v_2(0;s)<h(v_1(0;s)), \text { then } \lim _{t \rightarrow +\infty}(v_1,v_2)(t;s)=\mathbf{0}; \vspace{5pt} \\
		\text { if } v_2(0;s)>h(v_1(0;s)), \text { then } \lim _{t \rightarrow +\infty}(v_1,v_2)(t;s)=\mathbf{1};\vspace{5pt} \\
		\text { if } v_2(0;s)=h(v_1(0;s)), \text { then } \lim _{t \rightarrow +\infty}(v_1,v_2)(t;s)=\mathbf{e}^*_{s}.
	\end{array}$$
	Moreover, the curve $h(x)$ satisfies $\lim _{x \rightarrow x_{\infty}} h(x)=-\infty$ if  $x_{\infty}<+\infty$, and
	$$\lim_{t\to 0^+}h(t)=1,\quad h\left(v_{1,s}^*\right)=v_{2,s}^*.$$
\end{lem}

\begin{proof}[Proof of Lemma {\rm \ref{linearly unstable claim large T}}]
	We argue by contradiction, and assume that the conclusion of Lemma \ref{linearly unstable claim large T} is not true. Then, there are some sequences $(T_n)_{n \in \mathbb{N}}$ in $(0,+\infty)$, $(\bar{\mathbf{v}}_n)_{n \in \mathbb{N}}$ and $(\boldsymbol \psi_n)_{n \in \mathbb{N}}$ in $C^1(\R;\R^2)$ such that, $T_n \rightarrow +\infty$ as $n \rightarrow+\infty$ and, for each $n \in \mathbb{N}$, $\bar{\mathbf{v}}_n=(\bar{v}_{1, n}, \bar{v}_{2, n})$ satisfies \eqref{unstable  equilibria},
	$\boldsymbol \psi_n=(\psi_{1, n}, \psi_{2, n})$ satisfies \eqref{inter-eigen-n} together with the normalization \eqref{normalize unstable eigenfunction}, and  $\lambda_1(\bar{\mathbf{v}}_n,T_n)\geq 0$. 
	Similarly to the proof of Lemma \ref{linearly unstable claim}, one observes that the sequence $(\lambda_1\left(\bar{\mathbf{v}}_n, T_n\right))_{n \in \mathbb{N}}$ is bounded, and hence, up to extraction of some subsequence, there is a real number $\tilde{\lambda}_\infty \geq 0$ such that 
	$\lambda_1\left(\bar{\mathbf{v}}_n, T_n\right) \rightarrow \tilde{\lambda}_\infty$ as $n \rightarrow+\infty$.
	By the periodicity and \eqref{normalize unstable eigenfunction}, for each $n\in\N$, there exists $t_n\in[0,T_n]$ such that 
	$\boldsymbol \psi_n$ satisfies 
	\bqq\label{unstable eigenfunction find a t_n}
	\psi_{1, n}(t_n)+ \psi_{2, n}(t_n)=1.
	\eqq	
	Up to extraction of a subsequence, one can assume that $t_n/T_n \rightarrow s_{\infty} \in[0,1]$ as $n \rightarrow+\infty$. It then follows directly from Lemma \ref{uniform cover s} that $p_{i,T_n}(t+t_n)\rightarrow p_i(s_\infty)$ as $n \rightarrow+\infty$ locally uniformly in $t\in\R$. Furthermore, by the Arzela-Ascoli theorem,  possibly up to a further subsequence, one can find $\mathbf{v}_{\infty}=(v_{1,\infty},v_{2,\infty}) \in C^1(\R;\R^2)$ and $\boldsymbol \psi_\infty=(\psi_{1, \infty},\psi_{2, \infty})\in C^1(\R;\R^2)$ such that as $n\to+\infty$, 
	\begin{equation}\label{converge-vn-phsn}
	\bar{\mathbf{v}}_n(\cdot+t_n)\to \mathbf{v}_{\infty}(\cdot)\quad\hbox{and}\quad \boldsymbol \psi_n(\cdot+t_n) \to \boldsymbol \psi_\infty(\cdot)\quad\hbox{in}\quad C^1_{loc}(\R).
	\end{equation}   
	It is easily seen from \eqref{unstable  equilibria} and \eqref{inter-eigen-n} that, $\mathbf{v}_{\infty}(t)=(v_{1,\infty},v_{2,\infty})(t)$ solves
\bqq\label{eq: limits v_i}
\left\{\begin{array}{l}  
	\displaystyle \frac{d v_{1,\infty}}{dt}=g_1(s_\infty,v_{1,\infty},v_{2,\infty})   \quad \text{in }\, \mathbb{R}, \vspace{5pt} \\ \displaystyle \frac{d v_{2,\infty}}{dt}=g_2(s_\infty,v_{1,\infty},v_{2,\infty})  \quad \text{in }\, \mathbb{R},   \vspace{5pt} \\
	\mathbf{0} \leq \mathbf{v}_{\infty} \leq \mathbf{1}\quad  \text { in } \mathbb{R},
		\end{array}\right.
		\eqq
	and that $\boldsymbol \psi_\infty(t)=(\psi_{1, \infty},\psi_{2, \infty})(t)$ satisfies
	\bqq\label{eq: limits w_i}
 \frac{d \boldsymbol \psi_\infty}{dt}- A(\mathbf{v}_{\infty}) \boldsymbol \psi_\infty= \tilde{\lambda}_\infty \boldsymbol \psi_\infty \quad \text{in }\, \mathbb{R}, 
	\eqq
	where $(g_1,g_2)(s_{\infty},\cdot,\cdot)$ is defined as in \eqref{g-frozen}, and $A(\mathbf{v}_\infty)=(a_{i,j}(t))_{1\leq i,j\leq 2}$ is a matrix function given by
\begin{equation*}
	\left\{\begin{array}{l}
	a_{11}(t)=r_1(s_\infty)-2r_1(s_\infty)v_{1,\infty}(t)-k_1(s_\infty)p_2(s_\infty)(1-v_{2,\infty}(t)),\vspace{5pt} \\
	a_{12}(t)=k_1(s_\infty)p_2(s_\infty)v_{1,\infty}(t),\vspace{5pt} \\
	a_{21}(t)=k_2(s_\infty)p_1(s_\infty)(1-v_{2,\infty}(t)),\vspace{5pt} \\
	a_{22}(t)=-r_2(s_\infty)+2r_2(s_\infty)v_{2,\infty}(t)-k_2(s_\infty)p_1(s_\infty)v_{1,\infty}(t).
		\end{array}\right. 
\end{equation*}
Furthermore, by the normalization \eqref{normalize unstable eigenfunction} and \eqref{unstable eigenfunction find a t_n}, there holds
\begin{equation}\label{psiinfty-norm} 
\mathbf{0} \leq \boldsymbol \psi_\infty \leq \mathbf{1}\,  \text { in } \, \mathbb{R} \quad\hbox{and}\quad
\psi_{1,\infty}(0)+\psi_{2,\infty}(0)=1.
\end{equation}

Now, by the condition (A2), according to phase diagrams of system \eqref{eq: limits v_i} and Lemma \ref{separatrix}, 
$\mathbf{v}_{\infty}(t)$ must be either  one of the four equilibrium 
$\mathbf{0}$, $\mathbf{1}$, $\mathbf{e}_0$, $\textbf{e}^*_{s_\infty}$ (see \eqref{4s-equilibrium}) or 
one of the following heteroclinic orbits:
\begin{equation*}
\left.\begin{array}{l}
\Gamma_1=\left\{ \mathbf{v}_{\infty}\in C^1(\R;[0,1]^2): \, \lim_{t\to -\infty}\mathbf{v}_{\infty}(t)=\mathbf{e}_0,\, \,  \lim_{t\to +\infty}\mathbf{v}_{\infty}(t)=\mathbf{0}\right\},\vspace{5pt} \\
\Gamma_2=\left\{ \mathbf{v}_{\infty}\in C^1(\R;[0,1]^2): \, \lim_{t\to -\infty}\mathbf{v}_{\infty}(t)=\mathbf{e}_0,\, \,  \lim_{t\to +\infty}\mathbf{v}_{\infty}(t)=\mathbf{1}\right\},\vspace{5pt} \\
\Gamma_3=\left\{ \mathbf{v}_{\infty}\in C^1(\R;[0,1]^2): \, \lim_{t\to -\infty}\mathbf{v}_{\infty}(t)=\textbf{e}^*_{s_\infty},\, \,  \lim_{t\to +\infty}\mathbf{v}_{\infty}(t)=\mathbf{0}\right\},\vspace{5pt} \\
\Gamma_4=\left\{ \mathbf{v}_{\infty}\in C^1(\R;[0,1]^2): \, \lim_{t\to -\infty}\mathbf{v}_{\infty}(t)=\textbf{e}^*_{s_\infty},\, \,  \lim_{t\to +\infty}\mathbf{v}_{\infty}(t)=\mathbf{1}\right\},\vspace{5pt} \\
\Gamma_5=\left\{ \mathbf{v}_{\infty}\in C^1(\R;[0,1]^2): \, \lim_{t\to -\infty}\mathbf{v}_{\infty}(t)=\mathbf{e}_0,\, \,  \lim_{t\to +\infty}\mathbf{v}_{\infty}(t)=\textbf{e}^*_{s_\infty}\right\}.
		\end{array}\right.
	\end{equation*}
We will classify the above possibilities into four cases, and derive a contradiction in each case. 
	
	\textbf{Case (a)}:  $\mathbf{v}_{\infty}\in \{\mathbf{0},  \mathbf{1},\Gamma_1,\Gamma_2,\Gamma_3,\Gamma_4 \}$.
	
	In this case,  for any small $\epsilon>0$, there always exists a large $t'>0$ (if $\mathbf{v}_{\infty}\equiv \mathbf{0}$ or $\mathbf{v}_{\infty}\equiv \mathbf{1}$, then $t'$ can be an arbitrary time) such that either
	$\mathbf{0}\leq \mathbf{v}_{\infty}(t') \leq (\epsilon,\epsilon)$ or $(1-\epsilon,1-\epsilon)\leq \mathbf{v}_{\infty}(t') \leq \mathbf{1}$. 
    This together with \eqref{converge-vn-phsn} implies that for all large $n\in\N$, there holds  
	$$ \mathbf{0}\leq \bar{\mathbf{v}}_n(t_n+t') \leq (2\epsilon,2\epsilon) \quad\hbox{or}\quad  (1-2\epsilon, 1-2\epsilon)\leq \bar{\mathbf{v}}_n (t_n+t') \leq \mathbf{1}. $$
	 However, by Lemma \ref{linearly stable claim large T} and its subsequent discussion, one can find a small constant $\tilde{\epsilon}_1>0$ such that for each large $n\in\N$, 
	 $$\bar{\mathbf{v}}_n(t)\not \in \left\{(x_1,x_2): 0\leq x_i\leq \tilde{\epsilon}_1  \hbox{ for } i=1,2\right\}\cup \left\{(x_1,x_2): 1-\tilde{\epsilon}_1\leq x_i\leq 1  \hbox{ for } i=1,2\right\}  $$
	 for all $t\in\R$.  Therefore,  Case (a) is ruled out.

	\textbf{Case (b)}:  $\mathbf{v}_{\infty}(\cdot) \equiv \mathbf{e}_0$.
	
	By \eqref{psiinfty-norm}, $\boldsymbol \psi_\infty$ cannot be identically equal to $\mathbf{0}$, and hence, the uniqueness of solutions to the linear system \eqref{eq: limits w_i} implies that either $\psi_{1,\infty}(\cdot)>0$ in $\R$ or $\psi_{2,\infty}(\cdot)>0$ in $\R$. 
	Assuming without loss of generality that the former case happens, and considering the first component of  \eqref{eq: limits w_i} with $\mathbf{v}_{\infty}(\cdot)\equiv \textbf{e}_0$, we obtain that 
	$\psi'_{1,\infty}=(r_1(s_\infty)+\tilde{\lambda}_\infty))\psi_{1,\infty}$, namely, $\psi_{1,\infty}(t)=\psi_{1,\infty}(0){\rm exp}\{(r_1(s_\infty)+\tilde{\lambda}_\infty)t\}$. Since $\tilde{\lambda}_\infty\geq 0$, this contradicts the boundedness of $\psi_{1,\infty}(t)$ (see \eqref{psiinfty-norm}). In the other case that $\psi_{2,\infty}(\cdot)>0$ in $\R$, one can derive a similar contradiction by considering the second component of \eqref{eq: limits w_i}. Hence, Case (b) cannot happen.

\textbf{Case (c)}:  $\mathbf{v}_{\infty}(\cdot) \equiv \mathbf{e}_{s_\infty}^*$.

In this case, $A(\mathbf{v}_{\infty})$ is constant matrix with negative off-diagonal entries.  We rewrite this matrix by $A_*$. One readily sees that
\bqq\label{matrix-A*}
A_*= 
\left(\begin{array}{cc}
	-r_1(s_\infty) v_{1,s_\infty}^* &  k_1(s_\infty)p_2(s_\infty)v_{1,s_\infty}^*    \vspace{5pt} \\
	k_2(s_\infty)p_1(s_\infty)\left( 1-v_{2,s_\infty}^*\right)  &  -r_2(s_\infty)\left( 1-v_{2,s_\infty}^*\right) 
\end{array}\right).
\eqq
By condition (A2), it is straightforward to check that ${\rm det}(A_*)<0$. 
This implies that  $A_*$ has a positive eigenvalue and a negative eigenvalue, denote by $\lambda^+$ and $\lambda^-$, respectively. Let 
$\boldsymbol \psi^+=(\psi^+_1,\psi^+_2)$ and $\boldsymbol \psi^-=(\psi^-_1,\psi^-_2)$ be the corresponding eigenvectors. It is clear that $\psi^+_1\psi^+_2>0$ and $\psi^+_1\psi^+_2<0$. Up to multiplication, we may assume without loss of generality that
\begin{equation}\label{mul-eigenfunction}
	\psi^+_1>0,\quad \psi^+_2>0,\quad \psi^-_1>0,\quad \psi^-_2<0 \quad\hbox{and}\quad \| \boldsymbol \psi^{\pm} \|=1, 
\end{equation}  
where $\|\cdot\|$ denotes the $L^{\infty}$-norm in $\R^2$. Furthermore, it follows from \eqref{eq: limits w_i} that $\mathbf{w}(t):={\rm exp}(-\tilde{\lambda}_\infty t)\boldsymbol \psi_\infty(t)$ is an entire solution of the linear problem $d\mathbf{w}/dt=A_*\mathbf{w}$ in $\R$, and hence, by basic theory for linear ODEs with constant coefficients,  we have
$$\boldsymbol \psi_\infty(t)=h_1 \me^{\left(\tilde{\lambda}_\infty+\lambda^+\right)t}\boldsymbol \psi^+ + h_2 \me^{\left(\tilde{\lambda}_\infty+\lambda^-\right)t}\boldsymbol \psi^- \,\,\hbox{ for }\, t\in\R,  $$
where $h_1$, $h_2$ are two constants in $\R$. Since $\boldsymbol \psi_\infty(\cdot)$ is a nonnegative and nonzero function (see \eqref{psiinfty-norm}), due to \eqref{mul-eigenfunction}, there must hold $h_1>0$. Remember that $\tilde{\lambda}_\infty\geq 0$ and $\lambda^+>0>\lambda^+$. This implies that  $\| \boldsymbol \psi_\infty(t)\| \to +\infty$ as $t\to+\infty$, which contradicts the boundedness of $ \boldsymbol \psi_\infty(\cdot)$ again. Hence, Case (c) is ruled out, too.

It remains to exclude the case where $\mathbf{v}_{\infty}(t)$ is a heteroclinic orbit connecting the unstable equilibrium $\mathbf{e}_0$ to the saddle equilibrium $\mathbf{e}^*_{s_{\infty}}$, that is,

\textbf{Case (d)}:  $\mathbf{v}_{\infty}(\cdot)\in  \Gamma_5$.   

	In this case, it follows from Lemma \ref{separatrix} that 
	\begin{equation}\label{vinfty-prime}
		v'_{1,\infty}(\cdot)>0, \, \, v'_{2,\infty}(\cdot)<0\,\, \hbox{ in }\,\,\R, \quad\hbox{and}\quad \|\mathbf{v}'_{\infty}(t)\| \to 0\,\,\hbox{ as }\,\, t\rightarrow+\infty.
	\end{equation}	
 We further observe that 
    \begin{equation}\label{positive-psi-infty}
    \psi_{1,\infty}(\cdot)>0 \quad \hbox{and}\quad \psi_{2,\infty}(\cdot)>0 \,\,\hbox{ in } \,\, \R. 
    \end{equation}
    Indeed, if this is not true, then either $\psi_{1,\infty}(\cdot)>0$ and $\psi_{2,\infty}(\cdot)\equiv 0$ in $\R$ or $\psi_{1,\infty}(\cdot)\equiv 0$ and $\psi_{2,\infty}(\cdot)> 0$ in $\R$. Suppose that the former case happens, and consider the second component of system \eqref{eq: limits w_i}. One obtains that 
    $k_2(s_\infty)p_1(s_\infty)(1-v_{2,\infty}(t))\psi_{1,\infty}(t)\equiv 0$,  which is impossible. Proceeding similarly, one can derive a contradiction in the latter case. Therefore, \eqref{positive-psi-infty} is proved. 
    
    Next, we claim that  
	\begin{equation}\label{limit-psi-infty}
		\|\boldsymbol \psi_\infty(t)\|\rightarrow  0\,\,\hbox{ as }\,\, t\to+\infty.
	\end{equation}
 Assume by contradiction that there exists a sequence $(\tau_n)_{n\in \N} \subset (0,+\infty)$ such that as $n\to+\infty$,  
 $\tau_n\to+\infty$, and $\boldsymbol\psi_\infty(\tau_n)$ converges to a nonnegative and nonzero vector in $\R^2$. Then, by the Arzela-Ascoli theorem, one can find some function $\hat{\boldsymbol\psi}(\cdot) \in C^1(\R;\R^2)$ such that,  up to extraction of some subsequence, $\boldsymbol\psi_\infty(\cdot+\tau_n) \to \hat{\boldsymbol\psi}(\cdot) $ in $C^1_{loc}(\R)$ as $n\to+\infty$. 
 It is clear that  $\hat{\boldsymbol\psi}(\cdot)$ is nonnegatie and nonzero. Furthermore, since $\mathbf{v}_{\infty}(t)\to \textbf{e}^*_{s_\infty}$ as $t\to+\infty$, one readily checks that 
 $\hat{\boldsymbol\psi}(\cdot)$ is an entire solution of system \eqref{eq: limits w_i} with the matrix coefficients $A(\mathbf{v}_{\infty})$ replaced by $A_*$, where $A_*$ is the constant matrix given by \eqref{matrix-A*}. Then, the same reasoning as in excluding Case (c) shows that $\hat{\boldsymbol\psi}(\cdot)$ cannot be a bounded vector function. This is impossible, since 
 \eqref{psiinfty-norm} implies that $\mathbf{0} \leq \hat{\boldsymbol\psi}(\cdot)\leq \mathbf{1}$ in $\R$. Therefore, our claim \eqref{limit-psi-infty} is proved. 
 
We are now ready to derive a contradiction in Case (d). Recall that $\lambda^+>0>\lambda^-$ are the two eigenvalues of the matrix $A_*$, and that  $\boldsymbol \psi^+=(\psi^+_1,\psi^+_2)$ and $\boldsymbol \psi^-=(\psi^-_1,\psi^-_2)$ are the corresponding eigenvectors satisfying \eqref{mul-eigenfunction}. We choose a positive time $t^*>0$ such that 
\begin{equation}\label{choose-t*}
\me^{\lambda^+t^*}\min\left\{ \psi_1^+,\,\psi_2^+\right\} \geq 3.
\end{equation}
Let $\varepsilon_0\in (0,1)$ be a small constant such that 
\begin{equation}\label{choose-vare0}
\varepsilon_0  t^*  \me^{(\lambda^++\|A_*\|+1)t^*} \leq 1,   
\end{equation}
where $\|A_*\|$ is the $L^{\infty}$-norm of the matrix $A_*$. Since $\mathbf{v}_{\infty}(t)\to \textbf{e}^*_{s_\infty}$ as $t\to+\infty$,
and since $\mathbf{v}'_{\infty}(\cdot)$ and $\boldsymbol \psi_\infty(\cdot)$ satisfy \eqref{vinfty-prime} and \eqref{limit-psi-infty}, respectively, one can find some large $\tau^*>0$ such that 
\begin{equation}\label{AT-A-estimate}
	\| A(\mathbf{v}_{\infty}(t+\tau^*))-A_* \| \leq \varepsilon_0 \,\, \text{ for all }\, t\geq 0,
\end{equation}
and that 
\bqq\label{psi + v}
|\psi_{i,\infty}(t+\tau^*)| \leq | \psi_{i,\infty}(\tau^*)|\quad\hbox{and}\quad |v'_{i,\infty}(t+\tau^*)|\leq  |v'_{i,\infty}(\tau^*)|\,\, \text{ for all }\, t\geq 0,\,i=1,2.
\eqq	 

Now, let $\mathbf{w}(t):=(w_1,w_2)(t)$ and $\tilde{\mathbf{w}}(t):=(\tilde{w}_1,\tilde{w}_2)(t)$ be, respectively, the solutions of the following two Cauchy problems
\begin{equation}\label{eq-w-tildew}
		\left\{\begin{array}{l}
		\displaystyle \frac{d \mathbf{w}}{dt}= A_* \mathbf{w}, \quad 0<t\leq t^*, \vspace{5pt} \\
		\mathbf{w}(0)=	\boldsymbol \psi^{+},
	\end{array}\right.
\quad\hbox{and}\quad
\left\{\begin{array}{l}
\displaystyle \frac{d \tilde{\mathbf{w}}}{dt}= A(\mathbf{v}_{\infty}(t+\tau^*)) \tilde{\mathbf{w}}, \quad 0<t\leq t^*, \vspace{5pt} \\
\tilde{\mathbf{w}}(0)=\boldsymbol \psi^{+}. 
\end{array}\right.
\end{equation}
It is clear that $\mathbf{w}(t)=\me^{\lambda^+t}\boldsymbol\psi^{+}$ for $0\leq t\leq t^*$. Moreover, due to  \eqref{vinfty-prime} and \eqref{positive-psi-infty}, $\mathbf{v}_{\infty}'(t+\tau^*)$ and $\me^{-\tilde{\lambda}_\infty t}\boldsymbol \psi_\infty(t+\tau^*)$ are linearly independent, and they form a fundamental set of solutions to the 
 linear system $d\mathbf{z}/dt= A(\mathbf{v}_{\infty}(t+\tau^*)) \mathbf{z}$ for $t\in\R$. This immediately gives that
$$\tilde{\mathbf{w}}(t)=\tilde{h}_1 \me^{-\tilde{\lambda}_\infty t}\boldsymbol \psi_\infty(t+\tau^*) +\tilde{h}_2 \mathbf{v}_{\infty}'(t+\tau^*) \,\,\hbox{ for }\, 0\leq t\leq t^*,  $$
where $\tilde{h}_1$ and $\tilde{h}_2$ (which may depend on $\tau^*$) are two constants in $\R$. 
Furthermore, by using \eqref{mul-eigenfunction}, \eqref{psi + v} and the assumption that $\tilde{\lambda}_{\infty}\geq 0$, one can conclude that if $\tilde{h}_1\tilde{h}_2\geq 0$, then $|\tilde{w}_1(t)| \leq |\tilde{w}_1(0)| \leq \|\boldsymbol \psi^{+} \|=1 $, whence
$$\|\mathbf{w}(t)-\tilde{\mathbf{w}}(t)\| \geq \left| e^{ \lambda^+ t} \psi^+_1-\tilde{w}_1(t) \right| \geq e^{ \lambda^+ t} \psi^+_1-1; $$
and if $\tilde{h}_1\tilde{h}_2\leq 0$, then $|\tilde{w}_2(t)| \leq |\tilde{w}_2(0)| \leq \|\boldsymbol \psi^{+} \|=1 $, whence
$$\|\mathbf{w}(t)-\tilde{\mathbf{w}}(t)\| \geq \left| e^{ \lambda^+ t} \psi^+_2-\tilde{w}_2(t)\right| \geq e^{ \lambda^+ t} \psi^+_2-1.$$
Combining the above and choosing $t=t^*$,  we obtain that 
\begin{equation}\label{w-tildew-low}
\|\mathbf{w}(t^*)-\tilde{\mathbf{w}}(t^*)\| \geq  e^{ \lambda^+ t^*} \min\left\{\psi^+_1, \, \psi^+_2 \right\}-1\geq 2,
\end{equation}
where the last inequality follows immediately from \eqref{choose-t*}.

On the other hand, one readily checks from \eqref{eq-w-tildew} and \eqref{AT-A-estimate} that for any $0<t\leq t^*$, 
$$
\left.\begin{array}{ll}
\displaystyle  \left\| \mathbf{w}(t)-\tilde{\mathbf{w}}(t)\right\| \!\!\! & =  \displaystyle \left\|\int_{0}^{t} A(\mathbf{v}_{\infty}(\tau+\tau^*))\left( \mathbf{w}(\tau)-\tilde{\mathbf{w}}(\tau)\right) d\tau-\int_{0}^{t} \left(A(\mathbf{v}_{\infty}(\tau+\tau^*))-A_*\right)\mathbf{w}(\tau) d\tau \right\| \vspace{5pt} \\
 & \displaystyle \leq   \int_{0}^{t} \left\| A(\mathbf{v}_{\infty}(\tau+\tau^*))\right\| 
 \left\| \mathbf{w}(\tau)-\tilde{\mathbf{w}}(\tau)\right\| d\tau +  \varepsilon_0 t^* \me^{\lambda^+t^*}. 
\end{array}\right.$$
It then follows from the Gronwall inequality that 
$$\|\mathbf{w}(t)-\tilde{\mathbf{w}}(t) \| \leq   \varepsilon_0 t^*\me^{\lambda^+t^*} {\rm exp}\left(\int_{0}^{t} \left\| A(\mathbf{v}_{\infty}(\tau+\tau^*))\right\|  d\tau\right) 
 $$ 
for all $0<t\leq t^*$. In particular, choosing $t=t^*$, and using \eqref{choose-vare0} and \eqref{AT-A-estimate}, we obtain that
$$\|\mathbf{w}(t^*)-\tilde{\mathbf{w}}(t^*) \| \leq   \varepsilon_0  t^*  \me^{(\lambda^++\|A_*\|+1)t^* } \leq 1,
$$
which is a contradiction with the estimate in \eqref{w-tildew-low}. Consequently,  Case (d) is ruled out, and the proof of Lemma \ref{linearly unstable claim large T} is complete.		 	 	
\end{proof}

Finally, we take $T^*=\max\{\tilde{T}_1,\tilde{T}_2,\tilde{T}_3\}$. Based on Lemmas \ref{linearly stable claim large T} and \ref{linearly unstable claim large T}, Theorem \ref{exitstence wave large} is an easy consequence of Proposition \ref{existence-fix}. 
%%%%%%%%%%%%%%%%%%%%%%%%%%%%%%%%%%%%%%%%%%%%%%%%%%%%%%%%%%%%%%%%%%%%%%%%%%%%%%%%%%%%%%%%%%%%%%%%%%%%%%%%%%%%%%%%%%%

\subsection{Convergence of $c_T$ as $T\rightarrow+\infty$ }

We are now in a position to study the asymptotic behavior of the wave speed as $T\to+\infty$. Let $T^*>0$ be the constant provided by Theorem \ref{exitstence wave large}, and for each $T>T^*$, let $(\varphi_{1,T},\varphi_{2,T})( x-c_Tt,t)$ be the (unique) periodic traveling wave of \eqref{eq:change of variable} connecting $\mathbf{0}$ and $\mathbf{1}$ and satisfying the normalization \eqref{normalize small T wave}. 
In this section, we prove the following theorem.

\begin{theo}\label{converge-largeT}
	Let  {\rm (A2)} hold. Then the estimate \eqref{large speed cover} holds true.
\end{theo}

 We prove Theorem \ref{converge-largeT} by a sub- and super-solution method similar in spirit to that used in showing Theorem \ref{converge-samllT}. However, unlike the single-wave perturbation for rapid oscillations, the slowly varying case necessitates perturbations of the following family of homogeneous waves: 
 \bqq\label{eq: frozen-limit-change}
 \left\{\begin{array}{l}
 	d_1(s) \partial_{\xi\xi} \varphi_1(\xi;s)+c(s)\partial_\xi\varphi_1(\xi;s)+g_{1}(s,\varphi_1(\xi;s),\varphi_2(\xi;s))=0  \quad \hbox{in }\,\, \R, \vspace{5pt}\\  
 	d_2(s)   \partial_{\xi\xi} \varphi_2(\xi;s)+c(s)\partial_\xi\varphi_2(\xi;s)+g_{2}(s,\varphi_1(\xi;s),\varphi_2(\xi;s))=0 \quad \hbox{in }\,\, \R, \vspace{5pt}\\  
 	(\varphi_1,\varphi_2)(-\infty;s)=\mathbf{1} \quad\hbox{and}\quad(\varphi_1,\varphi_2)(+\infty;s)=\mathbf{0}, 
 \end{array}\right.
 \eqq
where $s$ is an arbitrary point in $\R$, $c(s)$ is the wave speed of the traveling wave for system \eqref{eq:frozen-wave} with coefficients frozen at the point $s$, and the pair of functions $(g_1,g_2)$ is defined in \eqref{g-frozen}. For each $s\in\R$, the existence of solution $\boldsymbol \varphi(\xi;s):=(\varphi_1,\varphi_2)(\xi;s)$ to system \eqref{eq: frozen-limit-change} is obvious, since it is directly transformed from the solution of \eqref{eq:frozen-wave} by  
taking the change of variable $(\varphi_1,\varphi_2)(\xi;s):=(\phi_1(\xi;s)/p_1(s), 1-\phi_2(\xi;s)/p_2(s))$. In the sequel, for definiteness, we normalized each $\boldsymbol \varphi(\xi;s)$ as follows
\bqq\label{normalize fix s wave}
\varphi_1(0;s)=\gamma_1 \,\, \hbox{ for each }\,\, s\in\R, 
\eqq
where $\gamma_1$ is a constant such that $0<\gamma_1<\min_{s\in\R}\left\lbrace  v_{1,s}^*,v_{2,s}^* \right\rbrace$, and $(v_{1,s}^*,v_{2,s}^*) $ is the unique positive saddle equilibrium of system \eqref{eq:fix-s-Lv} (note that $v_{1,s}^*\in (0,1)$ is continuous and $1$-periodic in $s$).  
Clearly, under the above normalization, for any $s\in\R$, the wave $\boldsymbol \varphi(\xi;s)$ is uniquely determined. It is also easily seen that the function $s\mapsto c(s)$ is $1$-periodic, and for each $i=1,2$, the function $(\xi,s)\mapsto \varphi_i(\xi,s) \in (0,1)$ is $1$-periodic in $s$ and decreasing in $\xi$.  

Before presenting our sub- and super-solutions, let us first summarize some uniform estimates of $\boldsymbol \varphi(\cdot;s)$ with respect to $s\in\R$. 

\begin{prop}\label{fix s continue}
	Let  {\rm (A2)}  hold. Then, for each $i=1,2$, the following statements hold true:
	\begin{itemize}
		\item [{\rm (i)}]  $\lim _{\xi \rightarrow+\infty} \varphi_i(\xi; s)=0$ and $\lim _{\xi \rightarrow-\infty} \varphi_i(\xi; s)=1$ uniformly in $s \in \mathbb{R}$; 
		
		\item [{\rm (ii)}] for any $\delta \in(0,1 / 2]$, there exists $\beta=\beta(\delta)>0$ $($independent of $s$$)$ such that
		$$
		\partial_{\xi} \varphi_i(\xi; s) \leq-\beta \,\,\text{ for all }\,\,(\xi, s) \in \mathbb{R}^2\,\, \text{ such that }\,\, \delta \leq \varphi_i(\xi; s) \leq 1-\delta;
		$$
		
		\item [{\rm (iii)}]  the function $s \mapsto c(s)$ is of class $C^1(\mathbb{R})$, and the function $(\xi, s) \mapsto \varphi_i(\xi; s)$ is of class $C_{\xi ; s}^{2 ; 1}\left(\mathbb{R}^2\right)$ and satisfies
		\begin{equation}\label{partial_s phi uniform bound} 
			\sup _{\xi \in \mathbb{R}, s \in \mathbb{R}}\left|\partial_s \varphi_i(\xi; s)\right|<+\infty.
		\end{equation}
	\end{itemize}
\end{prop}

The proof of Proposition \ref{fix s continue} follows from similar arguments to those used in showing \cite[Proposition 4.1]{dhl} for traveling waves of scalar bistable equations. We will outline the key steps in the appendix.

For each $T>T^*$, we define
$$
X_T(t)=\int_0^t c\left(\frac{\tau}{T}\right) d \tau \,\, \text { for }\,\, t \in \mathbb{R}.
$$
Clearly, $X_T(0)=0$, and since $c(\cdot)$ is of class $C^1(\mathbb{R})$, the function $t \mapsto X_T(t)$ is at least of class $C^2(\mathbb{R})$. Moreover, since $c(\cdot)$ is $1$-periodic, it is also easily seen that
$$
\lim _{t \rightarrow +\infty} \frac{X_T(t)}{t}=\lim _{t \rightarrow +\infty} \frac{T}{t} \int_0^{t / T} c(s) d s=c^*,
$$
where $c^*$ is the constant defined in \eqref{speed-limit-large}.

Recall that $\theta_+$, $\gamma_+$ and $\gamma_0$ are the positive constants provided by \eqref{bound-rak}, \eqref{bound-pT} and \eqref{gamma_0}, respectively, and that $\rho \in C^2(\R)$ is a nonnegative function satisfying \eqref{function-rho}. It is clear that $0< p_i(s)\leq \gamma_+$ for all $s\in\R$, $i=1,2$.

The following lemma gives a super-solution of system \eqref{eq:change of variable} when $T$ is large. 

\begin{lem}\label{large T supsub}
There exists $\varepsilon_0>0$ such that for any $\varepsilon\in (0,\varepsilon_0)$, there exists $\tilde{T}_4=\tilde{T}_4(\varepsilon)>0$ such that for any large $T > \tilde{T}_4$, the pair of functions $\mathbf{w}_T^+(t,x)=(w_{1,T}^+,w_{2,T}^+)(t,x)$ defined by 
\begin{equation*}
\left\{\begin{array}{l}
\displaystyle w_{1,T}^{+}(t, x)=\varphi_1\left(\xi_T; t/T\right)+q_T(t)\left(\rho(\xi_T)+\frac{\gamma_0}{8\theta_+\gamma_+}(1-\rho(\xi_T))\right)  \vspace{5pt}\\  
\displaystyle w_{2,T}^{+}(t, x)=\varphi_2\left(\xi_T ; t/T\right)+q_T(t)\left(\frac{\gamma_0}{8\theta_+\gamma_+}\rho(\xi_T)+(1-\rho(\xi_T))\right)
\end{array}\right.
\end{equation*}
is a super-solution of \eqref{eq:change of variable} for $(t,x) \in (0,+\infty) \times \R$, where $\xi_T=x-X_T(t)+\kappa_T(t)$, and $q_T(\cdot)$, $\kappa_T(\cdot)$ are $C^1([0,+\infty))$ functions satisfying
\bqq\label{defi q_T large T}
\left\{\begin{array}{l}
q_T(0)=\varepsilon, \quad q_T'(t)<0<q_T(t) \,\, \text{ for } \,\,t> 0, \vspace{5pt}\\  
\kappa_T(0)=0, \quad \kappa_T'(t)<0\,\, \text{ for }\,\, t> 0. \end{array}\right.
\eqq
\end{lem}

\begin{proof}
Denote $C_1=\gamma_0/(8\theta_+\gamma_+)$, and let $\sigma_0$ be a small constant given by
	\bqq\label{sigma0-largeT}
	\sigma_0=\min \left\{\frac{1}{4},\,\frac{\gamma_0}{4\theta_+\gamma_+},\, \frac{\gamma_0C_1}{2\theta_+\gamma_+}\right\}.
	\eqq
	By Proposition \ref{fix s continue} (i) and the fact that, for each $i=1,2$, the function  $\varphi_i(\xi;s)$ is continuous in $\R^2$, decreasing in $\xi$ and periodic in $s$, there exists $M>0$ sufficiently large such that 
\bqq\label{wave fix s +-}
\left\{\begin{array}{ll}
\displaystyle 0< \varphi_i(\xi;s)<\frac{\sigma_0}{2} & \text {for all }\,\, \xi \geq M,\,s\in\R, \vspace{5pt}\\  
\displaystyle \frac{\sigma_0}{4} \leq  \varphi_i(\xi;s) \leq 1-\frac{\sigma_0}{4} & \text{for all }\,\, -M \leq \xi \leq M,\,\, s \in \mathbb{R},  \vspace{5pt}\\  
\displaystyle  1-\frac{\sigma_0}{2}< \varphi_i(\xi;s)<1 & \text {for all } \,\, \xi \leq-M,\,s\in\R. 
\end{array}\right.
\eqq
It then follows from Proposition \ref{fix s continue} (ii) that there exists $\beta_0>0$ such that 
\bqq\label{beta_2}
\partial_{\xi} \varphi_i(\xi;s) \leq -\beta_0\,\,\hbox{ for all }\,\, -M\leq \xi\leq M,\,s\in\R.
\eqq 
We now set
	\bqq\label{delta_2}
	\varepsilon_0= \min \left\{\frac{\sigma_0}{2},\, \frac{\sigma_0}{2(C_1+1)},\,\frac{\beta_0}{C_1+1}\right\}.
	\eqq
Similarly to the strategy of the proof of Lemma \ref{small T supsub}, it suffices to prove that for any $\varepsilon\in (0,\varepsilon_0)$ and any large $T>0$, there exist $C^1([0,+\infty))$ functions $q_T(\cdot)$ and $\kappa_T(\cdot)$ satisfying \eqref{defi q_T large T} such that 
$$
\left\{\begin{array}{l}
	\mathcal{L}_1 (\textbf{w}^{+}_{T}):=\partial_t w^{+}_{1,T}-d_{1,T}(t)\partial_{x x} w^{+}_{1,T}-g_{1,T}\left(t, \textbf{w}^{+}_{T}\right)\geq 0, \vspace{5pt}\\ 
	\mathcal{L}_2 (\textbf{w}^{+}_{T}):=\partial_t w^{+}_{2,T}-d_{2,T}(t)\partial_{x x} w^{+}_{2,T}-g_{2,T}\left(t, \textbf{w}^{+}_{T}\right)\geq 0, 
\end{array}\right.$$
for $(t, x) \in (0,+\infty) \times \mathbb{R}$. 
We will determine the functions $q_T(\cdot)$ and $\kappa_T(\cdot)$, restricting of our verification to $\mathcal{L}_1 (\textbf{w}^{+}_{T})\geq 0$, since the verification of $\mathcal{L}_2 (\textbf{w}^{+}_{T})\geq 0$ is parallel. 

Since for each $s \in \mathbb{R}$, $(\varphi_1,\varphi_2)(\cdot;s)$ is an entire solution of \eqref{eq: frozen-limit-change}, it is straightforward to check that for $(t, x) \in(0, +\infty) \times \mathbb{R}$, 
\begin{equation*}
\begin{aligned}
	\mathcal{L}_1 (\textbf{w}^{+}_{T}) =&\kappa_T^{\prime} \partial_{\xi} \varphi_1+\frac{1}{T} \partial_s \varphi_1+q_T^{\prime}\left(\rho+(1-\rho)C_1\right)+q_T\rho'\left(\kappa_T^{\prime}-c\right)\left(1-C_1\right) \vspace{5pt}\\ 
	&\qquad -d_{1}q_T\rho''(1-C_1)-a_{1}p_{1,T}w_{1,T}^+(1-w_{1,T}^+)+a_{1}p_{1}\varphi_1(1-\varphi_1)  \vspace{5pt}\\ 
	&\qquad +k_{1}p_{2,T}w_{1,T}^+(1-w_{2,T}^+)-k_{1}p_{2}\varphi_1(1-\varphi_2),
\end{aligned}
\end{equation*}
where $\rho(\cdot)$, $\rho^{\prime}(\cdot)$, $\rho^{\prime \prime}(\cdot)$ are evaluated at $\xi_T$; $a_1(\cdot)$, $d_1(\cdot)$, $k_1(\cdot)$, $p_i(\cdot)$, $c(\cdot)$ are evaluated at $t/T$;  $p_{i,T}(\cdot)$, $q_T(\cdot)$, $q_T'(\cdot)$, $\kappa_T(\cdot)$, $\kappa_T'(\cdot)$ are evaluated at $t$; and $\varphi_i(\cdot;\cdot)$ together with its partial derivatives are evaluated at $(\xi_T;t/T)$.
By Lemma \ref{uniform cover s}, there exists a constant $C_2>0$ independent of $T$ such that 
$|p_{i,T}(t)-p_i(t/T)| \leq C_2/T$ for all $t>0$, $i=1,2$. Moreover, 
since $q_T(\cdot)$ is required to satisfy $0\leq q_T(t) \leq \varepsilon\leq \varepsilon_0$, and since $\rho(\cdot)$ satisfies \eqref{function-rho},  it follows easily from \eqref{delta_2} that $0<w_{i,T}^+(t,x)\leq \varphi_{1}(\xi_T;t/T)+\varepsilon_0(1+C_1) \leq 2$ for all $(t,x)\in (0,+\infty)\times\R$, whence we have
$|w_{1,T}^+(1-w_{1,T}^+)|\leq 2$, $|w_{1,T}^+(1-w_{2,T}^+)|\leq 2$. 
It then follows that 
$$\left\{\begin{array}{l}
\displaystyle \left|\left(a_{1}(t/T)p_{1,T}(t)-a_1(t/T)p_1(t/T)\right)w_{1,T}^+(1-w_{1,T}^+)   \right| \leq \frac{2\theta_+C_2}{T}, \vspace{5pt}\\ 
\displaystyle \left|\left(k_{1}(t/T)p_{2,T}(t)-k_1(t/T)p_2(t/T)\right)w_{1,T}^+(1-w_{2,T}^+)   \right| \leq \frac{2\theta_+C_2}{T},
\end{array}\right.$$
for all $(t,x)\in (0,+\infty)\times\R$. On the other hand, Proposition \ref{fix s continue} (iii) implies that, for each $i=1,2$,  
\begin{equation*}
	\sup _{\xi \in \mathbb{R}, s \in \mathbb{R}}\left|\partial_s \varphi_i(\xi; s)\right|\leq C_3
\end{equation*}
for some constant $C_3>0$. 
Combining the above, we obtain that
\begin{equation}\label{largeT-Lw}
	\begin{aligned}
		\mathcal{L}_1 (\textbf{w}^{+}_{T}) \geq & -\frac{C_4}{T}+ \kappa_T^{\prime} \partial_{\xi} \varphi_1+q_T^{\prime}\left(\rho+(1-\rho)C_1\right)+q_T\rho'\left(\kappa_T^{\prime}-c\right)\left(1-C_1\right) \vspace{5pt}\\ 
		&\qquad -d_{1}q_T\rho''(1-C_1)-a_{1}p_{1}w_{1,T}^+(1-w_{1,T}^+)+a_{1}p_{1}\varphi_1(1-\varphi_1)  \vspace{5pt}\\ 
		&\qquad +k_{1}p_{2}w_{1,T}^+(1-w_{2,T}^+)-k_{1}p_{2}\varphi_1(1-\varphi_2),
	\end{aligned}
\end{equation}
for all $(t,x)\in (0,+\infty)\times\R$, where $C_4=4\theta_+C_2+C_3$.

Below, we complete the proof by considering three cases.

\textbf{Case (a)}: $(t,x)\in (0,+\infty)\times \R$ such that $\xi_T\geq M$. 

In this case, we have $\rho= \rho^{\prime}=\rho^{\prime \prime}=0$. Since $\kappa_T(\cdot)$ is required to be decreasing, and since $\varphi_1(\xi;s)$ is decreasing in $\xi$, we have  $\kappa_T^{\prime} \partial_{\xi} \varphi_1\geq 0$. It then follows from \eqref{largeT-Lw} that 
\begin{equation}\label{largeT-LgeqM}
\left.\begin{array}{ll}
\displaystyle\mathcal{L}_1 (\textbf{w}^{+}_{T})\geq -\frac{C_4}{T} \!\!\!& \displaystyle +C_1q_T^{\prime}+C_1q_T\left( -a_{1}p_{1}\left(1-(w_{1,T}^++\varphi_1)\right)+k_{1}p_{2}\right) \vspace{5pt}\\ 
&\displaystyle -k_{1}p_{2}\left(\varphi_1q_T+C_1\varphi_2q_T+C_1q_T^2\right). 
\end{array}\right.\end{equation}
Notice that, by the first line of \eqref{wave fix s +-}, there holds 
$0<1-(w_{1,T}^++\varphi_1)<1$. This together with the first inequality of \eqref{gamma_0} implies that 
$$ C_1q_T\left(-a_{1}p_{1}\left(1-(w_{1,T}^++\varphi_1)\right)+k_{1}p_{2}\right)\geq  C_1q_T\left( -a_{1}p_{1}+k_{1}p_{2} \right) \geq \gamma_0C_1q_T. $$
Moreover, since 
$0\leq\varphi_2+q_T\leq \sigma_0/2+ \epsilon_0\leq \sigma_0$, by using \eqref{sigma0-largeT}, we have
$$ -k_{1}p_{2}\left(\varphi_2+q_T\right)C_1q_T \geq -\theta_+\gamma_+\sigma_0C_1q_T\geq -\frac{\gamma_0}{4}C_1q_T
$$
and
$$-k_{1}p_{2}\varphi_1q_T \geq -\theta_+\gamma_+ \frac{\sigma_0}{2}q_T\geq -\frac{\gamma_0}{4}C_1q_T.$$
Consequently, it is readily checked from \eqref{largeT-LgeqM} that 
\begin{equation}\label{largeT-geqM}
\mathcal{L}_1 (\textbf{w}^{+}_{T}) \geq -\frac{C_4}{T} + C_1q_T'+\frac{\gamma_0}{2}C_1q_T.
\end{equation}

\textbf{Case (b)}: $(t,x)\in (0,+\infty)\times \R$ such that $\xi_T\leq -M$.  

Proceeding similarly as above, we have $\rho=1$, $\rho^{\prime}=\rho^{\prime \prime}=0$, and
$$
\mathcal{L}_1 (\textbf{w}^{+}_{T})\geq -\frac{C_4}{T}+q_T^{\prime}-a_{1}p_{1}q_T\left(1-(w_{1,T}^++\varphi_1)\right)+k_{1}p_{2}(1-\varphi_2)q_T
-C_1q_Tk_{1}p_{2}\left(\varphi_1+q_T\right).
$$
It is clear that  $1-\varphi_2\geq0$ and, by the last line of \eqref{wave fix s +-}, there holds
$1-(w_{1,T}^++\varphi_1)<-1+\sigma_0$, and $0<\varphi_1+q_T\leq 1+\sigma_0$.
This together with  the  third inequality of  \eqref{gamma_0} implies that 
$$
\begin{aligned}
\mathcal{L}_1 (\textbf{w}^{+}_{T})&\geq -\frac{C_4}{T} +q_T^{\prime}+a_{1}p_{1}(1-\sigma_0)q_T-C_1(1+\sigma_0)q_Tk_{1}p_{2} \vspace{5pt}\\ 
& \geq -\frac{C_4}{T} +q_T^{\prime}+\gamma_0(1-\sigma_0)q_T-C_1(1+\sigma_0)q_Tk_{1}p_{2}.
\end{aligned}
$$
Furthermore, since $C_1=\gamma_0/(8\theta_+\gamma_+)$ and $\sigma_0\in (0,1/4)$, it follows that  
\begin{equation}\label{largeT-leq-M}
\mathcal{L}_1 (\textbf{w}^{+}_{T})\geq  -\frac{C_4}{T} +q_T^{\prime}+\frac{\gamma_0}{2}q_T.
\end{equation}

Now, we choose the function $q_T$ such that 
\begin{equation*}
-\frac{C_5}{T} +q_T^{\prime}+\frac{\gamma_0}{2}q_T=0\,\,\hbox{ for }\,\,t>0,\quad\hbox{and}\quad q_T(0)=\varepsilon,
\end{equation*}
where $C_5=\max\{C_4/C_1,C_4\}$, and choose 
\begin{equation}\label{choose-tT4}
	\tilde{T}_4=  \max\left\{  \frac{4C_5}{\varepsilon\gamma_0}  \,,T^* \right\}.
\end{equation}
It is straightforward to check that 
\bqq\label{function-q_T-largeT}
q_T(t)=\frac{2C_5}{T \gamma_0}+\left(\varepsilon-\frac{2C_5}{T \gamma_0}\right) \mathrm{e}^{-\gamma_0/2 t}\,\, \text{ for }\,\, t \geq 0,
\eqq
and it satisfies \eqref{defi q_T large T}. Furthermore, combining \eqref{largeT-geqM} and \eqref{largeT-leq-M}, one concludes that 
$\mathcal{L}_1 (\textbf{w}^{+}_{T})\geq 0$ for all $(t, x) \in (0, +\infty) \times \mathbb{R}$ such that $\xi_T \leq -M$  or $\xi_T \geq M$, provided that $\eta_T$ satisfies \eqref{defi q_T large T} too.

\textbf{Case (c)}: $(t,x)\in (0,+\infty)\times \R$ such that $-M<\xi_T<M$.

In this case, since $\kappa_T(\cdot)$ is required to be decreasing, it follows from \eqref{beta_2} that $\kappa'_T  \partial_{\xi}\varphi_1\geq -\beta_0 \kappa'_T$. Moreover, since $\rho(\cdot)$ satisfies \eqref{function-rho}, and since 
$0<q_T\leq \sigma_0/2$, $\sigma_0 \leq  \beta_0/(1+C_1)$ (see \eqref{delta_2}), there holds
$$q_T\rho'\left(1-C_1\right)\kappa_T^{\prime}\geq |q_T\rho'\left(1-C_1\right)|\kappa_T^{\prime}\geq \frac{\sigma_0}{2}(1+C_1)\kappa_T^{\prime}\geq \frac{\beta_0}{2}\kappa_T^{\prime},$$
 and 
$$-\left(c\rho'+d_{1,T}\rho''\right)(1-C_1)q_T\geq -\left(\max_{s\in\R}c(s)+\theta_+ \right)\left(1+C_1\right)q_T.$$ 
Using the fact that $w_{1,T}^++\varphi_1-1\geq -1$ and $\varphi_1\leq w_{1,T}^+\leq 2$, we compute that
$$\begin{aligned}
-a_{1}p_1\left(w_{1,T}^+(1-w_{1,T}^+)-\varphi_1(1-\varphi_1)\right)&
=-a_{1}p_{1}q_T\left((\rho+(1-\rho)C_1)(1-(w_{1,T}^++\varphi_1)\right)  \vspace{5pt}\\ 
&\geq-\theta_+\gamma_+\left(1+C_1\right) q_T,	
 \end{aligned}$$
and 
$$\begin{aligned}
 	k_{1}p_2\left(w_{1,T}^+(1-w_{2,T}^+)-\varphi_1(1-\varphi_2)\right)&\geq 
 	k_{1}p_{2}w_{1,T}^+\left( (1-w_{2,T}^+)-(1-\varphi_2)\right)  \vspace{5pt}\\ 
 	&\geq-2\theta_+\gamma_+ \left(1+C_1\right)q_T.
 \end{aligned}
$$
 
With the above estimates, we observe from \eqref{largeT-Lw} that 
$$\mathcal{L}_1 (\textbf{w}^{+}_{T})\geq -\frac{C_4}{T} -\frac{\beta_0}{2}\kappa_T' +\left(1+C_1\right)q_T^{\prime}-K_1q_T,$$
where 
$K_1=\left(\max_{s\in\R}c(s)+\theta_+ \right)\left(1+C_1\right)+3\theta_+\gamma_+\left(1+C_1\right)$.
Then, with $q_T(t)$ given by \eqref{function-q_T-largeT}, by choosing $\kappa(t)$ such that 
$$ -\frac{C_5(C_1+1)}{T} -\frac{\beta_0}{2}\kappa_T' +\left(1+C_1\right)q_T^{\prime}-K_1q_T=0 \,\,\hbox{ for }\,\,t>0,\quad\hbox{and}\quad \kappa_T(0)=0,$$
we obtain that $\mathcal{L}_1 (\textbf{w}^{+}_{T})\geq 0$ for $(t,x)\in (0,+\infty)\times \R$ such that $-M<\xi_T<M$ (here we use the fact $C_5(C_1+1)> C_4$). Direct computation gives that  
\bqq\label{function-kappa_T}
\kappa_T(t)=-K_2 \left(\frac{2C_5}{T \gamma_0}t+\frac{2}{\gamma_0} \left(\varepsilon-\frac{2C_5}{T \gamma_0}\right) \left(1-\mathrm{e}^{-\gamma_0/2 t} \right) \right) \quad \text {for } \,\, t \geq 0,\eqq
where 
$K_2=2\beta_0^{-1} \left(K_1+\gamma_0(C_1+1)/2 \right)$.
It is also easily checked that for any $T>\tilde{T}_4$,  $\kappa_T(t)$ satisfies \eqref{defi q_T large T}.  

Combining the above, one can conclude that for any $\varepsilon\in (0,\varepsilon_0)$ and any $T>\tilde{T}_4$, there holds $\mathcal{L}_1 (\textbf{w}^{+}_{T})\geq 0$ in $(0,+\infty)\times \R$. This completes the proof of Lemma \ref{large T supsub}.
\end{proof}

Similarly to Lemma \ref{large T supsub}, we have the following sub-solution for system \eqref{eq:change of variable} with large $T$.
 
 \begin{lem}\label{large T supsub2}
 Let $\varepsilon_0>0$ be provided by Lemma {\rm \ref{large T supsub}} and, for any $\varepsilon\in (0,\varepsilon_0)$, let $\tilde{T}_4>0$ be given by \eqref{choose-tT4}. Then for any $T>\tilde{T}_4$, the pair of functions 
 $\mathbf{w}_T^-(t,x)=(w_{1,T}^-,w_{2,T}^-)(t,x)$ defined by 
 \begin{equation*}
 	\left\{\begin{array}{l}
 		\displaystyle w_{1,T}^{-}(t, x)=\varphi_1\left(\tilde{\xi}_T; t/T\right)-q_T(t)\left(\rho\left(\tilde{\xi}_T\right)+\frac{\gamma_0}{8\theta_+\gamma_+}\left(1-\rho\left(\tilde{\xi}_T\right)\right)\right),  \vspace{5pt}\\  
 		\displaystyle w_{2,T}^{-}(t, x)=\varphi_2\left(\tilde{\xi}_T ; t/T\right)-q_T(t)\left(\frac{\gamma_0}{8\theta_+\gamma_+}\rho\left(\tilde{\xi}_T\right)+\left(1-\rho\left(\tilde{\xi}_T\right)\right)\right),
 	\end{array}\right.
 \end{equation*}
 is a sub-solution of \eqref{eq:change of variable} for $(t,x) \in (0,+\infty) \times \R$, where $\tilde{\xi}_T=x-X_T(t)-\kappa_T(t)$, and $q_T(\cdot)$, $\kappa_T(\cdot)$ are $C^1([0,+\infty))$ functions satisfying \eqref{defi q_T large T}, given by \eqref{function-q_T-largeT} and \eqref{function-kappa_T}, respectively.
 \end{lem} 
 
Building upon Lemmas \ref{large T supsub}-\ref{large T supsub2}, the proof of Theorem \ref{converge-largeT} follows a parallel structure to that of Theorem \ref{converge-samllT}. For clarity, we give its outline as follows. 
 
 \begin{proof}[Proof of Theorem  {\rm \ref{converge-largeT}}]
 For any $\varepsilon\in (0,\varepsilon_0)$ and any $T>\tilde{T}_4$, let $\mathbf{w}_T^+(t,x)=(w_{1,T}^+,w_{2,T}^+)(t,x)$ and $\mathbf{w}_T^-(t,x)=(w_{1,T}^-,w_{2,T}^-)(t,x)$ be, respectively, the super-solution and sub-solution of \eqref{eq:change of variable} in $(0,+\infty)\times\R$. The comparison principle yields that
 	$$
 	w_{i,T}^{-}(t, x) \leq w_{i,T}(t, x) \leq w_{i,T}^{+}(t, x) \,\, \text{ for all } \,\, (t,x)\in(0,+\infty)\times\R,\, i=1,2,
 	$$
 	where $\mathbf{w}_T(t,x):=(w_{1,T},w_{2,T})(t,x)$ is the solution of the Cauchy problem of \eqref{eq:change of variable} with initial data $(\varphi_1,\varphi_2)(x;0)$. By Proposition \ref{stability}, as $t\to +\infty$, this solution approaches the periodic traveling wave $(\varphi_{1,T},\varphi_{2,T})(x-c_Tt,t)$ in the following sense, for each $i=1,2$,
    $$
 	\sup _{x \in \mathbb{R}}\left|w_{i,T}(t, x)-\varphi_{i,T}\left(x-x_T^*-c_Tt,t\right)\right| \rightarrow 0 \,\, \text{ as }\,\, t \rightarrow+\infty
 	$$
 	for some constant $x_T^* \in \mathbb{R}$. 
 	Combining the above, one finds some large $t_0>0$ such that
 	$$
 	w_{i,T}^{-}(t,x)-\frac{\sigma_0}{2} \leq \varphi_{i,T}(x-x_T^*-c_T,t) \leq w_{i,T}^{+}(t,x)+\frac{\sigma_0}{2}\,\, \text { for all } \, t\geq t_0,\,x\in\R,
 	$$
 	where $\sigma_0 \in (0,1/4)$ is the constant provided by \eqref{sigma0-largeT}. 
 	By choosing $x=c_Tt+x_T^*$, $t=nT$, $i=1$, and using the normalization condition \eqref{normalize small T wave}, this in particular implies that 
 	 	$$
 	 w_{1,T}^{-}(nT,c_TnT+x_T^*)-\frac{\sigma_0}{2} \leq \frac{1}{2} \leq w_{1,T}^{+}(nT,c_TnT+x_T^*)+\frac{\sigma_0}{2}\,\, \text { for all large } \, n\in\N.
 	 $$
 	 Remember that $0\leq q_T(t) \leq \varepsilon_0 \leq \sigma_0/(2C_1+2)$ by \eqref{delta_2}. Then, by the definition of  $w_{1,T}^{+}$ and the periodicity of $\varphi_1$ with respect to its second variable, we obtain 
 	$$
 	\varphi_1\left(c_TnT+x_T^*-X_T(nT)-\kappa_T(nT); 0\right)-\sigma_0  \leq \frac{1}{2} \leq \varphi_1\left(c_TnT+x_T^*-X_T(nT)+\kappa_T(nT); 0\right)+\sigma_0. 
 	$$
 	Since $\sigma_0 \in (0,1/4)$ and since $\lim _{\xi \rightarrow+\infty} \varphi_1(\xi;0)=0$, $\lim _{\xi \rightarrow-\infty} \varphi_1(\xi; 0)=1$, there exists a constant $Z_0>0$ independent of $T$ such that
 	$$
 	\left|c_TnT+x_T^*-X_T(nT)\right| \leq\left|\kappa_T(nT)\right|+Z_0\,\, \text{ for all large }\, n\in\N.
 	$$
 	From the definition of $X_T(\cdot)$, one then infers that
 	\begin{equation*}
 	\left|  c_T-c^*\right| \leq \frac{\left|x_T^*\right|+Z_0}{nT}+\frac{\left|\kappa_T(nT)\right|}{nT} \,\, \text{ for all large }\, n\in\N.
 	\end{equation*}
 	Finally, passing to the limits as $n\to+\infty$, one easily sees from \eqref{function-kappa_T} that
 	$$
 	\left|  c_T-c^*\right| \leq \frac{2K_2C_5}{\gamma_0}\frac{1}{T} \,\,\hbox{ for all }\, T>\tilde{T}_4.
 	$$
 	Since $K_2$, $C_5$, $\gamma_0$ are all independent of $T$, the proof of Theorem \ref{converge-largeT} is thus complete. 
 \end{proof}

%%%%%%%%%%%%%%%%%%%%%%%%%%%%%%%%%%%%%%%%%%%%%%%%%%%%%%%%%%%%%%%%%%%%%%%%%%%%%%%%%%%%%%%%%%%%%%%%%%%%%%%%%%%%%%%%%%%%%%%
%%%%%%%%%%%%%%%%%%%%%%%%%%%%%%%%%%%%%%%%%%%%%%%%%%%%%%%%%%%%%%%%%%%%%%%%%%%%%%%%%%%%%%%%%%%%%%%%%%%%%%%%%%%%%%%%%%%%%%%%
\section{Sign of wave speeds: Proofs of Theorems \ref{Theorem  sign of speed} and \ref{theo-change-sign}}

\subsection{Proof of Theorem \ref{Theorem  sign of speed} }
Since $T>0$ is fixed in this subsection,  we simply the notations as follows: for each $i=1,2$, instead of $d_{i,T}$, $r_{i,T}$, $a_{i,T}$, $k_{i,T}$, $p_{i,T}$, we write $d_i$, $r_i$, $a_i$, $k_i$, $p_i$, respectively. Thanks to assumption \eqref{equiv-ad}, we further write
$$d(\cdot):=d_1(\cdot)=d_2(\cdot)\quad \hbox{and}\quad a(\cdot):=a_1(\cdot)=a_2(\cdot).$$
Moreover, let $(\varphi_{1}, \varphi_{2})(x-ct,t)$ denote a periodic traveling wave of \eqref{eq:change of variable} connecting $\textbf{0}$ and $\textbf{1}$. Clearly, the speed $c$ is identical to that of the periodic traveling wave 
$(\phi_{1}, \phi_{2})$ for the original system \eqref{eq:main}. Therefore, to prove Theorem \ref{Theorem sign of speed}, it suffices to determine the sign of $c$. We exploit a symmetry in the system by introducing the auxiliary profiles:
$$
\tilde{\varphi}_1(\xi, t)=1-\varphi_2(-\xi,t) \quad \text {and}\quad \tilde{\varphi}_2(\xi, t)=1-\varphi_1(-\xi,t)
\quad\hbox{for }\,\,t \in \mathbb{R},\,\xi \in \mathbb{R}.  $$
These functions will serve as comparison profiles, allowing us to relate the speed of the traveling wave $(\varphi_{1}, \varphi_{2})$  
to its reflected counterpart.
As a matter of fact, it is easily seen that $(\tilde{\varphi}_1,\tilde{\varphi}_2)(\xi,t)$ is $T$-periodic in $t\in\R$, decreasing in  $\xi\in\R$, that
\begin{equation*}%\label{tilde-varphi-connect}
	\left(\tilde{\varphi}_1,\tilde{\varphi}_2 \right)(-\infty, t)= \mathbf{1}\quad\hbox{and}\quad \left(\tilde{\varphi}_1,\tilde{\varphi}_2 \right)(+\infty, t)=\mathbf{0} \quad \hbox{uniformly in }\,\,t\in\R, 
\end{equation*}
and that $(\tilde{v}_1,\tilde{v}_2)(t,x):=(\tilde{\varphi}_1,\tilde{\varphi}_2)(x+ct,t)$ satisfies
\bqq\label{eq:main d1=d2}
\left\{\begin{array}{ll}
	\partial_t \tilde{v}_1=d(t)\partial_{x x} \tilde{v}_1+a(t)p_{2}(t)(1-\tilde{v}_{1})\tilde{v}_{1}-k_{2}(t){p}_{1}(t)\tilde{v}_{1}(1-\tilde{v}_{2}) & \hbox{in }\,\,\R^2, \vspace{5pt} \\
	\partial_t v_2=d(t)\partial_{x x} \tilde{v}_2-a(t) p_{1}(t)(1-\tilde{v}_{2})\tilde{v}_{2}+k_{1}(t)p_{2}(t)\tilde{v}_{1}(1-\tilde{v}_{2})  & \hbox{in }\,\,\R^2.
\end{array} \right.
\eqq
In other words, $(\tilde{\varphi}_1,\tilde{\varphi}_2)(x+ct,t)$ is a periodic traveling wave of system \eqref{eq:main d1=d2} connecting $\textbf{0}$ and $\textbf{1}$, and $-c$ is the associated wave speed.

\begin{proof}[Proof of Theorem {\rm \ref{Theorem  sign of speed} (i)-(ii)}]
We only show $c\geq 0$ under the assumption that  $k_{1}(\cdot) \leq k_{2}(\cdot)$, $r_{1}(\cdot) \geq r_{2}(\cdot)$ in $\R$, as we will sketch below that, the proof of $c\leq 0$ in the case where $k_{1}(\cdot) \geq k_{2}(\cdot)$, $r_{1}(\cdot) \leq r_{2}(\cdot)$ is analogous.  Assume by contradiction that $c<0$.  

Let $\mu>0$ be a positive constant satisfying \eqref{set-mu}, and for $i=1,2$, let $\psi_{i}^-(t)$ and $\psi_{i}^+(t)$ be, respectively, the unique solutions of \eqref{eigenfunction-0} and \eqref{eigenfunction-1} satisfying \eqref{normalize-supsub}. Moreover, let $\varepsilon_0>0$, and $K>0$ be the constants obtained in Lemma \ref{fix T sup sub}. 
Since both $(\varphi_{1}, \varphi_{2})$ and $(\tilde{\varphi}_1,\tilde{\varphi}_2)$ connect $\mathbf{0}$ and $\mathbf{1}$, and are decreasing in their first variables, one finds some  $\xi_0 \in \R$ such that
\begin{equation}\label{comp-initial}
	\left\{\begin{array}{l}
		\tilde{\varphi}_1(x, 0) \leq \varphi_1\left(x+\xi_0 , 0 \right)+\varepsilon_0\left(\rho(x+\xi_0 ) \psi_1^{+}(0)+(1-\rho(x+\xi_0 )) \psi_1^{-}(0)\right), \vspace{5pt} \\
		\tilde{\varphi}_2(x, 0) \leq \varphi_2\left(x+\xi_0, 0 \right)+\varepsilon_0\left(\rho(x+\xi_0 ) \psi_2^{+}(0)+(1-\rho(x+\xi_0 )) \psi_2^{-}(0)\right),
	\end{array}\right.
\end{equation}
for $x\in\R$, where $\rho(\cdot) \in C^2(\mathbb{R};[0,1])$ satisfies \eqref{function-rho}.  It follows directly from Lemma \ref{fix T sup sub} that the pair of functions $(v_1^+,v_2^+)(t,x)$ defined in \eqref{define-subsuper-solu}, with $\varepsilon=\varepsilon_0$, 
is a super-solution of \eqref{eq:change of variable} in $(t,x)\in (0,+\infty)\times\R$. Furthermore, by the assumptions $k_{1}(\cdot) \leq k_{2}(\cdot)$, $r_{1}(\cdot) \geq r_{2}(\cdot)$ and $a_1(\cdot)=a_2(\cdot)=a(\cdot)$, it is easily seen that 
$p_1(\cdot)\geq p_2(\cdot)$ (recall that $p_i(\cdot)$ is the unqiue positive solution to problem \eqref{T-periodic equilibria}), whence we have 
$$a(t)p_2(t)\leq a(t)p_1(t)\quad\hbox{and} \quad k_2(t)p_1(t)\geq  k_1(t)p_2(t).$$ 
This implies that  $(v_1^+,v_2^+)(t,x)$ is also a super-solution of system \eqref{eq:main d1=d2}. Remember that $(\tilde{\varphi}_1,\tilde{\varphi}_2)(x+ct,t)$ is an entire solution of \eqref{eq:main d1=d2}. Now, due to \eqref{comp-initial}, applying the comparison principle (Lemma \ref{comparison principle} (i)) to system \eqref{eq:main d1=d2}, we obtain that 
\begin{equation}\label{term comparison principle 1}
	(\tilde{\varphi}_1,\tilde{\varphi}_2)(x+ct,t)\leq  (v_1^+,v_2^+)(t,x)\, \text { in }\,  \R\,\hbox{ for all }\, t>0.
\end{equation}
In particular, taking $x=-cnT$ and $t=nT$ for $n\in\N$ in the above inequality gives that $\tilde{\varphi}_1(0, n T) \leq v^{+}_1(nT, -cn T)$.  Since $\tilde{\varphi}_1$ and $\varphi_1$ are $T$-periodic in their second variables, and $\psi_1^{\pm}$ are $T$-periodic, it follows that for each $n \in \mathbb{N}$,
\bqq\label{term comparison principle 2}
\tilde{\varphi}_1(0, 0) \leq \varphi_1\left( \xi_n,0\right)+\varepsilon_0 \mathrm{e}^{-\mu nT}\left(\rho(\xi_n) \psi_1^{+}(0)+(1-\rho(\xi_n)) \psi_1^{-}(0)\right),
\eqq
where $\xi_n=-2cnT+\xi_0-\varepsilon_0K(1-\mathrm{e}^{-\mu nT}) $. Notice that we have assumed $c<0$. It is clear that $\lim_{n\to+\infty}\xi_n=+\infty$. Then, passing to the limits as $n\to+\infty$ in \eqref{term comparison principle 2} yields that 
$ \tilde{\varphi}_1(0, 0)\leq  0$, which is impossible. 
Therefore, $c$ must be nonnegative.

Finally, in the case where $k_{1}(\cdot) \geq k_{2}(\cdot)$, $r_{1}(\cdot) \leq r_{2}(\cdot)$, it is easily checked that, under assumption \eqref{equiv-ad}, any sub-solution of \eqref{eq:change of variable} is automatically a sub-solution of \eqref{eq:main d1=d2}. Then with similar arguments as above (working this time with the sub-solution given in \eqref{define-subsuper-solu}), one can conclude that $c\leq 0$. This ends the proof of statements (i) and (ii).
\end{proof}

Next, we apply the strong comparison principle (Lemma \ref{comparison principle} (ii)) to show that the wave speed admits a strict sign.

\begin{proof}[Proof of Theorem {\rm \ref{Theorem  sign of speed} (iii)}]
We already know from statement (i) that if $k_1(\cdot) \leq k_2(\cdot)$ and $r_1(\cdot) \geq r_2(\cdot)$ in $\R$, then $c\geq 0$. It remains to show that, under the additional assumption that $k_1(\cdot) \not\equiv k_2(\cdot)$   or $r_1(\cdot)\not\equiv  r_2(\cdot)$, $c$ cannot be $0$.   

Assume by contradiction that $c=0$. Then, the same reasoning in showing \eqref{term comparison principle 1} yields that 
\begin{equation*}
	(\tilde{\varphi}_1,\tilde{\varphi}_2)(x,t)\leq  (v_1^+,v_2^+)(t,x) \, \text { in }\,  \R\,\hbox{ for all }\, t>0,
\end{equation*}
where $(v_1^+,v_2^+)(t,x)$ is the super-solution of \eqref{eq:change of variable} defined in \eqref{define-subsuper-solu} with $\varepsilon=\varepsilon_0$ and $\xi_0\in\R$ satisfying \eqref{comp-initial}. Taking $t=nT$ for $n\in\N$ and using the periodicity, we obtain 
	$$
\left\{\begin{array}{l}
	\tilde{\varphi}_1(x,0) \leq \varphi_1\left(\xi_n,0\right)+\varepsilon_0\mathrm{e}^{-\mu nT}\left(\rho(\xi_n)\psi_1^{+}(0)+(1-\rho(\xi_n)) \psi_1^{-}(0)\right), \vspace{5pt} \\
	\tilde{\varphi}_2(x,0) \leq \varphi_2\left(\xi_n,0\right)+\varepsilon_0\mathrm{e}^{-\mu nT}\left(\rho(\xi_n)\psi_2^{+}(0)+(1-\rho(\xi_n)) \psi_2^{-}(0)\right),
\end{array} \right.$$
for all $n\in\N$, $x\in\R$, where $\xi_n=x+\xi_0-\varepsilon_0K(1-\mathrm{e}^{-\mu nT})$.  Passing to the limits as $n\to+\infty$ gives
$$(\tilde{\varphi}_1, \tilde{\varphi}_2)(\cdot,0) \leq  (\varphi_1, \varphi_2)\left(\cdot+\xi_0-\varepsilon_0 K,0\right)\,\hbox{ in }\, \R.$$

Now, we set
$$
\xi_*:=\sup \left\{\xi \in \mathbb{R}:\,(\tilde{\varphi}_1, \tilde{\varphi}_2)(\cdot,0) \leq\left(\varphi_1, \varphi_2\right)(\cdot +\xi,0)\,\hbox{ in } \R\right\}.
$$
Since $(\varphi_1, \varphi_2)(\cdot,0)$ is decreasing in $\R$, it is clear that $\xi_*$ is a real number and $\xi_* \geq \xi_0-\varepsilon_0 K$. It is also easily seen that $(\tilde{\varphi}_1, \tilde{\varphi}_2)(\cdot,0) \leq (\varphi_1, \varphi_2)(\cdot +\xi_*,0)$ in $\R$. 
Furthermore, since either $k_1(\cdot) \not\equiv k_2(\cdot)$ or $r_1(\cdot)\not\equiv  r_2(\cdot)$, $(\varphi_1, \varphi_2)(x,t)$ cannot be a  solution of system \eqref{eq:main d1=d2}, and hence, there holds $\mathbf{0} \ll (\tilde{\varphi}_1, \tilde{\varphi}_2)(\cdot,0)<(\varphi_1, \varphi_2)(\cdot +\xi_*,0) \ll \mathbf{1}$ in $\R$. Remember that we assumed $k_1(\cdot) \leq k_2(\cdot)$ and $r_1(\cdot) \geq r_2(\cdot)$ in $\R$. Then,  $(\varphi_1, \varphi_2)$ is a strict super-solution of \eqref{eq:main d1=d2}, and using Lemma \ref{comparison principle} (ii), we obtain $(\tilde{\varphi}_1, \tilde{\varphi}_2)(\cdot,t)\ll (\varphi_1, \varphi_2)(\cdot+\xi_*,t)$ in $\R$ for each $t>0$. Furthermore, by the periodicity again, we have 
\bqq\label{strict compare}
	(\tilde{\varphi}_1, \tilde{\varphi}_2)(\cdot,0) \ll \left(\varphi_1, \varphi_2\right)(\cdot +\xi_*,0)\,\hbox{ in }\,\R.
\eqq

We are now ready to derive a contradiction. Denote $\nu_0=\min_{i=1,2} \{\psi_{i}^{+}(0),\psi_{i}^{-}(0)\}$, and take a large positive constant $C$ such that
	$$
	\delta_0:=\sup _{|x| \geq C-1}\left\{\left|\partial_x\varphi_1\left(x +\xi_*,0\right)\right|,\,\left|\partial_x\varphi_2\left(x +\xi_*,0\right)\right|\right\}\leq\frac{\nu_0}{2K}.
	$$
	Thanks to \eqref{strict compare} and the fact that  $(\varphi_1, \varphi_2)(\cdot,0)$ is decreasing in $\R$, one finds a small constant $0<\hat{\varepsilon} \leq \min \left\{1, \varepsilon_0 \nu_0 \delta_0^{-1}\right\}$ such that
	$$\left(\tilde{\varphi}_1, \tilde{\varphi}_2\right)(x,0) \leq\left(\varphi_1, \varphi_2\right)(x +\xi_*+\hat{\varepsilon},0)\, 
	\text { for  }\, x\in [-C, C].
	$$
	On the other hand, for $x \in(-\infty,-C) \cup(C,+\infty)$, by \eqref{strict compare}, the definition of $\delta_0$ and the monotonicity of $(\varphi_1, \varphi_2)(\cdot,0)$, one can check that $(\tilde{\varphi}_1, \tilde{\varphi}_2)(x,0) \leq (\varphi_1, \varphi_2)(x +\xi_*+\hat{\varepsilon},0)+\hat{\varepsilon} \delta_0$. 
	Recall that $\rho(\cdot) \in C^2(\mathbb{R};[0,1])$ satisfies \eqref{function-rho}. It then follows from the choice of $\nu_0$ that 
	\begin{equation}\label{hat-super-initial}
	\left\{\begin{array}{l}
		\tilde{\varphi}_1(x,0) \leq \varphi_1\left(x+\xi_*+\hat{\varepsilon},0\right)+\hat{\varepsilon} \delta_0 \nu_0^{-1}\left(\rho(x+\xi_*+\hat{\varepsilon})\psi_1^{+}(0)+(1-\rho(x+\xi_*+\hat{\varepsilon})) \psi_1^{-}(0)\right) \vspace{5pt} \\
		\tilde{\varphi}_2(x,0) \leq \varphi_2\left(x+\xi_*+\hat{\varepsilon},0\right)+\hat{\varepsilon} \delta_0 \nu_0^{-1}\left(\rho(x+\xi_*+\hat{\varepsilon})\psi_2^{+}(0)+(1-\rho(x+\xi_*+\hat{\varepsilon})) \psi_2^{-}(0)\right)
	\end{array} \right.
	\end{equation}
for all $x\in\R$. 
Let $(\hat{v}^{+}_1, \hat{v}^{+}_2): [0,+\infty)\times\R\to \R^2$ be the pair of functions defined as in \eqref{define-subsuper-solu} with $\varepsilon$ and $\xi_0$ replaced by $\hat{\varepsilon} \delta_0 \nu_0^{-1}$ and $\xi_*+\hat{\varepsilon}$, respectively. Since $0<\hat{\varepsilon} \delta_0 \nu_0^{-1} \leq \varepsilon_0$, it follows immediately from Lemma \ref{fix T sup sub} that $(\hat{v}^{+}_1, \hat{v}^{+}_2)$ is a super-solution of \eqref{eq:change of variable} for $(t,x)\in (0,+\infty)\times\R$, whence by the assumption  $k_1(\cdot) \leq k_2(\cdot)$ and $r_1(\cdot) \geq r_2(\cdot)$ in $\R$ again, it also serves as a super-solution of system \eqref{eq:main d1=d2}. Furthermore, since the inequalities in \eqref{hat-super-initial} are equivalent to the following:
$(\tilde{\varphi}_1, \tilde{\varphi}_2)(\cdot,0) \leq (\hat{v}_1^{+}, \hat{v}_2^{+})(0, \cdot)$ in $\R$, by using Lemma \ref{comparison principle} again, we obtain (remember that we assumed $c=0$)
$$\left(\tilde{\varphi}_1, \tilde{\varphi}_2\right)(\cdot,t) \leq \left(\hat{v}_1^{+}, \hat{v}_2^{+}\right)\left(t, \cdot \right) \, \text { in }\,  \R\,\hbox{ for all }\, t>0.$$
Choosing $t=nT$ for $n\in\N$ and using the periodicity gives that, for each $n\in\N$ and $i=1,2$, 
$$
\tilde{\varphi}_i(x,0) \leq  \varphi_i\left(\hat{\xi}_n,0\right)+ \hat{\varepsilon} \delta_0 \nu_0^{-1}\me^{-\mu nT}\left(\rho\left(\hat{\xi}_n\right)\psi_i^{+}(0)+\left(1-\rho\left(\hat{\xi}_n\right)\right) \psi_i^{-}(0)\right) \, \hbox{ for all }\,x\in \R,
$$
where  $\hat{\xi}_n=x+\xi_*+\hat{\varepsilon}-\hat{\varepsilon} \delta_0 \nu_0^{-1} K(1-\mathrm{e}^{-\mu nT})$. Since $0<\delta_0\leq\nu_0/(2K)$, sending to the limits as $n\to+\infty$ yields that
$$
\left(\tilde{\varphi}_1, \tilde{\varphi}_2\right)(\cdot,0) \leq(\varphi_1, \varphi_2)\left(\cdot+\xi_*+\hat{\varepsilon}-\hat{\varepsilon} \delta_0 v_0^{-1} K,0\right) \leq \left(\varphi_1, \varphi_2\right)\left(\cdot+\xi_*+ \frac{\hat{\varepsilon}}{2},0\right) \,\hbox{ in }\,\R,
$$
which is a contradiction with the definition of $\xi_*$. Therefore, the wave speed $c$ cannot be $0$, and consequently, we have $c>0$.

By a similar comparison argument, one sees that 
if $k_1(\cdot) \geq k_2(\cdot)$ and $r_1(\cdot) \leq r_2(\cdot)$ in $\R$, and if  $k_1(\cdot) \not\equiv k_2(\cdot)$   or $r_1(\cdot)\not\equiv  r_2(\cdot)$,  then there must hold $c<0$. This completes the proof of statement (iii) of Theorem \ref{Theorem  sign of speed}.   
\end{proof}

%%%%%%%%%%%%%%%%%%%%%%%%%%%%%%%%%%%%%%%%%%%%%%%%%%%%%%%%%%%%%%%%%%%%%%%%%%%%%%%%%%%%%%%%%%%%%%%%%%%%%%%%%%%%%%%%%%%%%%%%%%%%%%

\subsection{Proof of Theorem \ref{theo-change-sign}}

Theorem \ref{theo-change-sign} is proved by making use of 
the formulas for limiting speeds obtained in Theorems \ref{Theo small T} and \ref{Theo large T}, the explicit traveling wave solution for a special homogeneous system presented in \cite{M3}, and an approximation argument.

\begin{proof}
Let $(\tilde{r}_n(\cdot))_{n\in\N}\subset C^2(\R)$ be a sequence of positive, $1$-periodic functions such that
$3<\tilde{r}_n(\cdot)<13$ for each $n\in\N$, and that as $n\to+\infty$, $\tilde{r}_n(t)$ converges pointwise to a step function $\tilde{r}(t)$ given by  
\begin{equation*}
\tilde{r}(t)=\left\{\begin{array}{ll}
\displaystyle\frac{7}{2} &  \displaystyle \text {if }\, t \in \left(k, k+\frac{2}{3}\right], \vspace{6pt} \\
\displaystyle 12 & \displaystyle \text {if }\,  t \in \left(k+\frac{2}{3}, k+1\right], 
\end{array} \right. \text { for } \,  k \in \mathbb{Z}.
\end{equation*} 
Moreover, let $(\tilde{k}_n(\cdot))_{n\in\N}\subset C^2(\R)$ be a sequence defined by
$$\tilde{k}_n(\cdot)=2+\frac{4}{3} \tilde{r}_n(\cdot)\,\,\hbox{ for each }\, n\in\N. $$ 
For each $T>0$ and $n\in\N$, consider the following periodic competition system 
\begin{equation}\label{eq:exp}
	\left\{\begin{array}{ll}
		\displaystyle	\partial_t u_1= \partial_{xx}u_1+u_1\left(1- u_1 -\frac{1}{3}u_2\right), & \hbox{for }\,\, (t,x) \in \R^2,  \vspace{5pt}\\
		\displaystyle	\partial_t u_2= \partial_{xx}u_{2}+u_2\left(\tilde{r}_n\left(\frac{t}{T}\right)- u_2 -\tilde{k}_{n}\left(\frac{t}{T}\right) u_1\right), &\hbox{for }\,\, (t,x)\in \R^2.
	\end{array}\right.
\end{equation}  
Note that system \eqref{eq:exp} is a special case of  \eqref{eq:main} with coefficients chosen as
\begin{equation*}
	d_i(\cdot)\equiv a_i(\cdot)\equiv r_1(\cdot)\equiv 1  \,\hbox{ for }\, i=1,2,  \quad k_1(\cdot)\equiv\frac{1}{3}, 
\quad r_2(\cdot)\equiv \tilde{r}_n(\cdot)\quad \hbox{and} \quad k_2(\cdot) \equiv  \tilde{k}_n(\cdot).
\end{equation*}
We will show that, for all sufficiently large $n$,  \eqref{eq:exp} satisfies the conclusion of Theorem \ref{theo-change-sign}. 
To emphasize the dependence of $r_2$ and $k_2$ on the parameter $n$, we denote them by $r_2(\cdot;n)$ and $k_2(\cdot;n)$, respectively.

Let us first verify that for all large $n\in\N$, system \eqref{eq:exp} satisfies assumptions (A1) and (A2). To see (A1), observe that 
$\bar{d}_1=\bar{d}_2=\bar{a}_1=\bar{a}_2=\bar{r}_1=1$, $\bar{k}_1=1/3$, 
$$\bar{r}_2(n)=\int_{0}^1 r_2(t;n)dt=\int_{0}^1 \tilde{r}_n(t)dt,\quad \bar{k}_2(n)=\int_{0}^1 k_2(t;n)dt=\int_{0}^1 \tilde{k}_n(t)dt=2+\frac{4}{3}\bar{r}_2(n),  $$
and 
$$\bar{p}_1=\frac{\bar{r}_1}{\bar{a}_1}=1,\quad  \bar{p}_2(n)=\frac{\bar{r}_2(n)}{\bar{a}_2}=\bar{r}_2(n). $$
By the Lebesgue dominated convergence theorem, $\bar{r}_2(n)\to 19/3$ as $n\to+\infty$. This implies that, as $n\to+\infty$, 
\begin{equation*}
\left\{\begin{array}{l}
\displaystyle \bar{r}_1-\bar{k}_1\bar{p}_2(n)=1-\frac{1}{3}\bar{r}_2(n)\rightarrow -\frac{10}{9},		\vspace{5pt}\\
\displaystyle \bar{r}_2(n)-\bar{k}_2(n)\bar{p}_1=-2-\frac{1}{3}\bar{r}_2(n)\to -\frac{37}{9}, \vspace{5pt}\\
\displaystyle \bar{a}_1\bar{a}_2-\bar{k}_1\bar{k}_2(n)=\frac{1}{3}-\frac{4}{9}\bar{r}_2(n)\rightarrow -\frac{67}{27}.
\end{array}\right.
\end{equation*}  
Therefore, for all large $n\in\N$, system \eqref{eq:exp} satisfies (A1). 
Regarding (A2), since $p_1\equiv 1$, $p_2(s;n)=r_2(s;n)/a_2=r_2(s;n)$, and since $r_2(s;n)>3$ for all $s\in\R$ and $n\in\N$, it is easily seen that 
\begin{equation*}
	\left\{\begin{array}{l}
		\displaystyle r_1(s)-k_1(s)p_2(s;n)=1-\frac{1}{3}r_2(s;n)<0,		\vspace{5pt}\\
		\displaystyle r_2(s;n)-k_2(s;n)p_1(s)=-2-\frac{1}{3}r_2(s;n)<0, \vspace{5pt}\\
		\displaystyle a_1(s)a_2(s)-k_1(s)k_2(s;n)=\frac{1}{3}-\frac{4}{9}r_2(s;n)<0,
	\end{array}\right.
\end{equation*} 
for all $s\in\R$ and $n\in\N$. Consequently, for each $n\in\N$, system \eqref{eq:exp} satisfies (A2). 

It then follows from Theorems \ref{Theo small T}-\ref{Theo large T} that, there exist $N\in\N$, $0<T_*<T^*<+\infty$ such that for each $n\geq N$ and each $T\in (0,T_*)\cup (T^*,+\infty)$, system \eqref{eq:exp} admits a periodic traveling wave connecting two semi-trivial states with a unique wave speed, denoted by $c_{T,n}$. Furthermore, 
\begin{equation}\label{limit-speedsn}
\lim_{T\to 0^+} c_{T,n} =c_{0,n} \quad\hbox{and}\quad  \lim_{T\to +\infty} c_{T,n} = \int_0^1 c(s;n)ds,  
\end{equation}
where $c_{0,n}$  is the unique wave speed for the  homogenized system: 
\begin{equation*}
	\left\{\begin{array}{ll}
		\displaystyle	\partial_t u_1= \partial_{xx}u_1+u_1\left(1- u_1 -\frac{1}{3}u_2\right), & \hbox{for }\,\, (t,x) \in \R^2,  \vspace{5pt}\\
		\displaystyle	\partial_t u_2= \partial_{xx}u_{2}+u_2\left(\bar{r}_2(n)- u_2 -\bar{k}_2(n)u_1\right) , &\hbox{for }\,\, (t,x)\in \R^2,
	\end{array}\right.
\end{equation*}  
and for each $s\in\R$, $c(s;n)$ is the unique wave speed for the system with coefficients frozen at $s$: 
\begin{equation*}
	\left\{\begin{array}{ll}
		\displaystyle	\partial_t u_1= \partial_{xx}u_1+u_1\left(1- u_1 -\frac{1}{3}u_2\right), & \hbox{for }\,\, (t,x) \in \R^2,  \vspace{5pt}\\
		\displaystyle	\partial_t u_2= \partial_{xx}u_{2}+u_2\left(r_2(s;n)- u_2 -k_2(s;n)u_1\right), &\hbox{for }\,\, (t,x)\in \R^2.
	\end{array}\right.
\end{equation*}  
Notice that the above two homogeneous system are exactly the cases where explicit wave speeds and profiles can be obtained (see \cite[Section 3.1.1]{M3}). Specifically, for each $n\geq N$ and $s\in\R$, there holds
$$c_{0,n}= \left(-2+\frac{1}{3}\bar{r}_n\right) \left(\frac{2}{3}\bar{r}_n\right)^{-1/2}\quad\hbox{and}\quad c(s;n)=\left(-2+\frac{1}{3}r(s;n)\right)\left(\frac{2}{3}r(s;n)\right)^{-1/2}.
$$
Passing to the limits as $n\to+\infty$, we obtain 
$$\lim_{n\to+\infty} c_{0,n}=\frac{\sqrt{38}}{114}>0 \quad\hbox{and}\quad \lim_{n\to+\infty}\int_0^1 c(s;n)ds =\frac{\sqrt{2}}{6}-\frac{5\sqrt{21}}{63}<0.$$
This immediately implies that  $c_{0,n}>0$ and $\int_0^1 c(s;n)ds<0$ for all sufficiently large $n\in \N$. Hence, by \eqref{limit-speedsn}, for such $n$, system \eqref{eq:exp} is the desired example. The proof of Theorem \ref{theo-change-sign} is thus complete. 
\end{proof}

%%%%%%%%%%%%%%%%%%%%%%%%%%%%%%%%%%%%%%%%%%%%%%%%%%%%%%%%%%%%%%%%%%%%%%%%%%%%%%%%%%%%%%%%%%%%%%%%%%%%%%%%%%%%%%%%%%%%%%%%%%%
%%%%%%%%%%%%%%%%%%%%%%%%%%%%%%%%%%%%%%%%%%%%%%%%%%%%%%%%%%%%%%%%%%%%%%%%%%%%%%%%%%%%%%%%%%%%
\section{Appendix}
The appendix is devoted to the proof of Proposition \ref{fix s continue} under assumption (A2). Since the proof is parallel to that in the scalar case (see \cite{dhl}), we only give its outline and provide the details when considerable modifications are needed.  Recall that for each $s\in\R$, $(\boldsymbol \varphi(\xi;s),c(s))$ is the unique solution of \eqref{eq: frozen-limit-change} satisfying the normalization \eqref{normalize fix s wave}. 
Let us first show the continuity of $c(s)$ and $\boldsymbol \varphi(\xi;s)$ with respect to $s$, as stated in the following lemma.

\begin{lem}\label{speed conti}
	For any sequence $\left(s_n\right)_{n \in \mathbb{N}} \subset \mathbb{R}$ such that $s_n \rightarrow s_0$ as $n \rightarrow+\infty$ for some $s_0 \in \mathbb{R}$, we have $c\left(s_n\right) \rightarrow c\left(s_0\right)$ as $n \rightarrow+\infty$, and $\boldsymbol \varphi(\cdot;s_n) \rightarrow \boldsymbol \varphi(\cdot;s)$ in $C_{loc}^2(\mathbb{R})$ as $n \rightarrow+\infty$.
\end{lem}

\begin{proof}
First, we  show that the function $s\mapsto c(s)$ is bounded. It is well known (see e.g., \cite{aw}) that for each $i=1,2$ and $s\in\R$, the Fisher-KPP equation
$\partial_t w= d_i(s) \partial_{xx} w + r_i(s)w (1-w)$ in $\mathbb{R}^2$
admits a traveling wave $w(t,x)=\phi(x-ct)$ with $\phi(-\infty)=1$ and $\phi(+\infty)=0$ if and only if $c\geq 2\sqrt{d_i(s)r_i(s)}$. Moreover, a simple comparison argument yields that 
	$$-2\sqrt{d_2(s)r_2(s)} \leq c(s) \leq 2\sqrt{d_1(s)r_1(s)} \, \, \hbox{ for each }\,s\in\R. $$
	Then, by using \eqref{bound-rak}, we have $-2\theta_-\leq c(s)\leq 2\theta_+$ for  $s\in\R$.  Therefore, $c(s)$ is bounded in $s\in\R$. 
	
	Next, we show that $c\left(s_n\right) \rightarrow c\left(s_0\right)$ as $n \rightarrow+\infty$. By the periodicity, we may assume without loss of generality, that $\left(s_n\right)_{n \in \mathbb{N}} \subset[0,1]$ and $s_0 \in[0,1]$. 
	Assume by contradiction that $c\left(s_n\right)$ does not converge to $c\left(s_0\right)$ as $n \rightarrow+\infty$. Then, one finds some $c_{\infty} \in\left[-2\theta_-, 2\theta_+\right]$ such that $c_{\infty} \neq c\left(s_0\right)$ and, up to extraction of some subsequence, $c\left(s_n\right) \rightarrow c_{\infty}$ as $n \rightarrow+\infty$. 
	Remember that for each $n\in\R$, $(\varphi_1,\varphi_2)(\xi;s_n)$ is an entire solution of \eqref{eq: frozen-limit-change} with $s=s_n$. By standard elliptic estimates, there exists a function $\boldsymbol\varphi_{\infty}:=(\varphi_{1,\infty},\varphi_{2,\infty}) \in C^2(\R;\R^2)$ such that, 
    up to extraction of a further subsequence, 
\begin{equation}\label{a-varphi-n}
 \boldsymbol \varphi(\cdot;s_n) \rightarrow \boldsymbol\varphi_{\infty}(\cdot)\, \hbox{ in } \, C_{l o c}^2(\mathbb{R})\,\hbox{ as }\, n \rightarrow+\infty.
\end{equation} 
Clearly,  $\boldsymbol\varphi_{\infty}(\xi)$ satisfies
\bqq\label{fix s_0}
\left\{\begin{array}{l}
d_1(s_0)  \varphi''_{1, \infty}+c_\infty  \varphi'_{1, \infty}+g_{1}\left(s_0, \boldsymbol\varphi_{\infty}\right)=0 \quad \hbox{in }\,\,\R,   \vspace{5pt}\\
d_2(s_0)  \varphi''_{2, \infty}+c_\infty \varphi'_{2, \infty}+g_{2}\left(s_0, \boldsymbol\varphi_{\infty}\right)=0 \quad \hbox{in }\,\,\R,  \vspace{5pt}\\
 \mathbf{0}\leq  (\varphi_{1,\infty},\varphi_{2,\infty})(\cdot)\leq \mathbf{1} \quad \hbox{in }\,\,\R.
\end{array}\right.\eqq
It is also easily seen that for each $i=1,2$, $\varphi_{i,\infty}(\cdot)$ is nonincreasing.
This implies that the limits $\boldsymbol\varphi_{\infty}( \pm \infty)$ exist, and they must be equilibria of system \eqref{eq:fix-s-Lv} at $s=s_0$, and hence, $\boldsymbol\varphi_{\infty}( \pm \infty) \in \{\textbf{0},\textbf{1},\textbf{e}_{0}, \textbf{e}^*_{s_0}\}$. Moreover, due to \eqref{normalize fix s wave}, there holds 
\begin{equation}\label{normal-limit}
\varphi_{1,\infty}(0)=\gamma_1.
\end{equation}
Then, there must hold  $\varphi_{1,\infty}(+\infty)=0$, and one of the following three possibilities happens:
\begin{itemize}
	\item [{\rm (a)}] $\boldsymbol\varphi_{\infty}(+\infty)=\mathbf{0}$ and $\boldsymbol\varphi_{\infty}(-\infty)=\mathbf{1}$;
	\item [{\rm (b)}] $\boldsymbol\varphi_{\infty}(+\infty)=\mathbf{0}$ and $\boldsymbol\varphi_{\infty}(-\infty)=\textbf{e}^*_{s_0}$; 
	\item [{\rm (c)}] $\boldsymbol\varphi_{\infty}(+\infty)=\mathbf{e}_0$, $\boldsymbol\varphi_{\infty}(-\infty)=\mathbf{1}$, and $\varphi_{2,\infty} \equiv 1$.
\end{itemize}
We will derive a contradiction in each of these three cases.
	
If Case (a) happens, then it is easily seen from \eqref{fix s_0} that $\boldsymbol\varphi_{\infty}(x-c_{\infty}t)$ is a traveling wave of 
\bqq\label{origial-fix-s}
\left\{\begin{array}{l}
\partial_t v_1=d_{1}(s_0) \partial_{xx}v_1+g_{1}\left(s_0,v_1,v_2\right)=0 \quad \hbox{in }\,\,\R^2,   \vspace{5pt}\\
\partial_t v_2=d_{2}(s_0) \partial_{xx}v_{2}+g_{2}\left(s_0, v_1,v_2\right)=0 \quad \hbox{in }\,\,\R^2,  
\end{array}\right.\eqq
connecting $\mathbf{0}$ and $\mathbf{1}$. Yet, we see from \eqref{eq: frozen-limit-change} that $\boldsymbol\varphi(x-c(s_0)t;s_0)$ is such a traveling wave. By uniqueness, we have 
 $c_\infty= c(s_0)$, which contradicts the assumption that $c_\infty\neq  c(s_0)$. 

Suppose that Case (b) occurs. Then, $\boldsymbol\varphi_{\infty}(x-c_{\infty}t)$ is a traveling wave of \eqref{origial-fix-s} connecting $\mathbf{0}$  and $\textbf{e}^*_{s_0}$. Notice that system \eqref{origial-fix-s} restricted to the region  
$$
\Omega_-=\left\{\left(v_1, v_2\right) \in C\left(\mathbb{R}; \mathbb{R}^2\right) \mid \mathbf{0} \leq \left(v_1, v_2\right)\leq \textbf{e}^*_{s_0} \right\}
$$
admits a monostable structure, that is, $\mathbf{0}$ is asymptotically stable, $\textbf{e}^*_{s_0}$ is asymptotically unstable, and apart from $\mathbf{0}$ and $\textbf{e}^*_{s_0}$, \eqref{origial-fix-s} has no other constant steady state in  $\Omega_-$. Then, it is known (see e.g., \cite{lz07}) that there exists a maximal speed $c_- <0$ such that \eqref{origial-fix-s} has a traveling wave 
connecting $\mathbf{0}$  and $\textbf{e}^*_{s_0}$ with speed $c$ 
 if and only if $c \leq c_-$. This immediately implies that $c_{\infty} \leq c_-<0$. 
 
 To find a contradiction, we need to take the convergence of $\boldsymbol{\varphi}(\xi;s_n)$ along another normalization. For each $n \in \mathbb{N}$, since $\varphi_{1}(\xi;s_n) \rightarrow 0$ as $\xi \rightarrow +\infty$ and $\varphi_{1}(\xi;s_n) \rightarrow 1$ as $\xi \rightarrow-\infty$, and since $\varphi_{1}(\xi;s_n)$ is decreasing in $\xi \in \mathbb{R}$, there is, by continuity, a unique $r_n \in \mathbb{R}$ such that
 \bqq\label{norm r_n}
  \varphi_{1}\left(r_n;s_n\right)=\gamma_2,
 \eqq
 where $\gamma_2$ is a positive constant satisfying $\max_{s\in\R}\{ v_{1,s}^*,v_{2,s}^* \}<\gamma_2<1$. Clearly,  $\boldsymbol{\varphi}(\xi+r_n;s_n)$ remains an entire solution of \eqref{eq: frozen-limit-change} with $s=s_n$. By \eqref{normalize fix s wave} and the monotonicity of $\varphi_{1}(\xi;s_n)$, we have $r_n<0$.
 We further claim that $r_n \rightarrow-\infty$ as $n \rightarrow +\infty$. Otherwise, 
 as $n \rightarrow +\infty$, possibly up to a subsequence, there would hold  $r_n \rightarrow r_{\infty} \in (-\infty,0]$, whence by \eqref{a-varphi-n} and \eqref{norm r_n}, $\varphi_{1,\infty}\left(r_{\infty}\right)=\gamma_2>\max_{s\in\R}\{ v_{1,s}^*,v_{2,s}^*\} $. This would be  impossible, due to $\varphi_{1,\infty}(\xi) \in (0,v_{1,s_0}^*)$ for all $\xi\in\R$ in Case (b). Therefore, $r_n \rightarrow-\infty$ as $n \rightarrow +\infty$.
 
 Now, the same reasoning in showing \eqref{a-varphi-n} gives the existence of another entire solution of \eqref{fix s_0}, denoted by $\boldsymbol{\psi}=\left(\psi_1, \psi_2\right) \in C^2\left(\mathbb{R}; \mathbb{R}^2\right)$, such that
  up to extraction of some subsequence, $\boldsymbol{\varphi}\left(\cdot+r_n;s_n\right)$ converges to $\boldsymbol{\psi}(\cdot)$ in $C^2_{loc}(\R)$ as $n \rightarrow +\infty$. It is clear that  $\boldsymbol{\psi}(\xi)$  is nonincreasing in $\xi$, and by \eqref{norm r_n}, we have $\psi_1(0)=\gamma_2$. Since $r_n \rightarrow-\infty$ as $n \rightarrow +\infty$, and since $\varphi_{2}(\xi;s_n)$ is decreasing in $\xi \in \mathbb{R}$ and $\varphi_{2,\infty}(-\infty)=v_{2,s_0}^*$, it follows that for any $\xi \in \mathbb{R}$,
 $$
 \psi_2(\xi)=\lim _{n \rightarrow \infty} \varphi_{2}\left(r_n+\xi;s_n\right)\geq v_{2,s_0}^*.
 $$
Furthermore, since $\boldsymbol{\psi}(\pm \infty)$ are equilibria of system \eqref{eq:fix-s-Lv} with $s=s_0$, there must hold $\boldsymbol{\psi}(-\infty)=\textbf{1}$ and $\boldsymbol{\psi}(+\infty)=\textbf{e}^*_{s_0}$. In other words, $\boldsymbol{\psi}(x-c_{\infty}t)$ is a traveling wave of \eqref{origial-fix-s} connecting $\textbf{e}^*_{s_0}$ and $\textbf{1}$. Since system \eqref{origial-fix-s} restricted to the region  
$$
\Omega_+=\left\{\left(v_1, v_2\right) \in C\left(\mathbb{R}; \mathbb{R}^2\right) \mid  \textbf{e}^*_{s_0}  \leq (v_1,v_2)\leq  \mathbf{1} \right\}
$$
admits a monostable structure with $\textbf{e}^*_{s_0}$ being asymptotically unstable and  $\mathbf{1}$ being asymptotically stable, there exists a minimal speed $c_+>0$ such that \eqref{origial-fix-s} has a traveling wave with speed $c$ if and only if $c\geq c_+$. This implies that $c_{\infty} \geq c_+>0$, which is a contradiction with $c_{\infty} \leq c_-<0$. Therefore, Case (b) is ruled out.

It remains to exclude Case (c). In this case, it is easily seen that $\varphi_{1,\infty}(x-c_{\infty}t)$ is a traveling wave, connecting $0$ and $1$, of the Fisher-KPP equation $\partial_t w= d_1(s_0) \partial_{xx} w + r_1(s_0)w (1-w)$ in $\R^2$.  
Then, we have $c_\infty\geq2\sqrt{d_1(s_0)r_1(s_0)}>0$.
In order to find a contradiction, for each $n\in\N$, we choose $\tilde{r}_n\in\R$ such that 
 $$\varphi_{2}\left(\tilde{r}_n;s_n\right)=\gamma_1.$$ 
Proceeding similarly as in Case (b), one can conclude that $\tilde{r}_n\to +\infty$ as $n\to+\infty$, and that there exists $\boldsymbol{\tilde{\psi}}=(\tilde{\psi}_1, \tilde{\psi}_2) \in C^2(\mathbb{R};\mathbb{R}^2)$ which
is an entire solution of \eqref{fix s_0} such that, up to extraction of some subsequence, $\boldsymbol{\varphi}(\cdot+\tilde{r}_n;s_n)$ converges to $\boldsymbol{\tilde{\psi}}(\cdot)$ in $C^2_{loc}(\R)$ as $n \to \infty$.
One then infers that $\tilde{\psi}_1 \equiv 0$, $\tilde{\psi}_2(+\infty)=0$ and $\tilde{\psi}_2(-\infty)=1$. Namely, $\tilde{\psi}_2(x-c_{\infty}t)$ is a traveling wave, connecting $0$ and $1$, of the equation
$\partial_t w= d_2(s_0) \partial_{xx} w - r_2(s_0)w (1-w)$ in $\R^2$. Consequently, there holds $c_\infty\leq -2\sqrt{d_2(s_0)r_2(s_0)}<0$, which contradicts $c_\infty\geq 2\sqrt{d_1(s_0)r_1(s_0)}>0$. Hence,  Case (c) is ruled out, too.
Combining the above, we obtain that $c\left(s_n\right) \rightarrow c\left(s_0\right)$ as $n \rightarrow+\infty$.

Finally, since $\boldsymbol\varphi_{\infty}(x-c(s_0)t)$ is the unique traveling wave of \eqref{origial-fix-s} connecting $\mathbf{0}$ and $\mathbf{1}$ and satisfying \eqref{normal-limit}, it follows that $\boldsymbol\varphi_{\infty}(\cdot)\equiv \boldsymbol\varphi(\cdot;s_0)$. Furthermore, from the above proof, one sees that the whole sequence $(\boldsymbol \varphi(\cdot;s_n))_{n \in \mathbb{N}}$ converges to $\boldsymbol \varphi(\cdot;s_0)$ in $C_{loc}^2(\mathbb{R})$ as $n \rightarrow+\infty$. This completes the proof of Lemma \ref{speed conti}.
\end{proof}

Statements (i) and (ii) of Proposition \ref{fix s continue}  can be easily deduced from Lemma \ref{speed conti}, by invoking the monotonicity of $\boldsymbol \varphi(\xi;s)$ with respect to $\xi\in\R$ and its periodicity in $s\in\R$. 
To prove statement (iii), we need the following property that $\boldsymbol \varphi(\xi;s)$ converges to its limiting states as $\xi\to\pm \infty$ exponentially fast uniformly in $s\in\R$.

\begin{lem}\label{uniform exp s}
There exist some positive constants $\mu_1$, $\mu_2$, $M$, $C_1$ and $C_2$ $($all are independent of s$)$ such that 
\bqq\label{uniform exp cover}
	\left\{\begin{array}{ll}
	0<\varphi_i(\xi; s) \leq C_1 \mathrm{e}^{-\mu_1 \xi} & \text {for all }\,\, \xi \geq M,\, s \in \mathbb{R},\,i=1,2, \vspace{5pt}\\ 
	0<1-\varphi_i(\xi;s) \leq C_2 \mathrm{e}^{\mu_2 \xi} & \text {for all }\,\, \xi \leq-M,\, s \in \mathbb{R},\,i=1,2.
\end{array}\right.\eqq
\end{lem}

\begin{proof}
We only give the proof of the first inequality in \eqref{uniform exp cover}, since the analysis of the second one is analogous. 
Let $\varepsilon_0 >0$ be such that
\bqq\label{varepsilon_0 appendix}
\varepsilon_0=\min\left\lbrace \frac{1}{2},\, \frac{\gamma_0}{2\theta_+\gamma_+} \right\rbrace,
\eqq
where $\gamma_0$, $\theta_+$ and $\gamma_+$ are the positive constants provided by \eqref{gamma_0}, \eqref{bound-rak} and \eqref{bound-pT}, respectively. 
By Proposition \ref{fix s continue} (i), there exists a constant $M>0$ (independent of $s$) such that 
\begin{equation}\label{choose-M1-app}
0<\varphi_i(\xi;s)\leq \varepsilon_0\,\,\hbox{ for all } \xi \geq M,\, s \in \mathbb{R},\,i=1,2.
\end{equation}
	
We first prove that there exist $B_1>0$ and $\nu_1>0$ such that 
\begin{equation}\label{exp-varphi-1}
0< \varphi_1(\xi;s) \leq B_1 \mathrm{e}^{-\nu_1\xi} \,\,\hbox{ for all }\, \xi \geq M,\,s\in\R. 
\end{equation}	
By the first component of \eqref{eq: frozen-limit-change} and \eqref{choose-M1-app}, for  $s\in \mathbb{R}$ and  $\xi \geq M$, we have
$$ \begin{aligned}
0&\leq d_1(s) \partial_{\xi \xi} \varphi_{1}(\xi; s)+c(s) \partial_{\xi} \varphi_{1}(\xi; s)+r_{1}(s) \varphi_{1}(\xi; s)-k_{1}(s)p_{2}(s)(1-\varepsilon_0)\varphi_{1}(\xi; s) \vspace{5pt}\\
		&\leq d_1(s) \partial_{\xi \xi} \varphi_1(\xi; s)+c(s) \partial_{\xi} \varphi_1(\xi; s)-\frac{\gamma_0}{2} \varphi_1(\xi; s),
\end{aligned} $$
where the last inequality follows from \eqref{gamma_0} and \eqref{varepsilon_0 appendix}.	
We set
$$P_1(y;s)=d_1(s)y^2-c(s)y-\frac{\gamma_0}{2}\,\,\hbox{ for }\,\,y\in\R,\,s\in\R.$$
Since the functions $s \mapsto d_1(s)$ and $s \mapsto c(s)$ are continuous, positive and 1-periodic, we can choose a constant $\nu_1$ such that
	$$
	0<\nu_1 \leq \min _{s \in \mathbb{R}} \frac{c(s)+\sqrt{c^2(s)+2 \gamma_0 d_1(s)}}{2d_1(s)}.$$
Clearly, $P_1(\nu_1;s)\leq 0$ for all $s\in\R$. 
Define $\bar{v}_1(\xi)=B_1 \mathrm{e}^{-\nu_1 \xi}$ for $\xi\geq M$, where $B_1=\varepsilon_0 \mathrm{e}^{\nu_1 M}$. For each $s \in \mathbb{R}$, it is easily checked that $\varphi_1(M;s)\leq \varepsilon_0=\bar{v}_1(M)$, and that 
 $d_1(s) \bar{v}_1''+c(s) \bar{v}_1'-(\gamma_0\bar{v}_1)/2  \leq 0$ for $\xi \geq M$. Then, by the elliptic weak maximum principle, we obtain that $\varphi_1(\xi;s)\leq\bar{v}_1(\xi)$ for all $\xi \geq M$ and $s\in\R$. This completes the proof of \eqref{exp-varphi-1}. 

Next, we prove that there exist $B_2>0$ and $\nu_2>0$ such that 
\begin{equation}\label{exp-varphi-2}
	0< \varphi_2(\xi;s) \leq B_2 \mathrm{e}^{-\nu_2\xi} \,\,\hbox{ for all }\, \xi \geq M,\,s\in\R. 
\end{equation}	
Similarly as above, by the second component of \eqref{eq: frozen-limit-change} and \eqref{choose-M1-app}, \eqref{gamma_0}, \eqref{varepsilon_0 appendix}, for $s\in \mathbb{R}$ and  $\xi \geq M$, we have
	$$ \begin{aligned}
		0&\leq d_2(s) \partial_{\xi \xi} \varphi_{2}(\xi; s)+c(s) \partial_{\xi} \varphi_{2}(\xi; s)-r_2(s)(1-\varepsilon_0)\varphi_{2}(\xi; s)+k_{2}(s)p_{1}(s)\varphi_{1}(\xi; s) \vspace{5pt}\\
		&\leq d_2(s) \partial_{\xi \xi} \varphi_{2}(\xi; s)+c(s) \partial_{\xi} \varphi_{2}(\xi; s)-\frac{\gamma_0}{2}\varphi_{2}(\xi; s)+k_{2}(s)p_{1}(s)\varphi_{1}(\xi; s).
	\end{aligned}
	$$ 
It further follows from \eqref{exp-varphi-1} that 
$$0\leq d_2(s) \partial_{\xi \xi} \varphi_{2}(\xi; s)+c(s) \partial_{\xi} \varphi_{2}(\xi; s)-\frac{\gamma_0}{2}\varphi_{2}(\xi; s)+k_{2}(s)p_1(s) B_1 \mathrm{e}^{-\nu_1\xi}$$ 
for all $s\in \mathbb{R}$ and  $\xi \geq M$. Now, we set
$$P_2(y;s)=d_2(s)y^2-c(s)y-\frac{\gamma_0}{2}\,\,\hbox{ for }\,\,y\in\R,\,s\in\R,$$
and choose a positive constant $\nu_2$ such that
	$$0<\nu_2< \min\left\{\nu_1,\,\min _{s \in \mathbb{R}} \frac{c(s)+\sqrt{c^2(s)+2 \gamma_0 d_2(s)}}{2d_2(s)}\right\}.$$
It is clear that $P_2(\nu_2;s)<0$ for all $s\in\R$, and it is periodic, continuous in $s$. Define $\bar{v}_2(\xi)=$ $B_2 \mathrm{e}^{-\nu_2 \xi}$ for $\xi\geq M$, where $$B_2=\max\left\{ \varepsilon_0 \mathrm{e}^{\nu_2 M},\, \max_{s\in\R}  \frac{ B_1\theta_+\gamma_+}{|P_2(\nu_2;s)|}\right\}. $$ 
For each $s\in\R$, it is then straightforward to check that $\varphi_2(M;s)\leq\bar{v}_2(M)$, and that
$$d_2(s) \bar{v}_2^{\prime \prime}+c(s) \bar{v}_2^{\prime}-\frac{\gamma_0}{2} \bar{v}_2 +k_{2}(s) p_1(s)B_1 \mathrm{e}^{-\nu_1\xi}=P_2(\nu_2;s)B_2\mathrm{e}^{-\nu_2 \xi}+k_{2}(s)p_1(s) B_1 \mathrm{e}^{-\nu_1\xi} \leq 0$$
for all $\xi \geq M$. Furthermore, by the elliptic weak maximum principle again, we obtain that $\varphi_2(\xi;s)\leq\bar{v}_2(\xi)$ for all $\xi\geq M$ and $s\in\R$. This ends the proof of \eqref{exp-varphi-2}.
	
Finally, choosing $C_1=\max\{B_1,B_2\}$, and $\mu_1=\nu_2$, we obtain the first inequality of \eqref{uniform exp cover}. The second one can be argued in a similar way. The proof of 
Lemma \ref{uniform exp s} is thus complete. 
\end{proof}

As an easy consequence of Lemma \ref{uniform exp s}, we have $\partial_\xi \boldsymbol \varphi(\cdot;s)\in H^2(\R;\R^2)$ for each $s\in\R$, and 
\begin{equation}\label{L2norm-varphi}
\|1-\varphi_i(\cdot; s)\|_{L^2(-\infty, 0)}+\|\varphi_i(\cdot; s)\|_{L^2(0,+\infty)}+\left\|\partial_{\xi} \varphi_i(\cdot; s)\right\|_{H^2(\mathbb{R})} \leq C\,\,\hbox{ for }\,i=1,2,
\end{equation}
where $C$ is a positive constant independent of $s \in \mathbb{R}$.
As in \cite{dhl}, we will apply the implicit function theorem to show the $C^1$-smoothness of  $(\boldsymbol \varphi(\cdot; s), c(s))$ with respect to $s \in \mathbb{R}$. Let us first set some notations.  
In the sequel, we fix a real number $\beta>0$. For any $c \in \mathbb{R}$ and $s \in \mathbb{R}$, we define
$$
M_{c, s}(\mathbf{v})=D(s) \mathbf{v}^{\prime \prime}+c \mathbf{v}^{\prime}-\beta \mathbf{v} \quad \text {for} \quad \mathbf{v} \in H^2\left(\mathbb{R}; \mathbb{R}^2\right),
$$
where $D(s)$ is a diagonal matrix defined by
$$
D(s)=\left(\begin{array}{cc}
	d_1(s) & 0 \\
	0 & d_2(s)
\end{array}\right).
$$
Clearly, each $M_{c,s}$ maps $H^2(\mathbb{R};\R^2)$ into $L^2(\mathbb{R};\R^2)$. In the following lemma, we collect some basic properties of this operator. The proof follows the same lines as those used in \cite[Lemma 5.2]{dhl}.
 
\begin{lem}\label{M_c,s^-1}
Fix $\beta>0$. The following statements hold true.
\begin{enumerate}
\item [\rm{(i)}] For any $c \in \mathbb{R}$ and $s \in \mathbb{R}$, the operator $M_{c, s}: H^2(\mathbb{R};\R^2) \rightarrow L^2(\mathbb{R};\R^2)$ is invertible. Furthermore, for every $A>0$, there exists a constant $C>0$ \textup{(}depending on $\beta$ and $A$, but independent of $s$\textup{)} such that for any $c \in[-A, A]$, $s \in \mathbb{R}$ and $\textbf{h} \in L^2(\mathbb{R};\R^2)$,
\begin{equation*}
\left\|M_{c, s}^{-1}(\textbf{h})\right\|_{H^2(\mathbb{R})} \leq C\|\textbf{h}\|_{L^2(\mathbb{R})}.
\end{equation*}

\item [\rm{(ii)}]  
For any $\textbf{h} \in L^2(\mathbb{R};\R^2)$, there holds
\begin{equation*}
M_{c_n, s_n}^{-1}\left(\textbf{h}_n\right) \rightarrow M_{c, s}^{-1}(\textbf{h}) \quad \text {in } \,\,H^2(\mathbb{R};\R^2),
\end{equation*}
as $n \rightarrow+\infty$ for all sequences $(\textbf{h}_n)_{n \in \mathbb{N}} \subset L^2(\mathbb{R};\R^2)$, $(c_n)_{n \in \mathbb{N}} \subset \mathbb{R}$, $(s_n)_{n \in \mathbb{R}} \subset \mathbb{R}$ such that $\|\textbf{h}_n-\textbf{h} \|_{L^2(\mathbb{R})} \rightarrow 0$, $c_n \rightarrow c$ and $s_n \rightarrow s$ as $n \rightarrow+\infty$. Furthermore, the above convergence is uniform in $(\textbf{h}, c, s) \in B_A \times \mathbb{R}$ for any $A>0$, where $B_A$ is the ball given by $B_A=\{(\textbf{h}, c) \in L^2(\mathbb{R};\R^2) \times \mathbb{R}:\|\textbf{h}\|_{L^2(\mathbb{R})}+|c| \leq A \}$.
\end{enumerate}
\end{lem}

In addition to $\beta>0$, we also let $s_0\in\R$ be arbitrarily fixed. For any $(\textbf{v}, c, s) \in L^2(\mathbb{R};\R^2) \times \mathbb{R} \times \mathbb{R}$, we define
$$
\textbf{K}(\textbf{v}, c, s): \xi \mapsto K(\textbf{v}, c, s)(\xi):=D(s) \partial_{\xi \xi} \boldsymbol\varphi\left(\xi; s_0\right)+c \partial_{\xi} \boldsymbol\varphi\left(\xi; s_0\right)+\beta \textbf{v}(\xi)+\textbf{G}\left(s, \textbf{v}(\xi)+\boldsymbol\varphi\left(\xi; s_0\right)\right),
$$
where $\boldsymbol\varphi(\xi;s_0)=(\varphi_1,\varphi_2)(\xi; s_0)$ is the unique solution of \eqref{eq: frozen-limit-change} with $s=s_0$, and $\mathbf{G}:\mathbb{R} \times L^2(\mathbb{R};\R^2) \rightarrow L^2(\mathbb{R};\R^2)$ is a vector function defined by 
$$\mathbf{G}(s,\mathbf{v})=\left(g_1(s, \mathbf{v}), g_2(s, \mathbf{v})\right) $$
(recall that $g_1$, $g_2$ are the functions given by \eqref{g-frozen}). 
It is easily checked that $\textbf{K}(\textbf{0}, c(s_0), s_0)(\xi)=(\textbf{0},0)$ for all $\xi \in \mathbb{R}$. By \eqref{L2norm-varphi}, we know that $\varphi_i\left(\cdot; s_0\right) \in L^2(0,+\infty)$, $1-\varphi_i(\cdot; s_0) \in L^2(-\infty, 0)$, and $\partial_{\xi} \varphi_i(\cdot;s_0) \in H^2(\mathbb{R})$ for each $i=1,2$. Moreover,  the function $\xi \mapsto \textbf{G}(s, \textbf{v}(\xi)+\boldsymbol\varphi(\xi;s_0))$ belongs to $L^2(\mathbb{R};\R^2)$. Therefore, for any $(\textbf{v}, c, s) \in L^2(\mathbb{R};\R^2) \times \mathbb{R} \times \mathbb{R}$, we have $\textbf{K}(\textbf{v}, c, s) \in L^2(\mathbb{R};\R^2)$.

Now, for any $(\textbf{v}, c, s) \in H^2(\mathbb{R};\R^2) \times \mathbb{R} \times \mathbb{R}$, we set $\textbf{F}(\textbf{v}, c, s)=(F_1, F_2)(\textbf{v}, c, s)$ with
$$F_1(\textbf{v}, c, s)=\textbf{v}+M_{c, s}^{-1}(\textbf{K}(\textbf{v}, c, s))\quad  \hbox{and}\quad
F_2(\textbf{v}, c, s)=v_1(0),
$$
where $v_1(0)$ the first component of $\textbf{v}(0)$.
It is easily seen from Lemma \ref{M_c,s^-1} (i) that the function $\textbf{F}$ maps $H^2(\mathbb{R};\R^2) \times \mathbb{R} \times \mathbb{R}$ into $H^2(\mathbb{R};\R^2) \times \mathbb{R}$. Remember that $(\boldsymbol\varphi(\cdot;s),c(s))$ is the unique solution of \eqref{eq: frozen-limit-change} satisfying \eqref{normalize fix s wave}. Due to \eqref{L2norm-varphi}, one readily checks that for each $s\in\R$,  $\boldsymbol\varphi(\cdot; s)-\boldsymbol\varphi\left(\cdot; s_0\right) \in H^2(\mathbb{R},\R^2)$, and
\begin{equation}\label{F0}
\textbf{F}\left(\boldsymbol\varphi(\cdot; s)-\boldsymbol\varphi\left(\cdot; s_0\right), c(s), s\right)=(\textbf{0},0).
\end{equation}
On the other hand, for each $s\in\R$, denote by $A(s,\mathbf{v})$ the Jacobian matrix of the function $\mathbf{G}(s, \mathbf{v})$, namely,
for any $\mathbf{v}=\left(v_1, v_2\right) \in \mathbb{R}^2$, 
$$
A(s, \mathbf{v})=\left(\begin{array}{cc}
	r_1(s)-2r_1(s)v_1-k_1(s)p_2(s) (1-v_2) & k_1(s)p_2(s)  v_1 \vspace{5pt}\\
	k_2(s)p_1(s)  (1-v_2) & -r_2(s)+2r_2(s)v_2-k_2(s)p_1(s)  v_1
\end{array}\right).
$$
Thanks to Lemma \ref{M_c,s^-1} (ii), one can check that the function $\textbf{F}: H^2(\mathbb{R};\R^2) \times \mathbb{R} \times \mathbb{R} \rightarrow H^2(\mathbb{R};\R^2) \times \mathbb{R}$ is continuously Fr\'echet differentiable, and that the derivatives are given by
$$
\partial_{(\textbf{v}, c, s)} \textbf{F}(\textbf{v}, c, s)(\tilde{\textbf{v}}, \tilde{c}, \tilde{s})=\left(\begin{array}{l}
\partial_\textbf{v} F_1(\textbf{v}, c, s)(\tilde{\textbf{v}})+\partial_c F_1(\textbf{v}, c, s)(\tilde{c})+\partial_s F_1(\textbf{v}, c, s)(\tilde{s})  \vspace{5pt}\\
\partial_\textbf{v} F_2(\textbf{v}, c, s)(\tilde{\textbf{v}})+\partial_c F_2(\textbf{v}, c, s)(\tilde{c})+\partial_s F_2(\textbf{v}, c, s)(\tilde{s})
\end{array} \right)$$
for $(\tilde{\textbf{v}}, \tilde{c}, \tilde{s}) \in H^2(\mathbb{R};\R^2) \times \mathbb{R} \times \mathbb{R}$, where
$$
\left\{\begin{array}{l}
	\partial_\textbf{v} F_1(\textbf{v}, c, s)(\tilde{\textbf{v}})=\tilde{\textbf{v}}+M_{c, s}^{-1}\left[\left( A\left(s, \textbf{v}+\boldsymbol\varphi\left(\cdot; s_0\right)\right) +\beta I\right)  \tilde{\textbf{v}}\right], \vspace{5pt}\\
	\partial_c F_1(\textbf{v}, c, s)(\tilde{c})=-\tilde{c} M_{c, s}^{-1}\left\{\partial_{\xi}\left[M_{c, s}^{-1}(\textbf{K}(\textbf{v}, c, s))-\boldsymbol\varphi\left(\cdot; s_0\right)\right]\right\}, \vspace{5pt}\\
	\partial_s F_1(\textbf{v}, c, s)(\tilde{s})=-\tilde{s} M_{c, s}^{-1}\left\{D^{\prime}(s) \partial_{\xi \xi}\left[M_{c, s}^{-1}(\textbf{K}(\textbf{v}, c, s))-\boldsymbol\varphi\left(\cdot; s_0\right)\right]-\partial_s \textbf{G}\left(s, \textbf{v}+\boldsymbol\varphi\left(\cdot; s_0\right)\right)\right\},
\end{array}\right.
$$
and
$$
\partial_\textbf{v} F_2(\textbf{v}, c, s)(\tilde{\textbf{v}})=\tilde{v}_1(0), \quad \partial_c F_2(\textbf{v}, c, s)(\tilde{c})=0, \quad \partial_s F_2(\textbf{v}, c, s)(\tilde{s})=0.
$$
Here, $I$ is the identity matrix, and $\tilde{v}_1(0)$ is the first component of $\tilde{\textbf{v}}(0)$.
Furthermore, at the point $(\textbf{v}, c, s)=(\textbf{0}, c(s_0), s_0)$, the operator 
$$Q_{s_0}=\partial_{(\textbf{v}, c)} \textbf{F}(\textbf{0}, c(s_0), s_0): H^2(\mathbb{R};\R^2) \times \mathbb{R} \rightarrow H^2(\mathbb{R};\R^2) \times \mathbb{R}$$ is invertible, and we have the following property.

\begin{lem}\label{Q^-1}
There exists a constant $C>0$ independent of $s_0$ such that for any $(\tilde{\textbf{h}}, \tilde{d}) \in H^2(\mathbb{R};\R^2) \times \mathbb{R}$, 
\begin{equation*}
\left\|Q_{s_0}^{-1}(\tilde{\textbf{h}},  \tilde{d})\right\|_{H^2(\mathbb{R}) \times \mathbb{R}} \leq C \|(\tilde{\textbf{h}}, \tilde{d}) \|_{H^2(\mathbb{R}) \times \mathbb{R}},
\end{equation*}
where the space $H^2(\R;\R^2) \times \mathbb{R}$ is endowed with norm $\|(\tilde{\textbf{h}}, \tilde{d})\|_{H^2(\mathbb{R}) \times \mathbb{R}}=\|\tilde{\textbf{h}}\|_{H^2(\mathbb{R})}+|\tilde{d}|$.
\end{lem}
The proof of Lemma \ref{Q^-1} follows from that of \cite[Lemma 5.3]{dhl} with some obvious modifications; therefore we omit the details.

Based on the above preparations, we are now ready to complete the proof of Proposition \ref{fix s continue} (iii). 
	For any $s_0 \in \mathbb{R}$, we define
	$$
	P(s)=\left(\boldsymbol\varphi\left(\cdot; s\right)-\boldsymbol\varphi\left(\cdot; s_0\right), c(s)\right)\,\, \text { for }\,\, s \in \mathbb{R}.
	$$
	Clearly, $P$ maps $\mathbb{R}$ into $H^2(\mathbb{R};\R^2) \times \mathbb{R}$, and by \eqref{F0}, $\textbf{F}(P(s), s)=(\textbf{0},0)$ for all $s \in \mathbb{R}$. Applying the implicit function theorem to the function $\textbf{F}: H^2(\mathbb{R};\R^2) \times\mathbb{R} \times \mathbb{R} \rightarrow H^2(\mathbb{R};\R^2) \times \mathbb{R}$, one infers that, for each $s_0 \in \mathbb{R}$, there is $\delta>0$ such that the function $P:\left(s_0-\delta, s_0+\delta\right) \rightarrow H^2(\mathbb{R};\R^2) \times \mathbb{R}$ is continuously Fr\'echet differentiable. Denote the derivative operator at $s \in\left(s_0-\delta, s_0+\delta\right)$ by $O_s: \mathbb{R} \rightarrow H^2(\mathbb{R};\R^2) \times \mathbb{R}$. Thus, the function $s \mapsto O_s$ is continuous from $\left(s_0-\delta, s_0+\delta\right)$ to $\mathcal{L}\left(\mathbb{R}; H^2(\mathbb{R};\R^2) \times \mathbb{R}\right)$. At $s=s_0$, we have $O_{s_0}(\tilde{s})=-Q_{s_0}^{-1}(\partial_s \textbf{F}(\textbf{0}, c(s_0), s_0)(\tilde{s}),0)$ for every $\tilde{s} \in \mathbb{R}$. Furthermore, 
	 by Lemmas \ref{M_c,s^-1} and \ref{Q^-1}, there exists $C>0$, independent of $s_0 \in \mathbb{R}$, such that 
	 \begin{equation}\label{O_s_0 bound}
	 \left\|O_{s_0}(\tilde{s})\right\|_{H^2(\mathbb{R}) \times \mathbb{R}} \leq C|\tilde{s}|\,\,\hbox{ for all }\,\, \tilde{s} \in \mathbb{R}.
	 \end{equation} 
	 With a slight abuse of notation, we identify $O_{s_0}=(O_{s_0}^1,O_{s_0}^2) \in \mathcal{L}\left(\mathbb{R}; H^2(\mathbb{R};\R^2) \times \mathbb{R}\right)$ to an element of $H^2(\mathbb{R};\R^2) \times \mathbb{R}$. 
	 Then, one can check that the functions $s_0 \mapsto c(s_0)$ and $s_0 \mapsto \boldsymbol\varphi(\cdot; s_0)$ are of class $C^1$ in $\mathbb{R}$ with derivatives $\partial_s \boldsymbol\varphi(\cdot; s_0)=O_{s_0}^1$ and $c'(s_0)=O_{s_0}^2$ at each $s_0 \in \mathbb{R}$. It further follows from \eqref{O_s_0 bound} that $\|\partial_s \boldsymbol\varphi(\cdot;s_0)\|_{H^2(\mathbb{R})} \leq C$ for all $s_0 \in \mathbb{R}$.  This together with the Sobolev inequality implies \eqref{partial_s phi uniform bound}. The proof of Proposition \ref{fix s continue} (iii) is thus complete.

\section*{Acknowledgements} 
This work has received funding from NSFC (12471197), Guangdong Basic and Applied Basic Research Foundation (2023B1515020034) and Science and Technology Projects in Guangzhou (SL2024A04J00172).

\section*{Declarations}

{\bf Conflict of interest}:  There is no conflict of interest.

\vskip 0.2cm 

\noindent{\bf Data availability}: Data sharing is not applicable to this article as no datasets were generated or analyzed during the current study.


\begin{thebibliography}{AAA}	

%\bibitem{abc} N. D. Alikakos, P. W. Bates, X. Chen, {\it Periodic traveling waves and locating oscillating patterns in multidimensional domains}, Trans. Amer. Math. Soc., 351 (1999), 2777-2805.

\bibitem{aw} D. G. Aronson, H. F. Weinberger, {\it Multidimensional nonlinear diffusion arising in population genetics}, Adv. Math.,  30 (1978), 33-76.


\bibitem{B} X. Bao, Z.-C. Wang, {\it Existence and stability of time periodic traveling waves for a periodic bistable Lotka-Volterra competition system}, J. Differ. Equ., 255 (2013), 2402-2435.


%\bibitem{bls} X. Bao, W.-T. Li, W. Shen, {\it Traveling Wave Solutions of Lotka-Volterra Competition Systems with Nonlocal Dispersal in Periodic Habitats},  J. Differ. Equ., 260 (2016),  8590-8637.

\bibitem{blsw} X. Bao, W.-T. Li, W. Shen,  {\it Spreading speeds and linear determinacy of time dependent diffusive cooperative/competitive systems}, J. Differ. Equ., 265 (2018), 3048-3091.

\bibitem{ccs} M.-S. Chang, C.-C. Chen, S.-C. Wang, {\it Propagating direction in the two species Lotka-Volterra competition-diffusion system}, Discret. Contin. Dyn. Syst. B, 28 (2023), 59986014.

\bibitem{carrere} C. Carrère, {\it Spreading speeds for a two-species competition-diffusion system}, J. Differ. Equ., 264 (2018) 2133-2156.

\bibitem{C2} C. Conley, R. A. Gardner, {\it An application of the generalized Morse index to traveling wave solutions of a competitive reaction diffusion model}, Indiana Univ. Math. J., 33 (1984), 319-343.

%\bibitem{C3} B. Contri, {\it Pulsating fronts for bistable on average reaction-diffusion equations in a time periodic environment}, J. Math. Anal. Appl., 437 (2016), no.1, 90-132.

\bibitem{D1} W. Ding, {\it Average speeds of time almost periodic traveling waves for rapidly/slowly oscillating reaction-diffusion equations}, arXiv: 2406.06928, 2024.

\bibitem{dhl}
W. Ding,  F. Hamel, X. Liang, {\it Bistable traveling fronts in slowly oscillating one-dimensional environments}, J. Math. Pures Appl., 194 (2025), 103668. 

\bibitem{dhz} W. Ding, F. Hamel, X.-Q. Zhao, {\it Bistable traveling fronts for reaction-diffusion equations in a periodic habitat}, Indiana Univ. Math. J., 66 (2017), 1189-1265.
 
\bibitem{dhy} W. Ding, R. Huang, X. Yu, {\it Pulsating fronts of Lotka-Volterra competition system in rapidly varying media: Strong competition}, Disc. Cont. Dyn. Systems, 43 (2023), 2337-2370.

\bibitem{dl} W. Ding, X. Liang, {\it Sign of the traveling wave speed for the bistable competition-diffusion system in a periodic media}, Math. Ann., 385 (2023), 1001-1036.

\bibitem{dll} W. Ding, Z. Liang, W. Liu, {\it Continuity of pulsating wave speeds for bistable reaction-diffusion equations in spatially periodic media}, J. Math. Anal. Appl., 519 (2023), 126794.


\bibitem{dls}   L.-J. Du, W.-T. Li, W. Shen, {\it Propagation phenomena for time-space periodic monotone semiflows and applications to cooperative systems in multi-dimensional media}, J. Funct. Anal., 282 (2022), 109415.

\bibitem{dlw}  L.-J. Du, W.-T. Li, S.-L. Wu, {\it Propagation phenomena for a bistable Lotka-Volterra competition system with advection in a periodic habitat}, Z. Angew. Math. Phys., 67 (2020) 71.

\bibitem{ducrot} A. Ducrot, {\it A multi-dimensional bistable nonlinear diffusion equation in a periodic medium}, Math. Ann. 366 (2016), 783-818.

%\bibitem{fz1} S. Fan,  X.-Q. Zhao,  {\it Propagation dynamics of two species competition models in a periodic discrete habitat}, J. Differ. Equ., 377 (2023), 544-577.

\bibitem{fz2} S. Fan,  X.-Q. Zhao, {\it Bistable waves for two species competition systems in a periodic discrete habitat}, Disc. Cont. Dyn. Systems, 45 (2025), 4349-4387. 


\bibitem{fz} J. Fang, X.-Q. Zhao, {\it Bistable traveling waves for monotone semiflows with applications}, J. Eur. Math. Soc., 17 (2015), 2243-2288.

\bibitem{fyz} J. Fang, X. Yu, X.-Q. Zhao, {\it Traveling waves and spreading speeds for time-space periodic monotone 
systems}, J. Funct. Anal., 272 (2017), 4222-4262.

%\bibitem{F2} P. C. Fife, J. B. McLeod, {\it The approach of solutions of nonlinear diffusion equations to travelling front solutions}, Arch. Ration. Mech. Anal., 65 (1977), 335-361.

\bibitem{G1} R. A. Gardner, {\it Existence and stability of travelling wave solutions of competition models: a traveling wave theoretic approach}, J. Differ. Equ., 44 (1982), 343-364.

\bibitem{gl} J.-S. Guo, Y.-C. Lin, {\it The sign of the wave speed for the Lotka-Volterra competition diffusion system},  Commun. Pure Appl. Anal., 12 (2013), 2083-2090.

\bibitem{gno} J.-S. Guo, K.-I. Nakamura, T. Ogiwara, C.-H. Wu,  {\it The sign of traveling wave speed in bistable dynamics}, Discret. Contin. Dyn. Syst., 40 (2020), 3451-3466.

\bibitem{gw12} J.-S. Guo, C.-H. Wu, {\it  Recent developments on wave propagation in $2$-species competition systems}, Discrete Contin. Dyn. Syst. Ser. B, 17 (2012), 2713-2724.

\bibitem{girardin}
L. Girardin,  {\it Competition in periodic media: I-Existence of pulsating fronts}, Discret. Contin. Dyn. Syst. B, 22 (2017), 1341-1360.

\bibitem{girardin19} L. Girardin, {\it The effect of random dispersal on competitive exclusion-A review}, Math. Biosci., 318 (2019), 108271.

\bibitem{gl19} L. Girardin, K.-Y. Lam, {\it Invasion of an empty habitat by two competitors: spreading properties of monostable two-species competition-diffusion systems}, Proc. Lond. Math. Soc., 119 (2019), 1279-1335.

\bibitem{gn}
L. Girardin,  G. Nadin, {\it Competition in periodic media: II-Segregative limit of traveling fronts and ``Unity is not Strength''-type result},  J. Differ. Equ.,  265 (2018), 98-156.

%\bibitem{Heinze} S. Heinze, {\it Wave solutions to reaction-diffusion systems in perforated domains}, Z. Anal. Anwendungen, 20 (2001), 661-670.

\bibitem{hps} S. Heinze, G. Papanicolaou, A. Stevens, {\it Variational principles for propagation speeds in inhomogeneous media}, SIAM J.Appl. Math., 62 (2001), 129-148.

\bibitem{hess} P. Hess, {\it Periodic-Parabolic Boundary Value Problems and Positivity}, Pitman Res. Notes Math. Ser., vol. 247, Longman Scientific \& Technical, Wiley, Harlow, Essex, 1991.

\bibitem{hosono} Y. Hosono, {\it The minimal speed of traveling fronts for a diffusive Lotka-Volterra competition model}, Bull. Math. Biol. 60 (1998), 435-448.

\bibitem{K1} Y.Kan-on, {\it Parameter dependence of propagation speed of travelling waves for competition-diffusion equations}, SIAM J. Math. Anal., 26 (1995), 340-363.

\bibitem{Kan97} Y. Kan-on, {\it Fisher wave fronts for the Lotka-Volterra competition model with diffusion}, Nonlinear Anal. 28 (1997) 145-164.	

\bibitem{K2} Y. Kan-on, Q. Fang, {\it Stability of monotone travelling waves for competition-diffusion equations}, Japan J. Indust. Appl. Math., 13 (1996), 343-349.


\bibitem{llw} M. A. Lewis, B. Li, H. F. Weinberger, {\it Spreading speed and linear determinacy for two-species competition models}, J. Math. Biol., 45 (2002), 219-233.

\bibitem{lz07} X. Liang, X.-Q. Zhao, {\it Asymptotic speeds of spread and traveling waves for monotone semiflows with applications}, Comm. Pure Appl. Math., 60 (2007), 1-40.

\bibitem{lz} Y. Li, X.-Q. Zhao, {\it The propagation dynamics for three species competitive-cooperative reaction-diffusion systems}, Calc. Var. Partial Differential Equations, 64 (2025), 28. 

%\bibitem{lw} G. Lin, W.-T. Li, Asymptotic spreading of competition diffusion systems: the role of interspecific competitions, Eur. J. Appl. Math. 23 (2012), 669-689.

\bibitem{mho} M. Ma, Z. Huang, C. Ou, {\it  Speed of the traveling wave for the bistable Lotka-Volterra competition model}, Nonlinearity, 32 (2019), 3143-3162.

\bibitem{myh} M. Ma, J. Yue, Z. Huang, C. Ou, {\it Propagation dynamics of bistable traveling wave to a time-periodic Lotka-Volterra competition model: Effect of seasonality},  J. Dyn. Diff. Equat., 35 (2023), 1745-1767.


\bibitem{M4} M. Iida, T. Muramatsu, H. Ninomiya, E. Yanagida, {\it  Diffusion-induced extinction of a superior species in a competition system}, Japan J. Indust. Appl. Math., 15 (1998), 233.


\bibitem{pwr} L. Pang, S.-L. Wu, S. Ruan, {\it Long time behavior for a periodic Lotka-Volterra reaction-diffusion system with strong competition}, Calc. Var. Partial Differential Equations, 62 (2023), 99.

\bibitem{pwz} R. Peng, C.-H. Wu,  M. Zhou,  {\it Sharp estimates for the spreading speeds of the Lotka-Volterra diffusion system with strong competition}, Ann. Inst. H. Poincaré C Anal. Non Linéaire, 38 (2021), 507-547. 

\bibitem{M3} M. Rodrigo, M. Mimura, {\it Exact solutions of reaction-diffusion systems and nonlinear wave equations}, Japan J. Indust. Appl. Math., 18 (2001), 657-696.

%\bibitem{S} 
%H. L. Smith, {\it Monotone Dynamical Systems: An Introduction to the Theory of Competitive
%and Cooperative Systems}, Mathematical Surveys and Monographs, vol. 41, American Mathematical Society, Providence, RI, 1995.

\bibitem{tw} J.-C. Tsai, Y.-Y. Weng,  {\it Propagation direction of traveling waves for a class of bistable epidemic models},  J. Math. Biol., 81 (2020), 1465-1493. 

\bibitem{wo} H. Wang, C. Ou, {\it Propagation speed of the bistable traveling wave to the Lotka-Volterra competition system in a periodic habitat}, J. Nonlinear Sci., 30 (2020), 3129-3159.

\bibitem{wxz23}  C.-H. Wu, D. Xiao and M. Zhou, {\it Sharp estimates for the spreading speeds of the Lotka-Volterra
competition-diffusion system: the strong-weak type},  J. Math. Pures Appl., 172 (2023), 236-264. 

\bibitem{xiao} D. Xiao, {\it Sufficient conditions for determining the sign of the wave speed in the Lotka-Volterra competition system}, arXiv:2408.10481, 2024.

%\bibitem{xz24}  D. Xiao and M. Zhou, {\it Complete classification on traveling waves in monotone dynamical systems}, arXiv:2409.12463.

\bibitem{yz} X. Yu, X.-Q. Zhao, {\it Propagation phenomena for a reaction-advection-diffusion competition model in a periodic habitat}, J. Dyn. Diff. Equat., 29 (2017), 41-66.

\bibitem{ymo} J. Yue, M. Ma, C. Ou, {\it Traveling wave for a time-periodic Lotka-Volterra model with bistable nonlinearity}, Appl. Math. J. Chinese Univ. Ser. B, 37 (2022), 396-403.

\bibitem{zr} G. Zhao, S. Ruan, {\it Existence, uniqueness and asymptotic stability of time periodic traveling waves for a periodic Lotka-Volterra competition system with diffusion}, J. Math. Pures Appl., 95 (2011) 627–671.

\bibitem{zr2} G. Zhao, S. Ruan, {\it Time periodic traveling wave solutions for periodic advection-reaction-diffusion systems}, J. Differ. Equ., 257 (2014), 1078-1147.

\bibitem{zhaobook} X.-Q. Zhao, {\it  Dynamical Systems in Population Biolog}y, 2nd edn. Springer, New York, (2017).
\end{thebibliography}
\end{document}